\documentclass[11pt]{article}

\usepackage{
a4dutch,
xspace, calc,
amsthm, 
amssymb,amsbsy,amsfonts, 
amscd, diagrams, 
stmaryrd, 
curves,epsfig,wrapfig, 
index} 
\usepackage[centertags]{amsmath}
\usepackage[refpage]{nomencl} 

\author{Benjamin Friedrich}

\makeglossary
\makeindex

\newcounter{notnr}[section]
\newcommand{\notation}[3]{
  \stepcounter{notnr}\nomenclature[let:#1, sec:\thesection, nr:\alph{notnr} ]{#2}{#3}}
\newcommand{\ignore}[1]{#1} 
\newcommand{\myunderbrace}[2]{\underset{#1}{\underbrace{#2}}}
\newcommand{\myatop}[2]{\genfrac{}{}{0mm}{1}{#1}{#2}}
\newcommand{\ts}{\textstyle}
\newcommand{\fk}{\mathfrak}

\DeclareMathOperator{\Spec}{Spec}

\DeclareMathOperator{\Iso}{Iso}
\DeclareMathOperator{\Aut}{Aut}
\DeclareMathOperator{\im}{im}
\DeclareMathOperator{\rank}{rank}
\DeclareMathOperator{\coker}{coker}
\DeclareMathOperator{\Hom}{Hom}
\DeclareMathOperator{\Supp}{Supp}
\DeclareMathOperator{\tot}{tot}
\DeclareMathOperator{\id}{id}
\DeclareMathOperator{\linspan}{span}
\newcommand{\Li}{{\rm Li}}

\newtheorem{definition}{Definition}[subsection]
\newtheorem{lemma}[definition]{Lemma}
\newtheorem{proposition}[definition]{Proposition}
\newtheorem{theorem}[definition]{Theorem}
\newtheorem{maintheorem}[definition]{Main Theorem}
\newtheorem{corollary}[definition]{Corollary}
\newtheorem{fact}[definition]{Fact}
\newtheorem{claim}[definition]{Claim}
\newtheorem{conjecture}[definition]{Conjecture}
\newtheorem{preremark}[definition]{Remark}
\newtheorem{prenote}[definition]{Note}
\newtheorem{preexample}[definition]{Example}
\newenvironment{remark}{\begin{preremark}\upshape}{\end{preremark}}
\newenvironment{note}{\begin{prenote}\upshape}{\end{prenote}}
\newenvironment{example}{\begin{preexample}\upshape}{\end{preexample}}

\begin{document}

\setlength{\parindent}{0.0mm}
\advance\voffset by -1.00cm
\advance\textheight by  2\baselineskip
\advance\hsize by 1.00cm
\advance\textwidth by 1.00cm 

\newenvironment{enumroman}
{
\renewcommand{\labelenumi}{(\roman{enumi})}
\begin{enumerate}
}
{
\end{enumerate}
\renewcommand{\labelenumi}{(\arabic{enumi})}
}

\renewcommand{\a}{{\rm A}}
\newcommand{\A}{\ensuremath{\mathbb{A}}}
\renewcommand{\AA}{\ensuremath{\mathcal{A}}}
\newcommand{\AAA}{\ensuremath{\mathfrak{A}}}
\renewcommand{\b}{{\rm B}}
\newcommand{\B}{\ensuremath{\mathbb{B}}}
\newcommand{\BB}{\ensuremath{\mathcal{B}}}
\renewcommand{\c}{{\rm C}}
\newcommand{\C}{\ensuremath{\mathbb{C}}}
\newcommand{\CC}{\ensuremath{\mathcal{C}}}
\newcommand{\DDD}{\ensuremath{\mathfrak{D}}}
\newcommand{\E}{\ensuremath{\mathbb{E}}}
\newcommand{\EE}{\mathcal{E}}
\newcommand{\F}{\ensuremath{\mathbb{F}}}
\newcommand{\FF}{\ensuremath{\mathcal{F}}}
\newcommand{\f}{\ensuremath{\rm F}}
\newcommand{\G}{\ensuremath{\mathbb{G}}}
\newcommand{\GG}{\ensuremath{\mathcal{G}}}
\newcommand{\g}{\ensuremath{\rm G}}
\newcommand{\h}{\ensuremath{\text{\rm H}}}
\newcommand{\hh}{\ensuremath{\text{\rm h}}}
\renewcommand{\H}{\mathbb{H}}
\newcommand{\HH}{\mathcal{H}}
\newcommand{\HHH}{\boldsymbol{\mathcal H}}
\newcommand{\I}{\mathbb{I}}
\newcommand{\II}{\mathcal{I}}
\newcommand{\ii}{{\rm I}}
\newcommand{\K}{\ensuremath{\mathbb{K}}}
\newcommand{\KK}{\ensuremath{\mathcal{K}}}
\newcommand{\LL}{\ensuremath{\mathcal{L}}}
\newcommand{\m}{{\rm M}}
\newcommand{\N}{\ensuremath{\mathbb{N}}}
\newcommand{\NN}{\ensuremath{\mathcal{N}}}
\newcommand{\OO}{\ensuremath{\mathcal O}}
\renewcommand{\P}{\ensuremath{\mathbb{P}}}
\newcommand{\PP}{\ensuremath{\mathcal{P}}}
\newcommand{\Q}{\ensuremath{\mathbb Q}}
\newcommand{\QQ}{\ensuremath{\overline{\Q}}}
\newcommand{\QR}{\ensuremath{\widetilde{\Q}}}
\renewcommand{\r}{\ensuremath{\text{\rm R}}}
\newcommand{\R}{\ensuremath{\mathbb{R}}}
\newcommand{\RR}{\ensuremath{\mathcal{R}}}
\newcommand{\T}{T}
\newcommand{\U}{\ensuremath{\mathbb{U}}}
\newcommand{\UU}{\ensuremath{\mathcal{U}}}
\newcommand{\UUU}{\ensuremath{\mathfrak{U}}}
\newcommand{\XX}{\ensuremath{\mathcal{X}}}
\newcommand{\Z}{\ensuremath{\mathbb{Z}}}
\newcommand{\ZZ}{\ensuremath{\mathcal{Z}}}
\newcommand{\eps}{\epsilon}

\newcommand{\Triangle}{\ensuremath{\triangle^{\rm std}}}
\renewcommand{\epsilon}{\varepsilon}
\newcommand{\iso}{\cong}
\newcommand{\down}[1][\wr]{{\Big\downarrow\vcenter{%
  \rlap{$\scriptstyle #1 $}}}}
\newcommand{\up}[1][\wr]{{\Big\uparrow\vcenter{%
  \rlap{$ #1 $}}}}
\newcommand{\downiso}{\downarrow \! 
  \text{\raisebox{0.7mm}{$\scriptstyle \wr$}}}
\newcommand{\upiso}{\uparrow \! 
  \text{\raisebox{-0.1mm}{$\scriptstyle \wr$}}}
\renewcommand{\to}{\rightarrow}
\newcommand{\longto}{\longrightarrow}
\newcommand{\inclusion}{\hookrightarrow}
\newcommand{\injection}{\rightarrowtail}
\newcommand{\surjection}{\twoheadrightarrow}
\newcommand{\tensor}[1]{\otimes_{#1}}
\newcommand{\Union}{\bigcup}
\newcommand{\union}{\cup}
\newcommand{\del}{\partial}
\newcommand{\delbar}{\overline{\del}}
\renewcommand{\emptyset}{\varnothing}
\newcommand{\defeq}{\sim}
\newcommand{\st}{\,\vert \,} 
\newcommand{\Bul}{\bullet}
\newcommand{\bul}{^{\Bul}}
\renewcommand{\injlim}{\varinjlim}
\renewcommand{\projlim}{\varprojlim}
\newcommand{\circled}[1]{\ensuremath{\bigcirc\hspace{-7pt}\scriptstyle#1}}

\newcommand{\hDR}[1]{\h^{#1}_{\rm dR}}
\newcommand{\Xo}{\ensuremath{X_0}} \newcommand{\Do}{\ensuremath{D_0}}
\newcommand{\ko}{\ensuremath{k_0}}
\newcommand{\XC}{\ensuremath{X(\C)}}
\newcommand{\Xh}{\ensuremath{X^{\rm an}}}
\newcommand{\Yh}{\ensuremath{Y^{\rm an}}}
\newcommand{\Zh}{\ensuremath{Z^{\rm an}}}
\newcommand{\Uh}{\ensuremath{U^{\rm an}}}
\newcommand{\Vh}{\ensuremath{V^{\rm an}}}
\newcommand{\Dh}{\ensuremath{D^{\rm an}}} 
\newcommand{\Eh}{\ensuremath{E^{\rm an}}}
\newcommand{\Bh}{\ensuremath{B^{\rm an}}}
\newcommand{\Ph}{\ensuremath{P_{\rm an}}} 
\newcommand{\DDDh}{\DDD^{\rm}} 
\newcommand{\FFh}{\ensuremath{\FF_{\rm an}}} \renewcommand{\Re}{{\rm Re\,}} 
\renewcommand{\Im}{{\rm Im\,}} 
\newcommand{\andsym}{\,\wedge\,}
\newcommand{\zbar}{\overline{z}} \newcommand{\om}[1]{\Omega_{#1}\bul}
\newcommand{\wom}[1]{\widetilde{\Omega}_{#1}\bul}
\newcommand{\iast}{i_{\ast}} \newcommand{\sing}{_{\rm sing}}
\newcommand{\ssing}{\rm smooth} \newcommand{\Ch}{\check{\rm C}}
\newcommand{\CCh}{\check{\CC}}
\newcommand{\deRham}{de\hspace{0.3mm}Rham\xspace}
\newcommand{\DeRham}{De\hspace{0.3mm}Rham\xspace}
\newcommand{\osimplex}{\triangle\!\raisebox{0.5mm}{$^{\circ}$}\!}
\newcommand{\varosimplex}{\triangle\!\!\raisebox{0.15mm}{$^{\circ}$}\!}
\newcommand{\csimplex}{\triangle}

\newcommand{\ul}[1]{\ensuremath{\underline{#1}}}
\newcommand{\ol}[1]{\ensuremath{\overline{#1}}}
\newcommand{\wt}[1]{\ensuremath{\widetilde{#1}}}

\begin{titlepage}

\vspace{3cm}

\begin{center}
  \Large
  Universit{\"a}t Leipzig\\
  Fakult{\"a}t f{\"u}r Mathematik und Informatik \\
  Mathematisches Institut
\end{center}

\vspace{1cm}

\begin{center}
  \huge Periods and Algebraic \deRham Cohomology
\end{center}

\vspace{2cm}

\begin{center}
\Large
Diplomarbeit \\
im Studiengang Diplom-Mathematik
\end{center}

\vspace{12cm}

\large
\begin{tabular}{ll}
Leipzig, \hspace*{4cm} & vorgelegt von Benjamin Friedrich,\\
& geboren am 22. Juli 1979
\end{tabular}
\end{titlepage}

\begin{titlepage} 
\vspace*{3cm}

\section*{Abstract}

It is known that the algebraic \deRham cohomology group $\hDR{i}(X_0/\Q)$ of 
a nonsingular variety $X_0/\Q$ has the same rank as the rational
singular cohomology group $\h^i\sing(\Xh;\Q)$ of the complex
manifold $\Xh$ associated to the base change $X_0\times_{\Q}\C$. 
However, we do not have a natural isomorphism
$\hDR{i}(X_0/\Q)\iso\h^i\sing(\Xh;\Q)$. Any choice of
such an isomorphism 
produces certain integrals, so called periods,
which reveal valuable information about $X_0$. The aim of this thesis
is to explain these classical facts in detail. Based on an approach of
Kontsevich \cite[pp.~62--64]{kontsevich}, different definitions of a period are
compared and their properties discussed. Finally, the
theory is applied to some examples. These examples include a
representation of $\zeta(2)$ as a period 
and a variation of mixed Hodge structures used by Goncharov 
\cite{goncharov:dlog}.

\end{titlepage}

\setcounter{tocdepth}{2}
\begin{titlepage}
\enlargethispage*{3\baselineskip} 
\tableofcontents
\enlargethispage*{-3\baselineskip}
\end{titlepage}

\section{Introduction}

The prehistory of Algebraic Topology dates back to 
Euler, Riemann and Betti, 
who started the idea of attaching various invariants to a
topological space. 
With his simplicial (co)homology theory, 
Poincar{\'e} was the first to give an instance of
what in modern terms we would call a contravariant functor $\h\bul$
from the category of (sufficiently nice) topological spaces to
the category of cyclic complexes of abelian groups. 

Many of such functors have been found so far; the most common examples are the
standard cohomology theories (i.e. those satisfying the
Eilenberg-Steenrod axioms), which measure quite different phenomena
relating to diverse branches of mathematics. It is a beautiful basic
fact that all these standard cohomology theories agree (when
restricted to an appropriate subcategory).

This does not imply that we cannot hope for more. If the topological
space in question enjoys additional structure, one defines more
elaborate invariants which take values in an abelian category of
higher complexity. 
For example, Hodge theory gives us a functor from
the category of compact K{\"a}hler manifolds to the category of cyclic
complexes of pure Hodge structures.

In this thesis, we will concentrate on spaces originating from
Algebraic Geometry; these may be regarded as spaces carrying an
algebraic structure. 

Generalizing the concept of \deRham
theory to ``nice'' schemes over \Q\ gives us algebraic \deRham
cohomology groups where each is nothing but a full \Q-lattice inside the
\C-vector space of the corresponding classical \deRham cohomology
group. 
So, after tensoring with \C, algebraic \deRham 
cohomology agrees with all the
standard cohomology theories with complex coefficients.
However, a natural isomorphism between the original \Q-vector space
and a standard cohomology group with rational coefficients cannot
exist.

We will illustrate this phenomenon in the following example (see
Example \ref{noiso} for details). 
Let $\Xh:=\C^{\times}$ be the complex plane
with the point $0$ deleted and let $t$ be the standard
coordinate on $\Xh$. Then the first singular cohomology group of $\Xh$
is generated by the dual $\sigma^{\ast}$ of the unit circle
$\sigma:=S^1$
$$
\h^1\sing(\Xh;\Q) = \Q \, \sigma^{\ast}
\quad\text{and}\quad
\h^1\sing(\Xh;\C) = \C \, \sigma^{\ast};
$$
while for the first classical \deRham cohomology group, we have
$$
\hDR{1}(\Xh;\C)=\C\,\frac{dt}{t} .
$$
Under the comparison isomorphism
$$
\h^1\sing(\Xh;\C) \iso \hDR{1}(\Xh;\C)
$$
the generator $\sigma^{\ast}$ of $\h^1\sing(\Xh;\C)$ is mapped to
$$
\left( \int_{S^1}\frac{dt}{t}\right)^{-1} \frac{dt}{t} =
\frac{1}{2\pi i} \frac{dt}{t}.
$$
If we view $\Xh$ as the complex manifold associated to the
base change to \C\ of the algebraic variety $X_0:=\Spec \Q[t,t^{-1}]$
over \Q, we can also compute the algebraic \deRham cohomology group 
$\hDR{1}(X_0/\Q)$
of
$X_0$ and embed it into $\hDR{1}(\Xh;\C)$
$$
\hDR{1}(X_0/\Q) = \Q \frac{dt}{t} \subset \hDR{1}(\Xh;\C) = \C \frac{dt}{t}.
$$
Thus we get two \Q-lattices
inside $\hDR{1}(\Xh;\C)$, 
$\h^1\sing(\Xh;\Q)$ and $\hDR{1}(X_0/\Q)$,
which do not coincide. In fact, they differ
by the factor $2\pi i$ --- our first example of what we will call a
period. Other examples will produce period numbers like $\pi$, $\ln
2$, elliptic integrals, or $\zeta(2)$, which are interesting also from a number
theoretical point of view (cf. page \pageref{expper}).

There is some ambiguity about the precise definition of a period;
actually we will give four definitions in total:

\begin{enumroman}
\item
pairing periods 
  (cf. Definition \ref{per1} on page \pageref{per1})

\item
abstract periods 
  (cf. Definition \ref{per2} on page \pageref{per2})

\item
na{\"\i}ve periods 
  (cf. Definition \ref{per3} on page \pageref{per3})

\item
effective periods 
  (cf. Definition \ref{effperiodef} on page \pageref{effperiodef})

\end{enumroman}

For $X_0$ a nonsingular variety over \Q, we have a natural pairing
between the $i\,^{{\rm th}}$ algebraic \deRham cohomology of $X_0$ and the 
$i\,^{{\rm th}}$ singular homology group of the complex manifold $\Xh$
associated to the base change $X_0\times_{\Q}\C$
$$
\h^{\rm sing}_i(\Xh;\Q) \times \hDR{i}(X_0/\Q) \longto \C .
$$
The numbers which can appear in the image of this pairing (or its
version for relative cohomology) are called {\em pairing periods}; this is
the most traditional way to define a period.

In \cite[p.~62]{kontsevich}, Kontsevich gives the alternative
definition of {\em effective periods} which does not need algebraic \deRham cohomology and, at
least conjecturally, gives the set of all periods some extra algebraic
structure. We present his ideas in Subsection \ref{per4}.

{\em Abstract periods} describe just a variant of Kontsevich's definition. In fact, we
have a surjection from the set of effective periods to the set of
abstract ones (cf. page \pageref{evsurj}), 
which is conjectured to be an isomorphism.

{\em Na{\"\i}ve periods} are defined in an elementary way
and are used to
provide a connection between pairing periods and abstract periods.

In Kontsevich's paper \cite[p.~63]{kontsevich}, it is used that 
the notion of pairing and abstract period coincide.
The aim of this thesis is to show that
the following implications hold true (cf. Theorem \ref{perioperio})
$$
\text{abstract period} \Leftrightarrow 
\text{na{\"\i}ve period} \Rightarrow
\text{pairing period} .
$$

The thesis is organized as follows.
The discussion of the various definitions of a period makes up the
principal part of the work filling sections five to seven. 

{\bf Section two} gives an introduction to complex analytic
spaces. Additionally, we provide the connection to Algebraic Geometry
by defining the associated complex analytic space of a variety.

In {\bf Section three}, we define algebraic \deRham cohomology for
pairs consisting of a variety and a divisor on it. We also give some
working tools for this cohomology.

The aim of {\bf Section four} is to give a comparison theorem
(Theorem \ref{comp}) which states that algebraic \deRham 
and singular cohomology agree.

In {\bf Section five}, we present the definition of pairing, abstract,
and na{\"\i}ve periods and prove some of their properties.

{\bf Section six} provides some facts about the triangulation of
algebraic varieties.

{\bf Section seven} contains the main result (Theorem
\ref{perioperio}) about the implications between the various
definitions of a period mentioned above.
Furthermore, we give the definition of effective periods which
motivated the definition of abstract periods.

In the last section, {\bf Section eight}, we
consider five examples to give an application of the general theory.
Among them is a representation of $\zeta(2)$ as an abstract period 
and the famous double logarithm variation of mixed
Hodge structures used by Goncharov \cite{goncharov:dlog} 
whose geometric origin is emphasized.

\vspace*{3mm}
{\bf Conventions.} By a {\em variety}, we will always mean
\index{Variety} a reduced, quasi-projective scheme.
We will often deal with a variety $X_0$ defined over
some algebraic extension of \Q.
As a rule, skipping the subscript $0$ will 
always mean base change to \C
$$
X:= \text{``} X_0 \times_\Q \C \,\text{''} = 
X_0 \times_{\Spec \Q} \Spec \C.
$$
(An exception is section three, where arbitrary base fields are used.)
The complex analytic space associated to $X$ will be denoted by $\Xh$
(cf. Subsection \ref{defcomplexanspace}).

The sign conventions used throughout this thesis are listed in the appendix.

\vspace*{7mm}

\paragraph{Acknowledgments.}

I am greatly indebted to my supervisor Prof. A.~Huber-Klawitter
for her guidance and her encouragement. I very much appreciated the
informal style of our discussions in which she vividly pointed out to
me the central ideas of the mathematics involved.

I would also like to thank my fellow students R. Munck, M. Witte and
K. Zehmisch who read the manuscript and gave numerous
comments which helped to clarify the exposition.

 \section{The Associated Complex Analytic Space}

 \renewcommand{\d}{_{D^n}}

 Let $X$ be a variety over \C.
 The set $\vert X \vert$ 
 \notation{0}{$\vert X\vert$}{set of closed points of a variety $X$}
 of closed points of $X$ inherits the Zariski
 topology.
 However, we can also equip this set with the standard
 topology:
 for smooth $X$ this gives a complex manifold; in general we get a
 complex analytic space \Xh.

 The main reference for this section is \cite[B.1]{hartshorne}.

 \subsection{The Definition of the Associated Complex Analytic Space}
 \label{defcomplexanspace}

 We consider an example before giving the general definition of a
 complex analytic space.

 \begin{example} 
 \label{local}
   Let $D^n \subset \C^n$ be the polycylinder
   $$
   D^n:= \{ \underline{z} \in \C^n \st \vert z_i \vert < 1, i=1, \ldots, n \}
   $$
   and $\OO\d$ the sheaf of holomorphic functions on $D^n$. 
   For a set of holomorphic functions $f_1, \ldots, f_m \in
   \Gamma(D^n, \OO\d)$ we define
   \begin{equation}
   \begin{split}
   & \label{complxspacedef}
   \XX\d := \{ \ul{z} \in D^n \st f_1(\ul{z}) = \ldots = f_m(\ul{z})
   = 0 \} \\
   & \OO_{\XX\d} := \OO\d / (f_1, \ldots, f_m) .
   \end{split}
   \end{equation}
 \end{example}

 The locally ringed space $(\XX\d, \OO_{\XX\d})$ from this example is a complex
 analytic space. In general, complex analytic spaces are obtained by
 glueing spaces of the form (\ref{complxspacedef}).

 \begin{definition}[Complex analytic space, {\cite[B.1,
   p. 438]{hartshorne}}]
 \index{Complex analytic space|textbf}
   A locally ringed space $(\XX, \OO_{\XX})$ is called {\em complex
   analytic} if it is locally (as a locally ringed space) isomorphic to
   one of the form (\ref{complxspacedef}). A morphism of complex
   analytic spaces is a morphism of locally ringed spaces.
 \end{definition}

   For any scheme $(X,\OO_X)$ of finite type over \C\ we have an {\em
   associated complex analytic space} $(\Xh, \OO_{\Xh})$.
   \index{Complex analytic space!associated}

 \begin{definition}[Associated complex analytic space, {\cite[B.1,
   p. 439]{hartshorne}}]
 \label{asscplansp}
   Assume first that $X$ is affine. We fix an isomorphism
   $$
   X \iso \Spec \C [ x_1, \ldots , x_n] / (f_1, \ldots, f_m) 
   $$
   and then consider the $f_i$ as holomorphic functions on $\C^n$ in
   order to set
   \begin{align*}
   & \Xh := \{ \ul{z} \in \C^n \st f_1(\ul{z}) = \ldots =
   f_m(\ul{z})=0\} \\
   & \OO_{\Xh} := \OO_{\C^n} / (f_1, \ldots, f_m),
   \end{align*}
   where $\OO_{\C^n}$ denotes the sheaf of holomorphic functions on
   $\C^n$.

   For an arbitrary scheme $X$ of finite type over \C, we take a covering
   of $X$ by open affine subsets $U_i$. The scheme $X$ 
   is obtained by glueing the open sets $U_i$, so we can use the same
   glueing data to glue the complex analytic space $(U_i)^{\rm an}$
   into an analytic space $\Xh$. This is the associated complex
   analytic space of $X$.
 \end{definition}

 This construction is natural and we obtain a functor {\bf an} 
 \notation{a}{an}{functor yielding the associated complex analytic space}
 from the category of schemes of finite type over \C\ 
 to the category of complex analytic spaces.
 Note that its restriction  to the subcategory of smooth schemes maps
 into the category of complex manifolds as a consequence of the inverse
 function theorem (cf. \cite[Thm.~6, p.~20]{gunning}). 

 \begin{example}
 The complex analytic space associated to complex projective space
 is again complex projective space, but considered as a complex
 manifold.
 To avoid confusion in the subsequent sections, the notation $\C P^n$
 will be reserved for {\em complex projective space in the category of
 schemes}, whereas we write $\C \Ph^n$ for {\em complex projective space in
 the category of complex analytic spaces}.
 \notation{C}{$\C P^n$}{complex projective space in the category of
   schemes}
 \notation{C}{$\C \Ph^n$}{complex projective space in the category of
   complex analytic spaces}
 \end{example}

 For any scheme $X$ of finite type over \C, we have a natural map of
 locally ringed spaces
 \begin{equation}
 \label{phi}
 \phi : \Xh \rightarrow X
 \end{equation}
 which induces the identity on the set of closed points $\vert X \vert$
 of $X$.
 Note that $\phi^{\ast} \OO_X = \OO_{\Xh}$.

 \subsection{Algebraic and Analytic Coherent Sheaves}

 Let us consider sheaves of $\OO_X$-modules.
 The equality of functors
 $$
 \Gamma(X,?)=\Gamma(\Xh,\phi^{-1} ?)
 $$
 gives an equality of
 their right derived functors
 $$
 \h^i(X;?) = \r^i\Gamma(X;?)=\r^i\Gamma(\Xh;\phi^{-1}?) .
 $$
 Since $\phi^{-1}$ is an exact functor, the spectral sequence for the
 composition of the functors $\phi^{-1}$ and $\Gamma(\Xh;?)$
 degenerates and we obtain
 $$
 \h^i(\Xh;\phi^{-1}?)=\r^i\Gamma(\Xh;?)\circ\phi^{-1} =
 \r^i\Gamma(\Xh;\phi^{-1}?).
 $$
 Thus the natural map for $\FF$ a sheaf of $\OO_X$-modules
 $$
 \phi^{-1} \FF \rightarrow \phi^{\ast} \FF
 $$
 gives a natural map of cohomology groups
 \begin{equation}
 \label{gaganatmorph}
 \h^i(X;\FF) = \h^i(\Xh; \phi^{-1} \FF) \rightarrow 
 \h^i(\Xh;\phi^{\ast} \FF).
 \end{equation}

 Sheaf cohomology behaves particularly nice for coherent sheaves, this
 notion being defined as follows.
 \begin{definition}[Coherent sheaf]
 We define a {\em coherent sheaf} \index{Coherent sheaf} \FF\ on $X$
 (resp. \Xh) to be a sheaf of $\OO_X$-modules
 (resp. $\OO_{\Xh}$-modules) that is Zariski-locally (resp. locally in
 the standard topology) isomorphic to the cokernel of a morphism of
 free $\OO_X$-modules (resp. $\OO_{\Xh}$-modules) of finite rank
 \begin{equation}
 \label{coherdef}
 \OO_U^r \longrightarrow
 \OO_U^s \longrightarrow \FF_{|U} \longrightarrow 0,
 \quad\text{for}\quad
 U\subseteq X
 \text{ (resp. $U\subseteq \Xh$) }
 \quad\text{open},
 \end{equation}
 where $r,s \in \N$.
 \end{definition}

 For sheaves on $X$ this agrees with the definition given in
 \cite[II.5, p. 111]{hartshorne}. For this alternate definition, we
 need some notation: 
 If $U=\Spec A$ is an affine variety and M an $A$-module,
 we denote by $\wt{M}$ 
 \notation{0}{$\wt{M}$}{sheaf associated to a module $M$}%
 the sheaf on $U$ associated to $M$ 
 (i.e. the sheaf associated to the presheaf
 $V \mapsto \Gamma(V;\OO_{V}) \tensor{A} M \,\,
 for\,\, V \subseteq U \text{ open}$, 
 see \cite[II.5, p. 110]{hartshorne}).

 \begin{lemma}[cf. {\cite[II.5 Exercise 5.4, p. 124]{hartshorne}}]
 A sheaf \FF\ of $\OO_X$-modules is coherent if and  only if $X$ can be
 covered by open affine subsets $U_i = \Spec A_i$ such that
 $F_{|U_i} \iso \wt{M}_i$ for some finitely generated $A_i$-modules
 $M_i$.
 \end{lemma}

 \begin{proof}
 {\em ``if'':} 
 The $A_i$'s are Noetherian rings. Therefore any finitely
 generated $A_i$-module $M_i$ will be finitely presented
 $$
 A_i^r \longrightarrow A_i^s \longrightarrow M_i \longrightarrow 0 .
 $$
 Since localization is an exact functor, we get
 $$
 \OO_{U_i}^r \longrightarrow \OO_{U_i}^s \longrightarrow \wt{M}_i
 \longrightarrow 0 ,
 $$ 
 which proves the ``if''-part.
 \hfill $\blacksquare$

 {\em ``only if'':} 
 W.l.o.g. we may assume that the open subsets
 $U\subseteq X$ in (\ref{coherdef}) are affine $U=\Spec A$. 
 Then $\Gamma(U;\OO_U)=A$ and $\OO_U=\wt{A}$.
 Now the $A$-module 
 \begin{align*}
 M:= & \coker \bigl(\Gamma(U;\OO_U^r) \longrightarrow \Gamma(U;\OO_U^s)\bigr)\\
 = & \coker(A^r \longrightarrow A^s)
 \end{align*}
 is clearly finitely generated. Since
 $$
 A^r \longrightarrow A^s \longrightarrow M \longrightarrow 0
 $$
 gives
 $$
 \wt{A}^r \longrightarrow \wt{A}^s \longrightarrow \wt{M}
 \longrightarrow 0 ,
 $$
 we conclude
 $$
 \FF_{|U} = \coker(\OO_U^r\longrightarrow\OO_U^s) 
   = \coker(\wt{A}^r\longrightarrow \wt{A}^s)= \wt{M}.
 $$
 \end{proof}

 As an immediate consequence of the definition of a coherent sheaf \FF\
 on $X$, we see that the sheaf
 $$
 \FFh := \phi^{\ast} \FF
 $$
 will be coherent as well:
 If 
 $$ 
 \OO_U^r \longrightarrow \OO_U^s \longrightarrow \FF_{|U} \longrightarrow 0
 $$
 is exact, so is
 $$
 \OO_{\Uh}^r \longrightarrow \OO_{\Uh}^s \longrightarrow 
 \phi^{\ast} \FF_{|\Uh} \longrightarrow 0,
 $$
 since $\phi^{-1}$ is exact and tensoring is a right exact functor.

 There is a famous theorem by Serre usually referred to as {\small
   GAGA}, since it is contained in his paper ``{\em G{\'e}om{\'e}trie
   alg{\'e}brique et g{\'e}om{\'e}trie analytique}'' \cite{gaga}.

 \begin{theorem}[Serre, {\cite[B.2.1, p. 440]{hartshorne}}]
   \label{gagathm}
   Let $X$ be a {\bf projective} scheme over \C. Then the map
   $$
   \FF \mapsto \FFh
   $$
   induces an equivalence between the category of coherent sheaves on
   $X$ and the category of coherent sheaves on \Xh. Furthermore, the
   natural map {\rm (\ref{gaganatmorph}) }
   $$
   \h^i(X;\FF) \stackrel{\sim}{\longrightarrow}
   \h^i(\Xh;\FFh)
   $$
   is an isomorphism for all $i$.
 \end{theorem}

 Let us state a corollary of Theorem \ref{gagathm}, which is not
 included in \cite{hartshorne}.

 \begin{corollary}
 \label{gagahyper}
 In the situation of Theorem \ref{gagathm}, we also have a natural isomorphism for
 hypercohomology for all $i$
 $$
 \H^i(X;\FF^{\bullet}) \iso \H^i(\Xh;\FFh^{\bullet}),
 $$
 where $\FF^{\bullet}$ is a bounded complex of coherent sheaves on $X$.
 \end{corollary}

 Here we only require the boundary morphisms of $\FF\bul$ to be
 morphisms of sheaves of abelian groups. They do not need to be
 $\OO_X$-linear.

 Before we begin proving Corollary \ref{gagahyper},
 we need some homological algebra.

 \begin{lemma}
 \label{lemma3}
 Let \AAA\ be an abelian category and 
 \begin{equation}
 \label{lemma3m}
 \f\bul[0] \longto \g^{\Bul,\Bul},
 \end{equation}
 a morphism of double complexes of \AAA-objects (cf. the appendix), 
 where 
 \begin{itemize}
 \item
 $\f\bul[0]$ is a double complex concentrated in the zeroth
 row with $\f\bul$ being a complex vanishing below degree zero, 
 i.e. $\f^n=0$ for $n<0$, and
 \item $\g^{\Bul,\Bul}$ is a double complex living only in
   non-negative degrees.
 \end{itemize}
 If for all $q\in\Z$
 $$
 0 \to \f^q \to \g^{0,q} \to \g^{1,q} \to \, \cdots
 $$
 is a resolution of $\f^q$, then the map of total
 complexes induced by {\rm (\ref{lemma3m}) }
 $$
 \f\bul \stackrel{\sim}{\longto} \tot \g^{\Bul,\Bul}
 $$
 is a quasi-isomorphism.
 \end{lemma}

 \begin{proof}
 From (\ref{lemma3m}) we obtain a morphism of spectral sequences
 \begin{align*}
 \hh^p_{\rm II}\,\hh^q_{\rm I}\, \f\bul[0] & \Rightarrow \hh^n\,\f\bul \\
 \downarrow \hspace{0.5cm} & \hspace{0.7cm} \downarrow \\
 \hh^p_{\rm II}\,\hh^q_{\rm I}\, \g^{\Bul,\Bul} & \Rightarrow
 \hh^n\,\tot\g^{\Bul,\Bul} .
 \end{align*}
 Both spectral sequences degenerate because of
 $$
 \hh^p_{\rm II}\,\hh^q_{\rm I}\, \f\bul[0] =
 \begin{cases}
 \hh^p \, \f\bul & \text{if $q=0$,} \\
 0 & \text{else} 
 \end{cases}
 $$
 and
 $$
 \hh^p_{\rm II}\,\hh^q_{\rm I}\, \g^{\Bul,\Bul} =
 \begin{cases}
 \hh^p \, \f\bul & \text{if $q=0$,} \\
 0 & \text{else} 
 \end{cases}
 $$
 and so we get an isomorphism on the initial terms. Hence we also have
 an isomorphism on the limit terms and
 our assertion follows.
 \end{proof}

 \begin{remark}
 \label{lemma3dual}
 A similar statement holds with $\f\bul[0]$ considered as a double
 complex concentrated in the zeroth column.
 \end{remark}

 {\em Godement resolutions.}
 As we also need $\Gamma$-acyclic resolutions 
 that behave functorial
 in the proof of 
 Corollary \ref{gagahyper}, 
 we now describe the concept of Godement resolutions.
 \label{godement}%
 \index{Godement resolution}%
 For any sheaf \FF\ on $X$ (or \Xh) define
 $$
 g(\FF) := \prod_{x\in \vert X \vert} i_{\ast} \FF_x,
 $$
 where $i: x \inclusion X$ denotes the closed immersion of the point
 $x$. The sheaf $g(\FF)$ is flabby and we have a natural inclusion
 $$
 \FF \rightarrowtail g(\FF).
 $$
 Setting 
 $\GG^{i+1} := g \left(\coker (\GG^{i-1} \rightarrow \GG^i)\right)$
 with $\GG^{-1} := \FF$ and $\GG^0:=g(\FF)$ gives an exact sequence
 $$
 0 \longrightarrow \FF \longrightarrow \GG^0 \longrightarrow \GG^1
 \longrightarrow \dots\, .
 $$
 We define the {\em Godement resolution} of \FF\ to be 
 $$
 \GG_{\FF}^{\bullet} := 0 \longrightarrow \GG^0 \longrightarrow \GG^1 
 \longrightarrow \dots \, .
 $$
 \notation{G}{$\GG_{\FF}\bul$}{Godement resolution of a sheaf $\FF$}
 It is $\Gamma$-acyclic and functorial in \FF. Extending this
 definition to a bounded complex $\FF^{\bullet}$
 \begin{equation}
 \begin{split}
 & \GG_{\FF^{\bullet}}^{\bullet,q} := \GG_{\FF^q}^{\bullet} \\
 & \GG_{\FF^{\bullet}}^{\bullet} := {\rm tot}\, \GG_{\FF\bul}^{\Bul,\Bul}
 \end{split}
 \end{equation}
 \notation{G}{$\GG_{\FF\bul}\bul$}{Godement resolution of a complex of
 sheaves $\FF\bul$}
 yields a map of double complexes
 $$
 \FF\bul[0] \longto \GG^{\Bul,\Bul}
 $$
 and a quasi-isomorphism by Lemma \ref{lemma3}
 $$
 \FF\bul \stackrel{\sim}{\longrightarrow} \GG_{\FF\bul}\bul .
 $$
 Let $x\in |X|$ be a closed point of $X$. The map $\phi$ from
 (\ref{phi}) gives a commutative square
 \begin{gather*}
 \{x\}^{\rm an} \stackrel{i}{\inclusion} \Xh \\
 {\scriptstyle \phi} \downarrow \hspace{1.35cm} 
 \downarrow {\scriptstyle\phi} \,\, \\[-0.5mm]
 \{x\}\phantom{^{\rm an}} \stackrel{i}{\inclusion} X\phantom{^{\rm an}}
 \end{gather*}
 and a natural morphism (cf. \cite[II.5, p.~110]{hartshorne})
 $$
 \FF_x \longto \phi_{\ast} \phi^{\ast} \FF_x,
 $$
 where we consider the stalk $\FF_x$ as a sheaf on $\{x\}$.
 This map induces another map
 $$
 i_{\ast} \FF_x \longto i_{\ast} \phi_{\ast} \phi^{\ast} \FF_x
 = (i\circ\phi)_{\ast}\phi^{\ast}\FF_x
 = (\phi\circ i)_{\ast}\phi^{\ast}\FF_x
 = \phi_{\ast} i_{\ast} \phi^{\ast} \FF_x,
 $$
 and yet a third map (cf. loc. cit.)
 \begin{equation}
 \label{hyperswitch}
 \epsilon: \phi^{\ast} i_{\ast} \FF_x \longrightarrow i_{\ast}
 \phi^{\ast} \FF_x,
 \end{equation}
 which is an isomorphism, as can be seen on the stalks
 \begin{align*}
 & \epsilon_x : \phi^{\ast} \FF_x \stackrel{\sim}{\longrightarrow}
   \phi^{\ast}\FF_x, \\
 & \epsilon_y : \quad 0 \longrightarrow 0,\qquad \text{for $y \ne x$}.
 \end{align*}
 Consequently, 
 $$
 \textstyle
 \phi^{\ast}g(\FF)=
 \phi^{\ast}\prod_x i_{\ast} \FF_x =
 \prod_x \phi^{\ast} i_{\ast}\FF_x \stackrel{(\ref{hyperswitch})}{=}
 \prod_x i_{\ast} \phi^{\ast}\FF_x =
 \prod_x i_{\ast} (\phi^{\ast}\FF)_x =
 g(\phi^{\ast}\FF_x).
 $$
 Thus we get first 
 $$
 \phi^{\ast}\GG_{\FF^p}\bul = \GG_{\phi^{\ast}\FF^p}\bul ,
 $$
 then
 $$
 \phi^{\ast}\GG_{\FF\bul}\bul = \GG_{\phi^{\ast}\FF\bul}\bul .
 $$
 This gives us a natural map of hypercohomology groups
 \begin{multline}
 \label{gagahypernatmorph}
 \H^{i}(X;\FF\bul)=\hh^i\Gamma(X;\GG_{\FF\bul}\bul)=
 \hh^i\Gamma(\Xh;\phi^{-1}\GG_{\FF\bul}\bul)\\
 \longrightarrow\\
 \hh^i\Gamma(\Xh;\phi^{\ast}\GG_{\FF\bul}\bul)= 
 \hh^i\Gamma(\Xh;\GG_{\phi^{\ast}\FF\bul}\bul)=
 \hh^i\Gamma(\Xh;\GG_{\FF\bul_{\rm an}}\bul)=
 \H^i(\Xh;\FFh\bul).
 \end{multline}

 {\em Proof of Corollary \ref{gagahyper}.}
 We claim that the natural map (\ref{gagahypernatmorph}) is an
 isomorphism for all $i$ if $X$ is projective.

 If $\FF\bul$ has has length one, Theorem \ref{gagathm} tells us that
 this is indeed true. So let us assume that (\ref{gagahypernatmorph})
 is an isomorphism for all complexes of length $\le n$ and let
 $\FF\bul$ be a complex of coherent sheaves on $X$ of length $n+1$. 
 W.l.o.g. $\FF^{n+1} \neq 0$ but $\FF^p=0$ for $p>n+1$.
 We write 
 $\sigma_{\le n}\FF\bul$ 
 for the complex $\FF\bul$ cut off above degree $n$.
 The short exact sequence 
 $$
 0 \longrightarrow \FF^{n+1}[-n-1] \longrightarrow \FF\bul
 \longrightarrow \sigma_{\le n} \FF\bul \longrightarrow 0
 $$
 remains exact if we take the inverse image along $\phi$
 $$
 0 \longrightarrow \FFh^{n+1}[-n-1] \longrightarrow \FFh\bul
 \longrightarrow \sigma_{\le n} \FFh\bul \longrightarrow 0 .
 $$
 Using the naturality of (\ref{gagahypernatmorph}), we obtain the
 following ``ladder'' with commuting squares
 $$
 \begin{CD}
 \cdots \to \, @.
 \h^{-n-1+i}(X;\FF^{n+1}) @>>> 
 \H^i(X;\FF\bul) @>>> 
 \H^i(X;\sigma_{\le n}\FF\bul) 
 @. \, \to \cdots \\
 @.
 @VV{\wr}V
 @VVV
 @VV{\wr}V
 @. \\
 \cdots \to \, @.
 \h^{-n-1+i}(\Xh;\FFh^{n+1}) @>>> 
 \H^i(\Xh;\FFh\bul) @>>> 
 \H^i(\Xh;\sigma_{\le n}\FFh\bul) 
 @. \, \to \cdots
 \end{CD}
 $$
 and the induction step follows from the $5$-lemma. {\hfill $\Box$}

\section{Algebraic \deRham Theory}
\label{alg}

In this section, we define the algebraic \deRham cohomology
$\hDR{\Bul}(X,D/k)$ of a smooth variety $X$ over a field $k$ and a
normal-crossings-divisor $D$ on $X$ (cf. definitions \ref{algdr1},
\ref{algdr2}, and \ref{algdr3}). We also give some working tools for this
cohomology: a base change theorem (Proposition \ref{basechange}) and
two spectral sequences
(Corollary \ref{specseqsmooth} and Proposition \ref{specseqdiv}).

\subsection{Classical \deRham Cohomology}
\label{drclassic}

This subsection only serves as a motivation for the following giving an overview
of classical \makebox{\deRham}\ theory. For a complex manifold, 
analytic \deRham, 
complex, and singular cohomology are defined and shown to be equal to
classical \deRham cohomology. All material presented here will be
generalized to a relative setup later on.

\medskip
Let $M$ be a complex manifold. We recall two standard exact sequences,
\begin{enumroman}
\item 
the {\em analytic \deRham complex} of holomorphic differential forms on $M$
\index{DeRham complex@\DeRham complex!analytic}
\notation{O}{$\om{M}$}{analytic \deRham complex of a complex manifold $M$}
$$
0 \longrightarrow \C_M \longrightarrow \Omega^0_M
\stackrel{\del}{\longrightarrow} \Omega^1_M
\stackrel{\del}{\longrightarrow} \ldots , \text{ and}
$$
\item 
the {\em classical \deRham complex} of smooth \C-valued 
differential forms on $M$
\index{DeRham complex@\DeRham complex!classical}
\notation{E}{$\EE\bul_M$}{classical \deRham complex of a complex manifold $M$}
$$
0 \longrightarrow \C_M \longrightarrow \EE^0_M
\stackrel{d}{\longrightarrow} \EE^1_M
\stackrel{d}{\longrightarrow} \ldots , \phantom{\text{ and}}
$$
\end{enumroman}
where $\C_M$ is the constant sheaf with fibre \C\ on $M$.
\notation{C}{$\C_X$}{constant sheaf with fibre $\C$ on a space $X$}

In both cases exactness is a consequence of the respective Poincar\'e
lemmas \cite[4.18, p.~155]{warner} and \cite[p.~25]{griffiths_harris}.

Now consider the commutative diagram
\begin{gather*}
\C_M [0] = \C_M [0] \\
\downiso \hspace{1cm} \downiso \\
\om{M} \inclusion \EE_M\bul .
\end{gather*}
We can rephrase the exactness of the sequences (i) and (ii) by saying that the vertical maps
are quasi-isomorphisms. We indicate quasi-isomorphisms by a tilde. Hence the natural inclusion
$$
\om{M} \inclusion \EE_M\bul 
$$
is a quasi-isomorphism  as well. Therefore the
hypercohomology of the two complexes coincides
$$
\H\bul(M;\om{M}) = \H\bul(M;\EE_M\bul).
$$
The sheaves $\EE^p_M$ are fine, since they admit a partition of
unity. In particular they are acyclic for the global section functor
$\Gamma(M,?)$ and we obtain
$$
\H\bul(M;\EE_M\bul) = \hh\bul\Gamma(M;\EE_M\bul).
$$
The right-hand-side is usually called the {\em classical \deRham
cohomology} 
\index{DeRham cohomology@\DeRham cohomology!classical}
of $M$, denoted
$$
\hDR{\Bul}(M;\C).
$$
The equalities above give
$$
\hDR{\Bul}(M;\C) = \H\bul(M;\om{M}).
$$
We refer to the right-hand-side $\H\bul(M;\om{M})$ as 
{\em analytic \deRham cohomology},
\index{DeRham cohomology@\DeRham cohomology!analytic}%
for which we want to use the same symbol $\hDR{\Bul}(M;\C)$.%
\notation{H}{$\hDR{\Bul}(M;\C)$}{analytic \deRham cohomology of a %
complex manifold $M$}%
The hypercohomology $\H\bul(M;\om{M})$ turns out to be a good candidate for
generalizing \deRham theory to algebraic varieties.

Both variants of \deRham cohomology agree with {\em complex
  cohomology}
\index{Complex cohomology}
\label{complcoho}
$$
\h\bul(M;\C) := \h\bul(M;\C_M) = \H\bul(M;\C_M[0]) .
$$
\notation{H}{$\h\bul(M;\C)$}{complex cohomology of a space $M$}

\medskip
We have yet a third resolution of the constant sheaf $\C_M$ given by
the {\em complex of singular cochains}:
For any open set $U\subseteq M$, we write $\c^p\sing(U;\C)$ 
  \notation{C}{$\c^p\sing(M;\C)$}{\C-vector space of singular
    $p$-cochains of a
  space $M$}
for the vector space of singular $p$-cochains on $U$ with coefficients in
\C.
The sheaf of singular $p$-cochains is now defined as
\notation{C}{$\CC^p\sing(M;\C)$}{sheaf of singular $p$-cochains of a
  space $M$}
$$
\CC^p\sing(M;\C) : U \mapsto \c^p\sing(U;\C)
\quad\text{for}\quad U\subseteq M\quad\text{open}
$$
with the obvious restriction maps.
The sheaves $\CC^p\sing(M;\C)$ are flabby: The restriction maps are
the duals of injections between vector spaces of singular $p$-chains 
and the functor $\Hom_{\C}(?,\C)$ is exact. In
particular these sheaves are acyclic for the global section functor
$\Gamma(M;?)$ and we obtain
\label{singcoho}
\index{Singular cohomology}
$$
\H\bul(M;\CC\bul\sing(M;\C)) =
\hh\bul\Gamma(M;\CC\bul\sing(M;\C))=
\hh\bul \c\bul\sing(M;\C)=
\h\bul\sing(M;\C).
$$

By the following lemma, we conclude $\h\bul(M;\C)=\h\bul\sing(M;\C)$.

\begin{lemma}
\label{singres}
For any locally contractible, locally path-connected topological space
$M$ the sequence
$$
0 \longrightarrow \C_M \longrightarrow \CC^0\sing(M;\C)
\longrightarrow \CC^1\sing(M;\C) \longrightarrow \, \ldots \, .
$$
is exact.
\end{lemma}

\begin{proof}
Note first that 
$ \C_M = \hh^0\,\CC\bul\sing(M;\C)$,
since $M$ is locally path-connected.
For the higher cohomology sheaves
$ \hh^p\, \CC\bul\sing(M;\C),\, p > 0$,
we observe that any element $s_x$ of the stalk $\hh^p\,\CC\bul\sing(M;\C)_x$
at $x\in M$ not only
lifts to a section $s$ of $\c^p\sing(U;\C)$ for some contractible
open subset $x\in U \subset M$, but that we can assume $s$ to be a
cocycle by eventually shrinking $U$. Now this $s$ is also a coboundary
because of
$$
\hh^p \c\bul\sing(U;\C) = \h^p\sing(U;\C) = 0 \quad\text{for all}\quad  p>0.
$$
\end{proof}

\medskip
We summarize this subsection in the following proposition.

\begin{proposition}
\label{classiso}
Let $M$ be a complex manifold. Then we have a chain of natural
isomorphisms between the various cohomology groups defined in this 
subsection
$$
\hDR{\Bul}(M;\C) \iso
\H\bul(M;\om{M}) \iso
\h\bul(M;\C) \iso
\h\bul\sing(M;\C) .
$$
\end{proposition}

\subsection{Algebraic \deRham Cohomology}
\label{algDR}

Let $X$ be a smooth variety defined over a field $k$ and $D$ a divisor with
normal crossings on $X$; where having normal crossings
means, that locally $D$ looks like a collection of coordinate
hypersurfaces, or more precisely:

\begin{definition}[Divisor with normal crossings, 
  {\cite[p. 391]{hartshorne}}]
  \index{Divisor!with normal crossings}%
  \label{normalcrossings}%
  A divisor $D \subset X$ is said to have normal crossings, if each
  irreducible component of $D$ is nonsingular and whenever $s$
  irreducible components $D_1, \ldots, D_s$ meet at a closed point $P$, then
  the local equations $f_1, \ldots, f_s$ of the $D_i$ form part of a
  regular system of parameters $f_1,\ldots,f_d$ at $P$.
\end{definition}

It is proved in \cite[12, p.~78]{matsumura} that in this case the
$f_1,\ldots,f_d$ are linearly independent modulo $m_P^2$, where $m_P$
is the maximal ideal of the local ring $\OO_{X,P}$ at $P$.
By the inverse function theorem for holomorphic functions
\cite[Thm.~6, p.~20]{gunning}, we find in a neighbourhood of any $P\in\Xh$ a
holomorphic chart $z_1,\ldots,z_d$ such that $\Dh$ is given as the
zero-set $\{z_1\cdot\ldots\cdot z_s=0\}$.

We are now going to define algebraic \deRham cohomology groups
$$
\hDR{\Bul}(X/k),\quad \hDR{\Bul}(D/k) \quad \text{and} 
\quad \hDR{\Bul}(X,D/k) .
$$

\begin{remark}
  In \cite{hartshorne_71}, algebraic \deRham cohomology is defined for
  varieties with arbitrary singularities. However, the relative version of
  algebraic \deRham cohomology discussed in \cite{hartshorne_71} 
  deals with morphisms of varieties, not pairs of them.
\end{remark}

\subsubsection{The Smooth Case} 

\begin{definition}[Algebraic \deRham cohomology for a smooth variety]
\label{algdr1}
\index{DeRham cohomology@\DeRham cohomology!algebraic}
  We set
  $$
  \hDR{\Bul}(X/k) := \H\bul(X;\Omega_{X/k}\bul),
  $$
  where $\Omega_{X/k}\bul$ is the complex of algebraic differential
  forms on the smooth variety $X$ over $k$ (cf. \cite[II.8, p.~175]{hartshorne}).
  \index{DeRham complex@\DeRham complex!algebraic}
  \notation{O}{$\Omega_{X/k}\bul$}{algebraic \deRham complex of a
    variety $X$}
  \notation{H}{$\hDR{\Bul}(X/k)$}{algebraic \deRham cohomology of a smooth
    variety $X$}
\end{definition}

\subsubsection{The Case of a Divisor with Normal Crossings}

\begin{wrapfigure}{r}{4.2 cm}
\hbox{\vspace{-3mm}
\epsfig{file=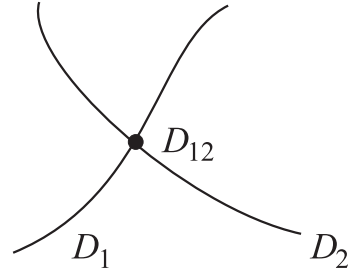,width=4.5 cm}}
\caption{Divisor with normal crossings}
\label{div}
\end{wrapfigure}
We write $D=\sum_{i=1}^r D_i$ as a sum of its irreducible components and
use the short-hand notation (cf. Figure \ref{div})
\label{shorthand}
$$
D_I := D_{i_0 \, \cdots \, i_p} := \bigcap_{k=0}^p D_{i_k} 
\quad\text{for}\quad I=\{i_0,\ldots,i_k\}.
$$
\notation{D}{$D_{i_0\,\ldots\,i_p}$}{intersection of
  $D_{i_0},\cdots,D_{i_p}$}
\notation{D}{$D_I$}{intersection of $D_i$ for $i\in I$}

Associated to the decomposition $D=\sum D_i$ is a simplicial scheme
$$
D\bul := 
\coprod D_a 
\begin{smallmatrix} 
\longleftarrow \\ 
\longleftarrow 
\end{smallmatrix}
\coprod D_{ab}
\begin{smallmatrix} 
\longleftarrow \\ 
\longleftarrow \\ 
\longleftarrow
\end{smallmatrix}
\coprod D_{abc} \cdots ,
$$
this notion being defined as follows. 
\paragraph{Simplicial sets and simplicial schemes.} 
\index{Simplicial scheme}
We consider a category
${\fk{Simplex}}$ 
\notation{S}{$\fk{Simplex}$}{category with objects $[m]$ and %
  $<$-preserving maps}
with objects $[m]:=\{0,\ldots,m\}$, $m\in\N_0$, and
$<$-order-preserving 
maps as morphisms. 
(Thus the existence of a map $f: [m] \to [n]$ implies $m\le n$.)

This category can be thought of as the
prototype of a simplicial complex: Let $K$ be a simplicial complex
and assign to $[m]$
the set of $m$-simplices of $K$. Any morphism $f: [m] \to [n]$ maps each
$n$-simplex to its $m$-dimensional ``$f$-face'' 
$$
f^{\ast}(e_0,\ldots,e_n) = (e_{f(0)},\ldots,e_{f(m)}).
$$
Thus $K$ can be described by a contravariant functor to 
the category $\fk{Sets}$ of sets
$$
\fk{Simplex} \longrightarrow \fk{Sets}.
$$
A contravariant functor from $\fk{Simplex}$ to the category of sets,
schemes, $\ldots$ is called a {\it simplicial set}, {\it simplicial
  scheme} and so forth \cite[Ch. 1, 2.2, p. 9]{gelfand_manin}.
(There, as a slight generalization, the maps $f: [m] \to
[n]$ are only required to be non-decreasing instead of 
being $<$-order-preserving, thus
allowing for ``degenerate simplices''.)

\medskip
\label{simpscheme}
In the above situation we define a simplicial scheme 
$D\bul:\fk{Simplex}\to \fk{Schemes}$ 
\notation{D}{$D\bul$}{simplicial scheme associated to a divisor $D$}
by the assignment
\begin{align*}
[m] & \mapsto \coprod_{1\le i_0 < \cdots < i_m \le r} D_{i_0 \, \cdots \,
  i_m} \\
\left( f: [m] \to [n] \right) & \mapsto 
  \left( D\bul(f) : \coprod D_{i_0 \, \cdots \, i_n} \to \coprod D_{i_0 \,
  \cdots \, i_m} \right);
\end{align*}
here $D\bul(f)$ is the sum of the natural inclusions
$$
D_{i_0 \, \cdots \, i_n} \inclusion D_{i_{f(0)} \, \cdots \, i_{f(m)}} .
$$
If $\delta^m_l : [m] \to [m+1]$ denotes the unique
$<$-order-preserving map, whose image does not contain $l$, we can
represent $D\bul$ by a diagram
\begin{equation}
\coprod_{a=1}^r D_a 
\begin{smallmatrix} 
\begin{CD} @<{D\bul(\delta^0_0)}<< \end{CD} \\
\begin{CD} @<{D\bul(\delta^0_1)}<< \end{CD}
\end{smallmatrix}
\coprod_{1\le a < b \le r} D_{ab}
\begin{smallmatrix} 
\begin{CD} @<{D\bul(\delta^1_0)}<< \end{CD} \\
\begin{CD} @<{D\bul(\delta^1_1)}<< \end{CD} \\
\begin{CD} @<{D\bul(\delta^1_2)}<< \end{CD} 
\end{smallmatrix}
\coprod_{1 \le a < b < c \le r} D_{abc} \cdots .
\tag{$D\bul$}
\end{equation}

Now the differentials come into play again. We consider the
contravariant functor $i_{\ast} \Omega\bul_{?/k}$ defined for
subschemes $Z$ of $\Supp(D)$
$$
i_{\ast}\Omega\bul_{?/k} : Z \mapsto i_{\ast}\Omega\bul_{Z/k} .
$$
where $i$ stands for the natural inclusion $Z \stackrel{i}{\inclusion} \Supp(D)$.
We enlarge the scope of $i_{\ast}\Omega\bul_{?/k}$
to schemes of the form $\coprod_a Z_a$, $Z_a\subseteq D$
by making the convention
$$
\left(i_{\ast} \Omega\bul_{?/k}\right)\left(\coprod_a Z_a\right)
:= \bigoplus_a\, i_{\ast} \Omega\bul_{Z_a/k} .
$$
Composing $D\bul$ with the contravariant functor $i_{\ast} \Omega\bul_{?/k}$ yields a diagram
$$
\bigoplus_{a=1}^r \iast \Omega\bul_{D_a/k} 
\begin{smallmatrix} 
\begin{CD} @>{d^0_0}>> \end{CD} \\
\begin{CD} @>{d^0_1}>> \end{CD} \\
\end{smallmatrix}
\bigoplus_{1\le a < b \le r} \iast\Omega\bul_{D_{ab}/k}
\begin{smallmatrix} 
\begin{CD} @>{d^1_0}>> \end{CD} \\
\begin{CD} @>{d^1_1}>> \end{CD} \\
\begin{CD} @>{d^1_2}>> \end{CD} \\
\end{smallmatrix}
\bigoplus_{1 \le a < b < c \le r} \iast\Omega\bul_{D_{abc}/k} \cdots\, ,
$$
where we have written $d^m_l$ for $(\iast\Omega\bul_{?/k} \circ D\bul)
(\delta^m_l)$. 
Summing up these maps $d^m_l$ with alternating signs
$$
d^m := \sum_{l=0}^m (-1)^l d^m_l
$$
gives a new diagram
$$
\bigoplus_{a=1}^r \iast\Omega\bul_{D_a/k}  
\begin{CD} @>{d^0}>> \end{CD}
\bigoplus_{1\le a < b \le r} \iast\Omega\bul_{D_{ab}/k}
\begin{CD} @>{d^1}>> \end{CD} \\
\bigoplus_{1 \le a < b < c \le r} \iast\Omega\bul_{D_{abc}/k} \cdots\, , 
$$
which turns out to be a double complex (cf. the appendix) as a little calculation shows
\begin{align*}
d^{m+1} \circ d^m = &
\sum_{l=0}^{m+1} \sum_{k=0}^m (-1)^{l+k} \, d^{m+1}_l \circ d^m_k \\
= & \sum_{0\le k < l \le m+1} (-1)^{l+k} \, d^{m+1}_l \circ d^m_k \, + 
\hspace{-1cm}
\underbrace{\sum_{0\le l \le k \le m} 
  (-1)^{l+k} \, d^{m+1}_l \circ d^m_k}_{
\substack{ 
\displaystyle
\hspace{1cm} = \sum_{0\le l \le k \le m} (-1)^{l+k} \, 
  d^{m+1}_{k+1} \circ d^m_l \quad\qquad \\
\displaystyle
\hspace{1cm} = \sum_{0\le k' < l' \le m+1} 
  (-1)^{k'+l'-1} \, d^{m+1}_{l'} \circ d^m_{k'} 
}} \\
= & \, 0 .
\end{align*}
We denote this double complex by $\Omega^{\Bul,\Bul}_{D\bul/k}$%
\notation{O}{$\Omega^{\Bul,\Bul}_{D\bul/k}$}{algebraic double %
\deRham complex of $D\bul$}%
and its total
complex (cf. the appendix) by 
$$
\wt{\Omega}\bul_{D/k} := \tot \Omega^{\Bul,\Bul}_{D\bul/k} .
$$
\label{defwom}
\notation{O}{$\wom{D/k}$}{algebraic \deRham complex of a divisor $D$} 
\begin{definition}[Algebraic \deRham cohomology for a divisor]
\label{algdr2}
\index{DeRham cohomology@\DeRham cohomology!algebraic}
  If $D$ is a divisor with normal crossings on a smooth variety $X$
  over $k$, we define
  its {\em algebraic \deRham cohomology} 
  by
  $$
  \hDR{\Bul}(D/k) := \H\bul(D;\wt{\Omega}\bul_{D/k}) .
  $$
\notation{H}{$\hDR{\Bul}(D/k)$}{algebraic \deRham cohomology of a divisor $D$}
\end{definition}

\begin{remark}
  Actually, 
  we do not need that $D\subseteq X$ has codimension $1$
  for this definition;
  but $D$ having normal crossings is essential.
\end{remark}

\subsubsection{The Relative Case}

The natural restriction maps
$$
\Omega\bul_{X/k} \to  (D_i \inclusion X)_{\ast} \Omega\bul_{D_i/k}
$$
sum up to a natural map of double complexes
\begin{equation}
\label{OmXD}
\Omega\bul_{X/k}[0] \to \left(\Supp(D) \inclusion X\right)_{\ast}
\Omega^{\Bul,\Bul}_{D\bul/k},
\end{equation}
where we view $\Omega\bul_{X/k}[0]$ as a double complex concentrated
in the zeroth column. Somewhat sloppy, we will often write
$i_{\ast}\FF$ instead of $\left(\Supp(D) \inclusion X\right)_{\ast}
\FF$.
Taking the total complex in (\ref{OmXD}) yields a natural map
$$
f: \om{X/k} \to i_{\ast} \wom{D/k},
$$
for whose mapping cone (cf. the appendix)
$$
M_f = \iast\wom{D/k}[-1]\oplus\om{X/k} 
\text{ with differential } 
\left(\begin{matrix}
-d_D & f \\
0 & d_X
\end{matrix}\right)
$$
we write $\wom{X,D/k}$.
\notation{O}{$\wom{X,D/k}$}{algebraic \deRham complex of a pair $(X,D)$}
\label{defwomrel}

\begin{definition}[Relative algebraic \deRham cohomology]
\label{algdr3}
\index{DeRham cohomology@\DeRham cohomology!algebraic}
  For $X$ a smooth variety over $k$ and $D$ a divisor with normal
  crossings on $X$, we define the {\em relative algebraic \deRham cohomology}
  of the pair $(X,D)$ by
  $$
  \hDR{\Bul}(X,D/k) := \H\bul(X; \wom{X,D/k}).
  $$
\notation{H}{$\hDR{\Bul}(X,D/k)$}{algebraic \deRham cohomology of a
  pair $(X,D)$}
\end{definition}

These definitions may seem a bit technical at first
glance. However, they will hopefully become more transparent in the
sequel, when the existence of a long exact sequence in algebraic
\deRham cohomology and various comparison isomorphisms are proved.

\subsection{Basic Lemmas}

In this subsection we gather two basic lemmas for the ease of reference.

\begin{lemma}
\label{injclosedimm}
If $i: Z \inclusion X$ is a closed immersion and $\II$ an injective
sheaf of abelian groups on $Z$, then $\iast\II$ is also injective.
\end{lemma}

\begin{proof}
Let $\AA \injection \BB$ be an injective morphism of sheaves of
abelian groups on $X$. Since $i^{-1}$ is an exact functor, the induced
map
$$
i^{-1}\AA \injection i^{-1}\BB
$$
will be injective as well. Now the adjoint property of $i^{-1}$
provides us with the following commutative square
\begin{gather*}
\Hom_X(\BB,\iast\II) \iso \Hom_Z(i^{-1}\BB,\II) \\
\downarrow \hspace{2.75cm} \downarrow \\
\Hom_X(\AA,\iast\II) \iso \Hom_Z(i^{-1}\AA,\II) .
\end{gather*}
Because $\II$ is injective, the map on the right-hand-side is
surjective. Hence so is the left-hand map, which in turn implies the
injectivity of $\iast\II$, since $\AA$ and $\BB$ were arbitrary.
\end{proof}

\begin{lemma}
\label{dirimhyper}
If $i: Z \inclusion X$ is a closed immersion and $\FF\bul$ a complex
of sheaves on $Z$, which is bounded below,  then 
there is a natural isomorphism
$$
\H\bul(X;\iast\FF\bul) \iso \H\bul(Z;\FF\bul).
$$
\end{lemma}

\begin{proof}
The point is that $\iast$ is an exact functor for closed immersions,
which can be easily checked on the stalks. Thus 
if $\FF\bul \injection \II\bul$
is a quasi-isomorphism between $\FF\bul$ and a complex $\II\bul$
consisting of injective sheaves, the induced map
$$
\iast\FF\bul\injection\iast\II\bul
$$
will be a quasi-isomorphism as well. But the sheaves $\iast\II^p$ are
injective by Lemma \ref{injclosedimm}, hence we can use
$\iast\II\bul$ to compute the hypercohomology of $\iast\FF\bul$
$$
\H\bul(X;\iast\FF\bul) = \hh\bul\Gamma(X;\iast\II\bul) .
$$
Now the claim follows from the identity of functors $\Gamma(X,\iast ?)
= \Gamma(Z;?)$
$$
\H\bul(X;\iast\FF\bul) = \hh\bul\Gamma(X;\iast\II\bul) =
\hh\bul\Gamma(Z;\II\bul)=\H\bul(Z;\FF\bul).
$$
\end{proof}

\subsection{The Long Exact Sequence in Algebraic \deRham Cohomology}
\label{longseqalgDR}

\begin{proposition}
We have a natural long exact sequence in algebraic \deRham
cohomology (as defined in Definition 
\ref{algdr1}, \ref{algdr2} and \ref{algdr3})
\begin{align*}
\cdots & \to \hDR{p-1}(D/k) \to \\
\hDR{p}(X,D/k) \to \, \hDR{p}(X/k) & \to \hDR{p}(D/k) \to \\
\hDR{p+1}(X,D/k) \to \, \cdots \, , \hspace{1cm} &
\end{align*}
where $X$ is a smooth variety over $k$ and $D$ a
normal-crossings-divisor on $X$.
\end{proposition}

\begin{proof} 
The short exact sequence of the mapping cone
$$
0 \to i_{\ast} \wom{D/k}[-1] \to \wom{X,D/k} \to \om{X/k} \to 0
$$
gives us a long exact sequence in hypercohomology
\begin{align*}
\cdots & \to \H^p(X;i_{\ast}\wom{D/k}[-1]) \to \\
\H^p(X;\wom{X,D/k}) \to \, \H^p(X;\om{X/k}) & \to \H^{p+1}(X;i_{\ast}\wom{D/k}[-1]) \to \\
\H^{p+1}(X;\wom{X,D/k}) \to \, \cdots \, , \hspace{1.5cm} &
\end{align*}
which we can rewrite as
\begin{align*}
\cdots & \to \H^{p-1}(X;i_{\ast}\wom{D/k}) \to \\
\hDR{p}(X,D/k) \to \, \hDR{p}(X/k) & \to \H^p(X;i_{\ast}\wom{D/k}) \to \\
\hDR{p+1}(X,D/k) \to \, \cdots \, . \hspace{1cm} & 
\end{align*}
By Lemma \ref{dirimhyper}, we see
$$
\H\bul(X;\iast\wom{D/k}) = \H\bul(D;\wom{D/k}) = \hDR{\Bul}(D/k)
$$
and our assertion follows.
\end{proof}

\subsection{Behaviour Under Base Change}

Let $X_0$ be a smooth variety defined over a field $k_0$ and $D_0$ a
divisor with normal crossings on $X_0$. For the base change to an
extension field $k$ of $k_0$ we write
$$
X:=X_0 \times_{k_0} k \text{ and } 
D:=D_0 \times_{k_0} k.
$$
Then we have the following proposition.

\begin{proposition}
\label{basechange}
In the above situation, there are natural
isomorphisms between algebraic \deRham cohomology groups
(as defined in definitions \ref{algdr1}, \ref{algdr2}, \ref{algdr3})
\begin{align*}
 \hDR{\Bul}(X/k) & \iso \hDR{\Bul}(X_0/k_0) \tensor{k_0} k, \\
 \hDR{\Bul}(D/k) & \iso \hDR{\Bul}(D_0/k_0) \tensor{k_0} k \quad \text{\rm and} \\
 \hDR{\Bul}(X,D/k) & \iso \hDR{\Bul}(X_0,D_0/k_0) \tensor{k_0} k .
\end{align*}
\end{proposition}

\begin{proof}
Denote by  $\pi: X \to X_0$ the natural projection. By
\cite[Prop.~II.8.10, p.~175]{hartshorne}, we have
$$
\om{X/k} = \pi^{\ast} \om{X_0/k_0} .
$$
Now the natural map
$$
\pi^{-1} \om{\Xo/k_0} \to \pi^{\ast} \om{\Xo/k_0} 
$$
factors through
$$ 
\pi^{-1} \om{\Xo/k_0} \to \pi^{-1} \om{\Xo/k_0} \tensor{k_0} k
$$
yielding a natural map
$$
\pi^{-1} \om{\Xo/k_0} \tensor{k_0} k \to \pi^{\ast} \om{\Xo/k_0} =
\om{X/k},
$$
which provides us with a morphism of spectral sequences
$$
\begin{matrix}
\hh^p \, \h^q(X; \pi^{-1} \om{\Xo/k_0} \tensor{k_0} k) & 
\Rightarrow & \H^n(X;\pi^{-1} \om{\Xo/k_0} \tensor{k_0} k) \\
\downarrow & & \downarrow \\
\hh^p \, \h^q(X;\pi^{\ast} \om{\Xo/k_0}) &
\Rightarrow & \H^n(X;\pi^{\ast} \om{\Xo/k_0}) \, .
\end{matrix}
$$ 
We can rewrite the first initial term using the exactness of
the functor $?\tensor{k_0} k$
\begin{align*}
\h^q(X;\pi^{-1} \Omega^p_{\Xo/k_0} \tensor{k_0} k) & = 
\h^q(X;\pi^{-1} \Omega^p_{\Xo/k_0}) \tensor{k_0} k \\
& = \h^q(X_0;\Omega^p_{\Xo/k_0}) \tensor{k_0} k .
\end{align*}
Since ``{\em cohomology commutes with flat base extension for
  quasi-coherent sheaves}'' {\cite[Prop.~III.9.3,
  p.~255]{hartshorne}}, the map
$$
\h^q(X_0;\Omega^p_{\Xo/k_0}) \tensor{k_0} k \longto
\h^q(X;\pi^{\ast} \Omega^p_{\Xo/k_0})
$$
is an isomorphism. Hence we also have an isomorphism on the limit
terms, for which we obtain
\begin{align*}
\H^n(X;\pi^{-1} \om{\Xo/k_0} \tensor{k_0} k) & =
\H^n(X;\pi^{-1} \om{\Xo/k_0} ) \tensor{k_0} k \\
& = \H^n(\Xo;\om{\Xo/k_0}) \tensor{k_0} k \\
& = \hDR{n}(X_0/k_0) \tensor{k_0} k 
\end{align*}
and
\begin{align*}
\H^n(X;\pi^{\ast} \om{X_0/k_0}) & = \H^n(X;\om{X/k}) \\
& = \hDR{n}(X/k).
\end{align*}
The other two statements can be proved analogously. (Observe that the
sheaves $\wt{\Omega}^p_{D_0/k_0}$ and $\wt{\Omega}^p_{X_0,D_0/k_0}$ are
quasi-coherent by \cite[Prop.~II.5.8(c), p.~115]{hartshorne}.)
\end{proof}

\subsection{Some Spectral Sequences}

Before we proceed in the discussion of algebraic \deRham cohomology, we
have to make a digression and provide some spectral sequences which we will
need later on.

\subsubsection{\v{C}ech Cohomology}
\label{chechcoho}

We adopt the notation from \cite[III.4]{hartshorne}:
Let $X$ be a topological space, 
$\UUU=(U_j)_{j\in J}$ 
\notation{U}{$\UUU$}{open covering $(U_j)_{j\in J}$}
an open covering of
$X$ and $\FF$ a sheaf of abelian groups on $X$. We assume $J$ to be
well-ordered. As usual, we use the short-hand notation
\notation{U}{$U_{i_0\,\ldots\,i_m}$}{intersection of
  $U_{i_0},\ldots,U_{i_m}$}
\notation{U}{$U_I$}{intersection of $U_i$ for $i\in I$}
$$
U_I := U_{i_0 \, \cdots \, i_m} := \bigcap_{k=0}^m U_{i_k} 
\quad\text{for}\quad I=\{i_0,\ldots,i_m\}.
$$
Recall that in this situation we have a {\em \v{C}ech functor}
defined on the category $\fk{Ab}(X)$ of sheaves of abelian groups on
$X$ (cf. \cite[III.4, p.~220]{hartshorne})

\notation{C}{$\CCh^p(\UUU,?)$}{\v{C}ech functor}
\index{Cech@\v{C}ech!functor}
\begin{align*}
\CCh^p(\UUU,?) : \fk{Ab}(X) & \to \fk{Ab}(X) \\
\FF & \mapsto \bigoplus_{\substack{|I|=p+1\\ I\subseteq J}} \iast \FF_{|U_I},
\end{align*}
where $i$ stands for the respective inclusions $U_I \inclusion X$,
and a differential
$$
d: \CCh^p(\UUU;\FF) \to \CCh^{p+1}(\UUU;\FF),
$$
which makes $\CCh\bul(\UUU,\FF)$ into a complex.
We explain this $d$ below.
We also need {\em \v{C}ech groups}
\notation{C}{$\Ch^p(\UUU;\FF)$}{\v{C}ech group of a sheaf $\FF$}
\index{Cech@\v{C}ech!group}
$$
\Ch^p(\UUU;\FF) := \Gamma(X;\CCh^p(\UUU;\FF))
$$
and of course {\em \v{C}ech cohomology}
\notation{H}{$\check{\rm H}^p(\UUU;\FF)$}{\v{C}ech cohomology of a sheaf
  $\FF$}
\index{Cech@\v{C}ech!cohomology}
$$
\check{\rm H}^p(\UUU;\FF) := \hh^p \, \Ch\bul(\UUU;\FF).
$$
Clearly, we can consider the \v{C}ech functor also for closed
coverings; but then \v{C}ech cohomology will not behave
particularly well.

\medskip
We can give the following fancy definition of the \v{C}ech complex
$\CCh\bul(\UUU;\FF)$.
The assignment
\begin{align*}
\UUU\bul: \fk{Simplex} & \to \fk{Schemes} \\
[m] & \mapsto \coprod_{i_0 < \cdots < i_m} U_{i_0 \, \cdots \, i_m} \\
\left( f: [m] \to [n] \right) & \mapsto 
\left( \UUU\bul(f) : \coprod U_{i_0 \, \cdots \, i_n} \to 
  \coprod U_{i_0 \, \cdots \, i_m} \right)
\end{align*}
provides us with a simplicial scheme $\UUU\bul$ (cf. page
\pageref{simpscheme}), 
\notation{U}{$\UUU\bul$}{simplicial scheme associated to an open covering $\UUU$}
and thus we have a diagram
$$
\coprod_a U_a 
\begin{smallmatrix} 
\begin{CD} @<{\UUU\bul(\delta^0_0)}<< \end{CD} \\
\begin{CD} @<{\UUU\bul(\delta^0_1)}<< \end{CD}
\end{smallmatrix}
\coprod_{a<b} U_{ab}
\begin{smallmatrix} 
\begin{CD} @<{\UUU\bul(\delta^1_0)}<< \end{CD} \\
\begin{CD} @<{\UUU\bul(\delta^1_1)}<< \end{CD} \\
\begin{CD} @<{\UUU\bul(\delta^1_2)}<< \end{CD} 
\end{smallmatrix}
\coprod_{a<b<c} U_{abc} \cdots .
$$
Denote by
$$
\eps: \coprod_{|I|=m+1} U_I \to X
$$
the natural morphism induced by $U_I \inclusion X$. Then 
$$
\eps_{\ast} \eps^{-1}\FF = \bigoplus_{|I|=m+1}  \iast
\FF_{|U_I},
$$ where $i: U_I \inclusion X$.
Composing $\UUU\bul$ with the contravariant functor
$$
\eps_{\ast} \eps^{-1} \FF : \coprod_{|I|=m+1} U_I \mapsto 
\bigoplus_{|I|=m+1}  \iast
\FF_{|U_I}
$$
yields a diagram
$$
\bigoplus_{a} \iast \FF_{xU_a} 
\begin{smallmatrix} 
\begin{CD} @>{d^0_0}>> \end{CD} \\
\begin{CD} @>{d^0_1}>> \end{CD}
\end{smallmatrix}
\bigoplus_{a<b} \iast \FF_{|U_{ab}}
\begin{smallmatrix} 
\begin{CD} @>{d^1_0}>> \end{CD} \\
\begin{CD} @>{d^1_1}>> \end{CD} \\
\begin{CD} @>{d^1_2}>> \end{CD} 
\end{smallmatrix}
\bigoplus_{a<b<c} \iast \FF_{|U_{abc}} \cdots 
$$
with maps
$$
d^m_l := \left( \eps_{\ast} \eps^{-1} \FF \circ \UUU\bul \right)
(\delta^m_l) .
$$
Summing up these maps with alternating signs
$$
d^m := \sum_{l=0}^m (-1)^l \, d^m_l
$$
gives us the \v{C}ech complex $\CCh\bul(\UUU;\FF)$.

\medskip
We have the following lemma.

\begin{lemma}[{\cite[III.4.2, p. 220]{hartshorne}}]
\label{chechres}
Let $X$ be a topological space, $\FF$ a sheaf of abelian groups on $X$,
and $\UUU=(U_j)_{j\in J}$ an open covering of $X$. 
Then the \v{C}ech complex $\CCh\bul(\UUU;\FF)$, as defined above, is a
resolution for $\FF$, i.e. the sequence
$$
0 \to \FF \injection \CCh^0(\UUU;\FF) \stackrel{d^0}{\longrightarrow} 
  \CCh^1(\UUU;\FF) \stackrel{d^1}{\longrightarrow} \, \cdots 
$$
is exact. 
In particular, we have a natural isomorphism
$$
\h\bul(X;\FF) \iso \H\bul(X;\CCh\bul(\UUU;\FF)) .
$$
\end{lemma}

\subsubsection{A Spectral Sequence for Hypercohomology of an Open Covering}

We generalize the well-known spectral sequence of an open covering
$\UUU$ of a space $X$ 
(where $\HH^q(\FF)$ is the presheaf $V \mapsto \h^q(V;\FF_{|V})$)
$$
\check{\h}^p(\UUU;\HH^q(\FF)) \Rightarrow \h^n(X;\FF)
$$
to the case of a complex $\FF\bul$ of abelian sheaves and
its hypercohomology.

Let $X$ be a topological space, $\UUU=(U_j)_{j\in J}$ an open covering
of $X$ 
and $\FF\bul$ a complex of sheaves of abelian groups on $X$, which is bounded
below.
We define a presheaf
$$
\HHH^p(\FF\bul) : V \mapsto \H^p(V;\FF\bul_{|V})
\quad\text{for}\quad V\subseteq X\quad\text{open}
$$
\notation{H}{$\HHH^p(\FF\bul)$}{sheafified hypercohomology of a
  complex of sheaves $\FF\bul$}
as a ``sheafified'' version of hypercohomology.

\begin{proposition}
\label{specseqopcov}
In the situation above, 
the following spectral sequence converges
$$
E_2^{p,q} := \check{\h}^p(\UUU;\HHH^q(\FF\bul)) \Rightarrow E_{\infty}^n :=  \H^n(X;\FF\bul) .
$$
\end{proposition}

\begin{proof}
We choose a quasi-isomorphism $\FF\bul\stackrel{\sim}{\injection} \II\bul$
with $\II\bul$ being a complex of flabby sheaves. Note that the double
complex $\CCh\bul\II\bul := \CCh\bul(\UUU;\II\bul)$ consists of flabby
sheaves, as well.

Now we consider the two spectral sequences of a double complex
(cf. \cite[Ch.~1, 3.5, p.~20]{gelfand_manin})
for the double complex 
$\Ch\bul\II\bul:=\Ch\bul(\UUU;\II\bul)=\Gamma(X;\CCh\bul\II\bul)$
\begin{align}
\label{specseq1}
& \hh_{\rm I}^p \, \hh_{\rm II}^q \, \Ch\bul\II\bul \Rightarrow 
\hh^n \tot \Ch\bul\II\bul \quad \text{and} \\
\label{specseq2}
& \hh_{\rm II}^p \, \hh_{\rm I}^q \, \Ch\bul\II\bul \Rightarrow 
\hh^n \tot \Ch\bul\II\bul .
\end{align}

\paragraph{The first spectral sequence (\ref{specseq1}).}
Since $\CCh\bul\II^p$ is a flabby resolution of $\II^p$ (by Lemma
\ref{chechres}), we have
$$
\hh_{\rm II}^q\Ch\bul\II^p = \hh^q\Gamma(X;\CCh\bul\II^p) 
= \h^q(X;\II^p) 
= \begin{cases}
    \Gamma(X;\II^p) & \text{ if $q=0$,} \\
    0 & \text{ else.}
  \end{cases}
$$
Therefore
\begin{align*}
\hh_{\rm I}^p \hh_{\rm II}^q \Ch\bul\II\bul & =
  \begin{cases}
    \hh^p\Gamma(X;\II\bul) & \text{ if $q=0$,} \\
    0 & \text{ else }
  \end{cases} \\[0.5mm]
&  = 
  \begin{cases}
    \H^p(X;\FF\bul) & \text{ if $q=0$,} \\
    0 & \text{ else.}
  \end{cases}
\end{align*}
Thus the first spectral sequence degenerates and we obtain a natural
isomorphism
$$
\H\bul(X;\FF\bul) \iso \hh\bul \tot \Ch\bul\II\bul .
$$

\paragraph{The second spectral sequence
(\ref{specseq2}).}
From
\begin{align*}
\Ch^p\II^q & = 
  \Gamma\bigl(X; \bigoplus_{|I|=p+1} \iast \II^q_{|U_I} \bigr) \\    
& = \bigoplus_{|I|=p+1} \Gamma\bigl(X;  \iast \II^q_{|U_I} \bigr) \\    
& = \bigoplus_{|I|=p+1} \Gamma\bigl(U_I;  \II^q_{|U_I} \bigr) ,
\end{align*}
we see
\begin{align*}
\allowdisplaybreaks[4]
\hh_{\rm I}^q\Ch^p\II\bul & =
\bigoplus_{|I|=p+1} \hh^q\Gamma\bigl(U_I;\II\bul_{|U_I} \bigr) 
\displaybreak[0] \\    
& = \bigoplus_{|I|=p+1} \H^q\bigl(U_I; \FF\bul_{|U_I} \bigr)
\displaybreak[0] \\    
& = \bigoplus_{|I|=p+1} \Gamma\bigl(U_I;\HHH^q (\FF\bul)_{|U_I} \bigr) 
\displaybreak[0] \\
& = \bigoplus_{|I|=p+1} \Gamma\bigl(X;\iast \HHH^q (\FF\bul)_{|U_I} \bigr) 
\displaybreak[0] \\    
& = \Gamma\bigl(X;\bigoplus_{|I|=p+1} \iast \HHH^q (\FF\bul)_{|U_I} \bigr) 
\displaybreak[0] \\ 
& = \Gamma\left(X;\CCh^p\bigl(\UUU;\HHH^q (\FF\bul)\bigr)\right) 
\displaybreak[0] \\
& = \Ch^p\bigl(\UUU;\HHH^q(\FF\bul)\bigr).
\end{align*} 
This gives
$$
\hh_{\rm II}^p \, \hh_{\rm I}^q \, \Ch\bul\II\bul =
\hh^p \, \Ch\bul\left(\UUU;\HHH^q\left(\FF\bul\right)\right) 
 = \check{\h}^p\left(\UUU;\HHH^q\left(\FF\bul\right)\right)
$$
and concludes the proof.
\end{proof}

\subsubsection{A Spectral Sequence for Algebraic \deRham Cohomology of
  a Smooth Variety}

We are especially interested in Proposition \ref{specseqopcov} for
$\FF\bul = \om{X/k}$.

\begin{corollary}
\label{specseqsmooth}
Let $X$ be a smooth variety over a field $k$ and $\UUU$ an open
covering of $X$. Then we have a convergent spectral sequence for
algebraic \deRham cohomology of smooth varieties 
(as defined in Definition \ref{algdr1})
$$
E_2^{p,q} := \hh^p \bigoplus_{|I|=\,\Bul+1} \hDR{q}(U_I/k) \Rightarrow 
E_{\infty}^n := \hDR{n}(X/k) .
$$
\end{corollary}

\begin{proof}
Let us compute the initial terms of the spectral sequence
$$
\check{\h}(\UUU;\HHH^q(\FF\bul)) \Rightarrow \H^n(X;\FF\bul)
$$
of Proposition \ref{specseqopcov} for $\FF\bul = \om{X/k}$
\begin{align*}
\check{\h}^p(\UUU;\HHH^q(\om{X/k})) 
& = \hh^p\Gamma\bigl(X;\CCh\bul\bigl(\UUU;\HHH^q(\om{X/k})\bigl)\bigr) \displaybreak[0] \\
& = \hh^p\Gamma\Bigl(X;\bigoplus_{|I|=\,\Bul+1} \iast \HHH^q(\om{X/k})_{|U_I}\Bigr) \displaybreak[0] \\
& = \hh^p \bigoplus_{|I|=\,\Bul+1}
 \Gamma\bigl(X;\iast \HHH^q(\om{X/k})_{|U_I}\bigr) \displaybreak[0] \\
& = \hh^p \bigoplus_{|I|=\,\Bul+1}
 \Gamma\bigl(U_I;\HHH^q(\om{X/k})_{|U_I}\bigr) \displaybreak[0] \\
& = \hh^p \bigoplus_{|I|=\,\Bul+1}
 \H^q(U_I;\om{X/k}{}_{|U_I}) \displaybreak[0] \\
& = \hh^p \bigoplus_{|I|=\,\Bul+1}
 \H^q(U_I;\om{U_I/k}) \displaybreak[0] \\
& = \hh^p \bigoplus_{|I|=\,\Bul+1}
 \hDR{q}(U_I/k) .
\end{align*}
For the limit term we get immediately
$$
\H\bul(X;\om{X/k}) = \hDR{\Bul}(X/k)
$$
and the corollary is proved.
\end{proof}

\subsubsection{A Spectral Sequence for Algebraic \deRham Cohomology of
  a Divisor with Normal Crossings}

For a divisor $D$ with normal crossings, we have a spectral sequence
expressing $\hDR{\Bul}(D/k)$ in terms of $\hDR{\Bul}(D_{i_0 \, \cdots
  \, i_m}/k)$.

Recall that (cf. page \pageref{defwom})
$$
\wom{D/k} = \tot \Omega^{\Bul,\Bul}_{D\bul/k} = \tot \bigoplus_{|I|} \iast
\om{D_I/k} .
$$
Since we are interested in hypercohomology, we replace the $\om{\ast}$
by their Godement resolutions (cf. page \pageref{godement}), where we shall use the following
abbreviations
\label{godalg}
\begin{itemize}
\item
$ 
\GG^p_{D_I/k} := \GG^p_{\om{D_I/k}}
$,
\item
$
\GG^{p,q}_{D\bul/k} := 
  \GG^p_{\Omega^{\Bul,q}_{D/k}} \\
\phantom{\GG^{p,q}_{D\bul/k}}\, = 
  \bigoplus_{|I|=q+1} \GG^p_{\iast \om{D_I/k}} \\
\phantom{\GG^{p,q}_{D\bul/k}}\, \stackrel{(\ast)}{=} 
\bigoplus_{|I|=q+1} \iast \GG^p_{D_I/k},
$

(we have equality at $(\ast)$ since $\iast$ is exact for closed
immersions $i$)
\end{itemize}
and
\begin{itemize}
\item
$
\GG\bul_{D/k}  := 
  \GG\bul_{\wom{D/k}} \\
\phantom{\GG\bul_{D/k}}\, = 
  \GG\bul_{\tot \Omega^{\Bul,\Bul}_{D\bul/k}} \\
\phantom{\GG\bul_{D/k}} \, = 
  \tot \GG^{\Bul,\Bul}_{D\bul/k} .
$
\end{itemize}
Furthermore, we write
$$
{\rm G}^p_{D_I/k} := \Gamma\bigl(D_I; \GG^p_{D_I/k}\bigr) , \, 
{\rm G}^{p,q}_{D\bul/k}:= \Gamma\bigl(D; \GG^{p,q}_{D\bul/k}\bigr) \text{ and }
{\rm G}^p_{D/k} :=\Gamma\bigl(D; \GG^p_{D/k}\bigr)
$$
for the groups of global sections of these sheaves.

We consider the second spectral sequence for the double complex
${\rm G}^{p,q}_{D\bul/k}$
$$
\hh^p_{\rm II} \, \hh^q_{\rm I} \, {\rm G}^{\Bul,\Bul}_{D\bul/k}
\Rightarrow
\hh^n \tot {\rm G}^{\Bul,\Bul}_{D\bul/k} .
$$
For the limit terms, we obtain
$$
\hh^n \tot {\rm G}^{\Bul,\Bul}_{D\bul/k} = \hh^n {\rm G}\bul_{D/k} =
\hh^n\Gamma(D;\GG\bul_{D/k}) = \H^n(D;\wom{D/k}) = \hDR{n} (D/k) .
$$
Let us compute the initial terms
\begin{align*}
\hh^q {\rm G}^{\Bul,p} _{D\bul/k} 
& = \hh^q\Gamma(D;\bigoplus_{|I|=p+1} \iast \GG\bul_{D_I/k} ) \\
& = \hh^q \bigoplus_{|I|=p+1} \Gamma(D; \iast \GG\bul_{D_I/k} ) \\
& = \bigoplus_{|I|=p+1} \hh^q \Gamma(D; \iast \GG\bul_{D_I/k} ) \\
& = \bigoplus_{|I|=p+1} \hh^q \Gamma(D_I; \GG\bul_{D_I/k} ) \\
& = \bigoplus_{|I|=p+1} \H^q (D_I; \om{D_I/k} ) \\
& = \bigoplus_{|I|=p+1} \hDR{q} (D_I/k) .
\end{align*}

Thus we have proved the following proposition.

\begin{proposition}
\label{specseqdiv}
Let $X$ be a smooth variety over a field $k$ and $D$ a divisor with
normal crossings on $X$. Then we have a convergent spectral sequence
for algebraic \deRham cohomology of a divisor (as defined in
Definition 
\ref{algdr2})
$$
E^{p,q}_2 := \hh^p \bigoplus_{|I|=\,\Bul+1} \hDR{q}(D_I/k) \Rightarrow
E^n_{\infty} := \hDR{n}(D/k) .
$$
\end{proposition}

\section{Comparison Isomorphisms}

If we consider a smooth variety $X$ defined over \C\ and a
divisor $D$ with normal crossings on $X$, then their algebraic \deRham cohomology as
defined in Subsection \ref{algDR} turns out to be naturally isomorphic to the
singular cohomology of the associated complex analytic spaces \Xh\ and
\Dh.

Our proof proceeds in three steps: $(1)$ First we mimic the definition of
algebraic \deRham cohomology in the complex analytic setting using
holomorphic differential forms instead of algebraic ones, thus
defining analytic
\deRham cohomology groups
$$
\hDR{\Bul}(\Xh;\C), \quad \hDR{\Bul}(\Dh;\C) \quad \text{and}\quad
\hDR{\Bul}(\Xh,\Dh;\C) .
$$
These are then shown to be isomorphic to their algebraic counterparts
using an application of Serre's ``{\small GAGA}-type'' results by Grothendieck.

$(2)$ In the second step, analytic \deRham cohomology is proved to coincide
with complex cohomology, $(3)$ which in turn is isomorphic to singular
cohomology, as we show in the third step. (In the smooth case this has
already been shown in Proposition \ref{classiso}.)

\subsection{Situation}
\label{situation}

Throughout this section, $X$ will be a smooth variety over $\C$ and
$D$ a divisor with normal crossings on $X$. As usual, we denote the
complex analytic spaces associated to $X$ and $D$ by $\Xh$ and $\Dh$, 
respectively (cf. Definition \ref{asscplansp}).

\subsection{Analytic \deRham Cohomology}

If we replace in Subsection \ref{algDR} every occurrence of the complex of
algebraic differential forms $\om{Y/\C}$ on a complex variety $Y$ by the
complex of holomorphic differential forms $\om{\Yh}$ on the associated
complex analytic space \Yh, we obtain complexes
\label{defwoman}
\notation{O}{$\Omega^{\Bul,\Bul}_{\Dh{}\bul}$}{analytic double \deRham complex of $\Dh{}\bul$}
\notation{O}{$\wom{\Dh}$}{analytic \deRham complex of a divisor $\Dh$}
\notation{O}{$\wom{\Xh,\Dh}$}{analytic \deRham complex of a pair $(\Xh,\Dh)$}
$$
\om{\Xh},\quad \Omega^{\Bul,\Bul}_{\Dh{}\bul},\quad \wom{\Dh} \quad \text{and} \quad \wom{\Xh,\Dh}
$$
and thus are able to define:

\begin{definition}[Analytic \deRham cohomology]
\label{andr}
\index{DeRham cohomology@\DeRham cohomology!analytic|textbf}
\notation{H}{$\hDR{\Bul}(\Dh;\C)$}{analytic \deRham cohomology of a
    divisor $\Dh$}
  \notation{H}{$\hDR{\Bul}(\Xh,\Dh;\C)$}{analytic \deRham cohomology of a
    pair $(\Xh,\Dh)$}
  For \Xh\ and \Dh\ as in Subsection \ref{situation}, we define {\em analytic \deRham cohomology} by
  \begin{align}
   \hDR{\Bul}(\Xh;\C) & := \H\bul(\Xh;\om{\Xh}), \tag{$\ast$} \\
   \hDR{\Bul}(\Dh;\C) & := \H\bul(\Dh;\wom{\Dh}) \quad \text{and}
   \notag \\
   \hDR{\Bul}(\Xh,\Dh;\C) & := \H\bul(\Xh;\wom{\Xh,\Dh}). \notag
  \end{align}
Here {\rm ($\ast$)} is the definition of analytic \deRham cohomology  for
complex manifolds already considered on page \pageref{complcoho},
which we listed again for completeness.
\end{definition}

Note that $(\om{Y/\C} )_{\rm an} = \om{\Yh}$, hence
\begin{equation}
\begin{split}
\label{omalgan}
(\om{X/\C})_{\rm an} & = \om{\Xh}, \\
(\Omega^{\Bul,\Bul}_{D\bul/\C})_{\rm an} & = \Omega^{\Bul,\Bul}_{\Dh{}\bul}, \\
(\wom{D/\C})_{\rm an} & = \wom{\Dh} \quad \text{and} \\
(\wom{X,D/\C})_{\rm an} & = \wom{\Xh,\Dh} .
\end{split}
\end{equation}

\begin{remark}
We could easily generalize Definition \ref{andr} to arbitrary complex
ma\-ni\-folds (and normal-crossings-divisors on them) not necessarily
associated to something algebraic, but we
will not need this.
\end{remark}

\subsection{The Long Exact Sequence in Analytic \deRham Cohomology}

The results of Subsection \ref{longseqalgDR} carry over one-to-one to the
complex analytic case.

\begin{proposition}
\label{longexseqandr}
We have naturally a long exact sequence in analytic \deRham cohomology
(as defined in Definition \ref{andr})
\begin{align*}
\cdots & \to \hDR{p-1}(\Dh;\C) \to \\
\hDR{p}(\Xh,\Dh;\C) \to \, \hDR{p}(\Xh;\C) & \to \hDR{p}(\Dh;\C) \to \\
\hDR{p+1}(\Xh,\Dh;\C) \to \, \cdots \, , \hspace{1.5cm} &
\end{align*}
where \Xh\ and \Dh\ are 
as in Subsection \ref{situation}.
\end{proposition}

Moreover, by Equation (\ref{gagahypernatmorph}) 
on page \pageref{gagahypernatmorph}
and Equation (\ref{omalgan}), 
we have a natural 
morphism from the
long exact sequence in algebraic \deRham cohomology to the one in
analytic \deRham cohomology
\begin{equation}
\label{alganlongseq}
\begin{CD}
\cdots \to \, @. \hDR{p}(X,D/\C) @>>> \hDR{p}(X/\C) @>>> \hDR{p}(D/\C)
@. \, \to \cdots \phantom{\, .} \\
@. @VVV @VVV @VVV @. \\ 
\cdots \to \, @. \hDR{p}(\Xh,\Dh;\C) @>>> \hDR{p}(\Xh;\C) @>>>
\hDR{p}(\Dh;\C) @. \, \to \cdots \, .
\end{CD}
\end{equation}
We want to show that all the vertical maps are isomorphisms.

\subsection{Some Spectral Sequences for Analytic \deRham
  Cohomology}

The proof of Proposition \ref{specseqsmooth} applies in the analytic case as well
and we obtain 

\begin{proposition}
\label{specseqsmoothan}
We have a convergent spectral sequence for analytic \deRham cohomology
(as defined in Definition \ref{andr})
$$
E^{p,q}_2:= \hh^p \bigoplus_{|I|=\,\Bul+1} \hDR{q}(\Uh_I;\C) \Rightarrow 
E^n_{\infty} := \hDR{n} (\Xh;\C),
$$
where $\Xh$ is associated to a smooth variety $X$ over \C\ and
$\UUU=\left( U_j \right)_{j\in J}$
is an open covering of $X$.
\end{proposition}

Furthermore, the natural map of hypercohomology groups
(cf. (\ref{gagahypernatmorph}) and (\ref{omalgan}))
$$
\H\bul(U;\om{U/\C}) \to \H\bul(\Uh;\om{\Uh})
\quad
\text{for\quad$U\subseteq X$ open},
$$
provides us with a natural map of spectral
sequences
\begin{equation*}
\label{alganspecseq}
\begin{split}
\check{\h}^p(\UUU;\HHH^q(\om{X/\C})) \Rightarrow 
\H^n(X;\om{X/\C}) \hspace*{0.5cm} \\
\downarrow \hspace{3cm} \downarrow \hspace*{1.4cm} \\
\check{\h}^p(\UUU^{\rm an};\HHH^q(\om{\Xh})) \Rightarrow
\H^n(\Xh;\om{\Xh}),
\end{split}
\end{equation*}
which we can rewrite using Corollary \ref{specseqsmooth} and
Proposition \ref{specseqsmoothan}
\begin{equation}
\label{specseqsmoothmorph}
\begin{split}
\hh^p \bigoplus_{|I|=\,\Bul+1} \hDR{q}(U_I/\C) \Rightarrow 
\hDR{n} (X/\C) \hspace*{0.5cm} \\[-3mm]
\down[] \hspace{3cm}  \down[] \hspace*{1cm} \\
\hh^p \bigoplus_{|I|=\,\Bul+1} \hDR{q}(\Uh_I;\C) \Rightarrow 
\hDR{n} (\Xh;\C).
\end{split}
\end{equation}

Likewise, we can transfer Proposition \ref{specseqdiv} to the analytic
case.

\begin{proposition}
\label{specseqdivan}
We have a convergent spectral sequence
$$
E^{p,q}_2 := \hh^p \bigoplus_{|I|=\,\Bul+1} \hDR{q}(\Dh_I;\C) \Rightarrow 
E^n_{\infty} := \hDR{n}(\Dh;\C) ,
$$
where $\Dh$ is associated to a normal-crossings-divisor $D$ on a
smooth variety $X$ defined over \C\ as in Subsection \ref{situation}.
\end{proposition}

Again, we have a natural map of spectral sequences. We write (cf. page
\pageref{godalg})
$$
\GG^{p,q}_{\Dh{}\bul} := \left( \GG^{p,q}_{D\bul/\C}\right)_{\rm an}
$$
and
$${\rm G}^{p,q}_{\Dh{}\bul} := \Gamma(\Dh;\GG^{p,q}_{\Dh{}\bul}) .
$$
Now the natural map of zeroth homology (cf. Equation
(\ref{gaganatmorph}) for $i=0$)
$$
{\rm G}^{\Bul,\Bul}_{D\bul/\C} \to {\rm G}^{\Bul,\Bul}_{\Dh{}\bul}
$$
gives us a natural map of spectral sequences
\begin{equation*}
\begin{split}
\hh^p_{\rm I} \, \hh^q_{\rm II} \, {\rm G}^{\Bul,\Bul}_{D\bul/\C} 
\Rightarrow
\hh^n \tot {\rm G}^{\Bul,\Bul}_{D\bul/\C} \hspace*{0.1cm} \\
\downarrow \hspace{2.5cm} \downarrow \hspace*{1.2cm} \\
\hh^p_{\rm I} \, \hh^q_{\rm II} \, {\rm G}^{\Bul,\Bul}_{\Dh{}\bul} 
\Rightarrow
\hh^n \tot {\rm G}^{\Bul,\Bul}_{\Dh{}\bul} .
\end{split}
\end{equation*}
With propositions \ref{specseqdiv} and \ref{specseqdivan}, this map
takes the form
\begin{equation}
\begin{split}
\label{specseqdivmorph}
\hh^p \bigoplus_{|I|=\,\Bul+1} \hDR{q}(D_I/\C) \Rightarrow 
\hDR{n}(D/\C) \hspace*{0.5cm} \\[-3mm]
\down[] \hspace{2.5cm} \down[] \hspace*{1.2cm} \\
\hh^p \bigoplus_{|I|=\,\Bul+1} \hDR{q}(\Dh_I;\C) \Rightarrow 
\hDR{n}(\Dh;\C) .
\end{split}
\end{equation}

\subsection{Comparison of Algebraic and Analytic \deRham
  Cohomology}

The proof of $\hDR{\Bul}(X/\C) \iso \hDR{\Bul}(\Xh;\C)$,
i.e. $\H\bul(X;\om{X/\C}) \iso \H\bul(\Xh;\om{\Xh})$ is easy in the
{\bf projective case}, because of (\ref{omalgan}) and Corollary \ref{gagahyper}.

The {\bf affine case} is contained in a theorem of Grothendieck.

\begin{theorem}[Grothendieck, {\cite[Thm. 1, p. 95]{grothendieck_66}}]
\label{grothendieck}
Let $X$ be a smooth affine variety over \C. Then the natural map
between algebraic and analytic \deRham cohomology 
(as defined in definitions \ref{algdr1}  and  \ref{andr})
$$
\hDR{\Bul}(X/\C) \stackrel{\sim}{\longrightarrow} \hDR{\Bul}(\Xh;\C)
$$ 
is an isomorphism.
\end{theorem}

We give only some comments on the proof. 
Let $j:X \inclusion \ol{X}$ be the projective
closure of $X$ and let $Z:=\ol{X} \setminus X$. Furthermore, we set
\begin{align*}
\om{\ol{X}/\C}(\ast Z) & := \underset{n}{\injlim} \, \om{\ol{X}/\C}(n Z)
\quad \text{and} \\
\om{\ol{X}^{\rm an}}(\ast \Zh) & := \underset{n}{\injlim} \,
\om{\ol{X}^{\rm an}}(n \Zh),
\end{align*}
where $\FF(nZ')$ denotes the sheaf $\FF$ twisted by the $n$-fold of
the divisor $Z'$.
In his paper \cite{grothendieck_66}, Grothendieck 
considers the following chain of isomorphisms:
\begin{align*}
\allowdisplaybreaks[4]
\hDR{\Bul}(X/\C) = \, & \H\bul(X;\om{X/\C}) 
\displaybreak[0] \\[-1mm]
  & \hspace{1cm} \up \quad \text{$\Omega^p_{X/\C}$ are
    $\Gamma(X;?)$-acyclic \cite[Thm. III.3.5, p. 215]{hartshorne}} 
\displaybreak[0] \\
& \hh\bul\Gamma(X;\om{X/\C}) 
\displaybreak[0] \\[-1mm]
  & \hspace{1cm} \Arrowvert 
\displaybreak[0] \\
& \hh\bul\Gamma(\ol{X};j_{\ast} \om{X/\C}) 
\displaybreak[0] \\[-1mm]
  & \hspace{1cm} \Arrowvert \\
& \hh\bul\Gamma(\ol{X};\om{\ol{X}/\C}(\ast Z)) 
\displaybreak[0] \\[-1mm]
  & \hspace{1cm} \down \quad \text{\cite[Prop. III.2.9, p. 209]{hartshorne}}
\displaybreak[0] \\
& \underset{n}{\injlim} \, \hh\bul\Gamma(\ol{X};\om{\ol{X}/\C}(nZ)) 
\displaybreak[0] \\[-3mm]
  & \hspace{1cm} \down \quad \text{Thm. \ref{gagathm} ({\small
      GAGA})} \displaybreak[0] \\
& \underset{n}{\injlim} \, 
  \hh\bul\Gamma(
    \ol{X}^{\rm an};
    \om{\ol{X}^{\rm an}}(n \Zh)) 
\displaybreak[0] \\[-3mm]
  & \hspace{1cm} \down \quad 
  \text{\cite[Lemma 6]{atiyah_hodge}} 
\displaybreak[0] \\
& \hh\bul\Gamma(\ol{X}^{\rm an};\om{\ol{X}^{\rm an}}(\ast \Zh)) 
\displaybreak[0] \\[-1mm]
  & \hspace{1cm} \down \quad \text{\bf main step} 
\displaybreak[0] \\
& \hh\bul\Gamma(\ol{X}^{\rm an};j_{\ast}\om{\Xh}) 
\displaybreak[0] \\[-1mm]
  & \hspace{1cm} \Arrowvert 
\displaybreak[0] \\
& \hh\bul\Gamma(\Xh;\om{\Xh}) 
\displaybreak[0] \\[0mm]
  & \hspace{1cm} \down \quad 
      \text{\parbox{10cm}{$\Omega^p_{\Xh}$ are
      $\Gamma(\Xh;?)$-acyclic 
      since $\Xh$ is Stein \\
      \cite[IV.1.1, Def.~1, p.~103;
      V.1, Bem., p.~130]{grauert_remmert}}} 
\displaybreak[0] \\
\hDR{\Bul}(\Xh;\C) = \, & \H\bul(\Xh;\om{\Xh}) .
\end{align*}
In the case that $Z$ has normal crossings, the main step is due to
Atiyah and Hodge \cite[Lemma 17]{atiyah_hodge} and is proved by
explicit calculations. 
Grothendieck reduces the general case to this one by using the resolution
of singularities according to Hironaka \cite{hironaka_64}.

\paragraph{The smooth case.}
Choose an open affine covering $\UUU=\left( U_j \right)_{j\in J}$ of $X$. 
Since $X$ is separated, all the $U_I$ are affine. Hence by Theorem
\ref{grothendieck}  we get an
isomorphism on the initial terms in (\ref{specseqsmoothmorph}), and
therefore also on the limit terms
$$
\hDR{\Bul}(X/\C) \iso \hDR{\Bul}(\Xh;\C) .
$$

\paragraph{The case of a normal-crossings-divisor.}
Since $D$ is a divisor with normal crossings, all the $D_I$ are smooth.
Thus we see from the discussion of the smooth case, that we have
an isomorphism on the initial terms in (\ref{specseqdivmorph}), hence also on
the limit terms
\begin{align*}
\hDR{\Bul}(D/\C) \iso
\hDR{\Bul}(\Dh;\C) .
\end{align*}

\paragraph{The relative case.}
By the $5$-lemma, this case follows immediately from the two cases
considered above using diagram (\ref{alganlongseq}).

\medskip
We summarize our results so far.

\begin{proposition}
\label{compalgan}
For $X$ a smooth variety over \C\ and $D$ a normal-crossings-divisor on
$X$, we have a natural isomorphism between algebraic \deRham cohomology
of $X$ and $D$ (see definitions \ref{algdr1}, \ref{algdr2}, \ref{algdr3}) 
and analytic \deRham cohomology (see Definition \ref{andr}) of the associated complex analytic
spaces $\Xh$ and $\Dh$
\begin{equation*}
\begin{CD}
\cdots \to \, @. \hDR{p}(X,D/\C) @>>> \hDR{p}(X/\C) @>>> \hDR{p}(D/\C)
@. \, \to \cdots \phantom{\, .} \\
@. @VV{\wr}V @VV{\wr}V @VV{\wr}V @. \\ 
\cdots \to \, @. \hDR{p}(\Xh,\Dh;\C) @>>> \hDR{p}(\Xh;\C) @>>>
\hDR{p}(\Dh;\C) @. \, \to \cdots \, .
\end{CD}
\end{equation*}
\end{proposition}

\subsection{Complex Cohomology}
\label{red0}

\begin{definition}[Complex cohomology]
\label{cplxcoho}
\index{Complex cohomology|textbf}
Let $\Xh$ and $\Dh$ be as in Subsection \ref{situation}.
In the absolute case, {\em complex cohomology} is simply sheaf cohomology of
the constant sheaf with fibre \C
$$
\h\bul(\Xh;\C) := \h\bul(\Xh;\C_{\Xh}) \quad \text{and} \quad
\h\bul(\Dh;\C) := \h\bul(\Dh;\C_{\Dh}) .
$$
We define a relative version by
$$
\h\bul(\Xh,\Dh;\C) := \h\bul(\Xh;j_! \C_{\Uh}),
$$
where $U:=X \setminus D$ and $j: U \inclusion X$.
\notation{H}{$\h\bul(\Xh,\Dh;\C)$}{complex cohomology of a pair $(\Xh,\Dh)$}
\end{definition}

\begin{remark}
Taking $\Q$ instead of \C\ in the above definition would give us {\em rational
cohomology groups} $\h\bul(\Xh;\Q)$, $\h\bul(\Dh;\Q)$ and $\h\bul(\Xh,\Dh;\Q)$. 
\end{remark}

The short exact sequence
\begin{equation}
\label{complcohoshortseq}
0 \to j_! \C_{\Uh} \stackrel{\alpha}{\longto} \C_{\Xh} \to \iast \C_{\Dh} \to 0
\end{equation}
with $j:\Uh\inclusion \Xh$ and $i:\Dh\inclusion \Xh$,
yields a long exact sequence in complex cohomology
\begin{align*}
\cdots & \to \h^{p-1}(\Dh;\C) \to \\
\h^{p}(\Xh,\Dh;\C) \to \, \h^{p}(\Xh;\C) & \to \h^{p}(\Dh;\C) \to \\
\h^{p+1}(\Xh,\Dh;\C) \to \, \cdots \, . \hspace{0.85cm} &
\end{align*}

The mapping cone $M_{\C}$ (cf. the appendix) of the natural restriction map 
$\C_{\Xh} \to \iast \C_{\Dh}$
\notation{M}{$M_{\C}$}{mapping cone of $\C_{\Xh} \to \iast \C_{\Dh}$}
gives us another short exact sequence
\begin{equation}
\label{cplxmapcon}
0 \to \iast \C_{\Dh}[-1] \to M_{\C} \to \C_{\Xh}[0] \to 0.
\end{equation}
Consider the map
$$
\gamma:=(0,\alpha): j_! \C_{\Uh} [0] \to 
M_{\C} = \iast \C_{\Dh}[-1] \oplus \C_{\Xh}[0],
$$
where 
$\alpha: j_! \C_{\Uh} [0] \to \C_{\Xh} [0]$
was defined in (\ref{complcohoshortseq}) and
$0:j_! \C_{\Uh} \to \iast\C_{\Dh}[-1]$ 
denotes the zero map.
In order to show that $\alpha$ is a morphism of complexes, 
we only have to check the compatibility of the differentials, but this
follows from the fact that (\ref{complcohoshortseq}) is a complex. The
exactness of the sequence in (\ref{complcohoshortseq}) translates into
$\gamma$ being a quasi-isomorphism.
This gives us the following diagram
$$
\begin{matrix}
\h^{p-1}(\Dh;\C) 
\hspace{-1.7mm}&\hspace{-1.7mm} \to \hspace{-1.7mm}&\hspace{-1.7mm} 
\h^p(\Xh,\Dh;\C)
\hspace{-1.7mm}&\hspace{-1.7mm} \to \hspace{-1.7mm}&\hspace{-1.7mm} 
\h^p(\Xh;\C) 
\hspace{-1.7mm}&\hspace{-1.7mm} \to \hspace{-1.7mm}&\hspace{-1.7mm} 
\h^p(\Dh;\C) \\
\Arrowvert 
\hspace{-1.7mm}&\hspace{-1.7mm} \hspace{-1.7mm}&\hspace{-1.7mm}  
\Arrowvert 
\hspace{-1.7mm}&\hspace{-1.7mm} \hspace{-1.7mm}&\hspace{-1.7mm}
\Arrowvert 
\hspace{-1.7mm}&\hspace{-1.7mm} \hspace{-1.7mm}&\hspace{-1.7mm} 
\Arrowvert \\
\H^p(\Xh;\iast \C_{\Dh}[-1]) 
\hspace{-1.7mm}&\hspace{-1.7mm} \to \hspace{-1.7mm}&\hspace{-1.7mm} 
\H^p(\Xh;M_{\C}) 
\hspace{-1.7mm}&\hspace{-1.7mm} \to \hspace{-1.7mm}&\hspace{-1.7mm} 
\H^p(\Xh;\C_{\Xh}[0]) 
\hspace{-1.7mm}&\hspace{-1.7mm} \to \hspace{-1.7mm}&\hspace{-1.7mm} 
\H^{p+1}(\Xh;\iast\C_{\Dh}[-1]).
\end{matrix}
$$
Unfortunately, this diagram does not commute (there is a sign
mismatch), so we consider the natural transformation
$$
\eps^p = (-1)^p \, {\rm id} \quad \text{for} \quad p\in\Z
$$
and write down a new diagram
$$
\begin{matrix}
\h^{p-1}(\Dh;\C) 
\hspace{-1.7mm}&\hspace{-1.7mm} \to \hspace{-1.7mm}&\hspace{-1.7mm} 
\h^p(\Xh,\Dh;\C)
\hspace{-1.7mm}&\hspace{-1.7mm} \to \hspace{-1.7mm}&\hspace{-1.7mm} 
\h^p(\Xh;\C) 
\hspace{-1.7mm}&\hspace{-1.7mm} \to \hspace{-1.7mm}&\hspace{-1.7mm} 
\h^p(\Dh;\C) \\[1mm]
{\scriptstyle \eps^{p-1}} \downiso 
\hspace{-1.7mm}&\hspace{-1.7mm} \hspace{-1.7mm}&\hspace{-1.7mm}
{\scriptstyle \eps^p} \downiso 
\hspace{-1.7mm}&\hspace{-1.7mm} \hspace{-1.7mm}&\hspace{-1.7mm}
{\scriptstyle \eps^p} \downiso 
\hspace{-1.7mm}&\hspace{-1.7mm} \hspace{-1.7mm}&\hspace{-1.7mm} 
{\scriptstyle \eps^p} \downiso \\[1mm]
\H^p(\Xh;\iast \C_{\Dh}[-1]) 
\hspace{-1.7mm}&\hspace{-1.7mm} \to \hspace{-1.7mm}&\hspace{-1.7mm} 
\H^p(\Xh;M_{\C}) 
\hspace{-1.7mm}&\hspace{-1.7mm} \to \hspace{-1.7mm}&\hspace{-1.7mm} 
\H^p(\Xh;\C_{\Xh}[0]) 
\hspace{-1.7mm}&\hspace{-1.7mm} \to \hspace{-1.7mm}&\hspace{-1.7mm} 
\H^{p+1}(\Xh;\iast
\C_{\Dh}[-1]) \, .
\end{matrix}
$$
This diagram commutes as we see by replacing the sheaves involved in
(\ref{complcohoshortseq}) by injective resolutions $\AA\bul$, $\BB\bul$ and $\CC\bul$
and applying Lemma \ref{homalg} below to the complexes $\a\bul$, $\b\bul$ and
$\c\bul$ of their global sections.

\begin{lemma}
\label{homalg}
Let 
$$
0\longto \a\bul \stackrel{\alpha}{\longto} \b\bul
  \stackrel{\beta}{\longto} \c\bul \longto 0 
$$
be a short exact sequence of complexes of abelian groups and
$$
0 \longto \c\bul[-1] \stackrel{\iota}{\longto} \m\bul
\stackrel{\pi}{\longto} \b\bul \longto 0
$$
the mapping cone of $\beta: \b\bul \to \c\bul$. 
Then we have a quasi-isomorphism
$$
\gamma:=(0,\alpha): \a\bul \stackrel{\sim}{\longto}
\m\bul=\c\bul[-1]\oplus\b\bul
$$
and a commutative diagram
$$
\begin{CD}
\hh^{p-1} \c\bul & @>{\delta}>> & \hh^p \a\bul &
@>{\alpha_{\ast}}>> & \hh^p\b\bul & @>{\beta_{\ast}}>> &
\hh^p \c\bul \\
@V{\eps^{p-1}}V{\wr}V & \hspace{-17mm} \fbox{$\scriptstyle 1$} &
  @V{\eps^p\circ\gamma_{\ast}}V{\wr}V & \hspace{-15mm} \circlearrowleft & 
@V{\eps^p}V{\wr}V & \hspace{-15mm} \fbox{$\scriptstyle 2$} & @V{\eps^p}V{\wr}V \\
\hh^{p-1} \c\bul & @>{\iota_{\ast}}>> & \hh^p \m\bul &
@>{\pi_{\ast}}>> & \hh^p\b\bul & @>{\delta'}>> &
\hh^p \c\bul \, .
\end{CD}
$$
\end{lemma}

\begin{proof}
The assertion about the quasi-isomorphism $\gamma: \a\bul \to \m\bul$
follows from the exactness of $0\to\a\bul\to\b\bul\to\c\bul\to 0$. The
middle square commutes because of $\pi\circ\gamma = \alpha$.
So we are left to show the commutativity of the squares
\fbox{$\scriptstyle 1$} and \fbox{$\scriptstyle 2$}.
\begin{description}
\item[Ad \fbox{$\scriptstyle 2$}\hspace{0.3mm}.)] 
From the diagram
$$
\begin{matrix}
\c^p & \stackrel{\iota}{\longto} & \m^{p+1} & &  \\
& & \uparrow {\scriptstyle d_{\m}} & & \\
& & \m^p & \stackrel{\pi}{\longto} & \b^p 
\end{matrix}
$$
and the definition of the connecting morphism $\delta'$ we see for $b\in\b^p$ 
$$
\delta'[b] = [ \beta(b) ] = \beta_{\ast}[b],
$$
hence
$$
(\delta'\circ\eps^p)[b] = (\eps^p\circ\beta_{\ast})[b].
$$
\item[Ad \fbox{$\scriptstyle 1$}\hspace{0.3mm}.)]
Pick an element $c\in\c^{p-1}$ and choose a preimage $b\in\b^{p-1}$. With
$$
\begin{matrix}
\a^p & \stackrel{\alpha}{\longto} & \b^p & & \\
& & \uparrow {\scriptstyle d_{\b}} & & \\
& & \b^{p-1} & \stackrel{\beta}{\longto} & \c^{p-1} 
\end{matrix}
$$
we get
$$
\delta[c] = [\alpha^{-1}\circ d_{\b} b  ] \in \hh^p \a\bul .
$$
Applying $\gamma_{\ast}$ gives
$$
(\gamma_{\ast}\circ\delta)[c] = 
(0,\alpha)_{\ast}[\alpha^{-1}\circ d_{\b} b] = 
[0,d_{\b} b] \in \hh^p \m \bul = \hh^p (\c\bul[-1] \oplus \b\bul)
$$
or after subtracting the coboundary $d_{\m}(0,b)=(\beta(b),d_{\b} b)
= (c,d_{\b} b) $
$$
(\gamma_{\ast}\circ\delta)[c] = [-c,0] = - \iota_{\ast} [c],
$$
hence
$$(\eps^p\circ\gamma_{\ast}\circ\delta)[c]=(\iota_{\ast}\circ\eps^{p-1})[c].
$$
\end{description}
\end{proof}

\subsection{Comparison of Analytic {\deRham Cohomology} and Complex Cohomology}

We apply the \v{C}ech functor (cf. Subsection \ref{chechcoho}) for the closed covering 
$\DDDh := (\Dh_i)_{i=1}^r$ 
\notation{D}{$\DDDh$}{closed covering $(\Dh_i)_{i=1}^r$}
to the constant sheaf $\C_{\Dh}$, and thus get
$ \CCh^q(\DDDh;\C):=\CCh^q(\DDDh;\C_{\Dh})$. 
Now the natural inclusions $\C_{\Dh_I} \inclusion \Omega^0_{\Dh_I}$
sum up to a natural morphism of sheaves on $\Dh$
\begin{equation}
\label{morcomplan}
\begin{split}
\CCh^q(\DDDh;\C) & = \bigoplus_{|I|=q+1} \iast \C_{\Dh_I} \\
& \hspace{1cm} \downarrow \\
\Omega^{0,q}_{\Dh{}\bul} & = \bigoplus_{|I|=q+1} \iast \Omega^0_{\Dh_I} ,
\end{split}
\end{equation}
where $i:\Dh_I\inclusion\Dh$.
Since the complexes $\om{\Dh_I}$ are resolutions of the sheaves $\C_{\Dh_I}$ 
and $\iast$ and $\oplus$ are exact
functors,
this shows that $\Omega^{\Bul,q}_{\Dh{}\bul}$ is a resolution of 
$\CCh^q(\DDDh;\C)$.

If we consider $\CCh\bul(\DDDh;\C)[0]$ as a double complex
(cf. the appendix)
concentrated in the zeroth row, we can rewrite (\ref{morcomplan})
as
$$
\CCh\bul(\DDDh;\C)[0] \to \Omega^{\Bul,\Bul}_{\Dh{}\bul}
$$
or, after taking the total complex, as 
$$
\CCh\bul(\DDDh;\C) \to \wom{\Dh} .
$$
This morphism is a quasi-isomorphism by Lemma \ref{lemma3}
(and Remark \ref{lemma3dual}) and fits into the following commutative diagram
\begin{align*}
\begin{matrix}
\C_{\Xh}[0] & \to & \iast \C_{\Dh}[0] \\
\Arrowvert & & \downarrow{\scriptstyle \alpha} \\
\C_{\Xh}[0] & \to & \iast \CCh\bul(\DDDh;\C) \\
\downiso & & \downiso \\
\om{\Xh} & \to & \iast \wom{\Dh} \, .
\end{matrix}
\end{align*}
The map $\alpha$ is a quasi-isomorphism as a consequence of the
following proposition.
\begin{proposition}
\label{closedres}
We have an exact sequence 
$$
0 \to \C_{\Dh} \to \CCh^0(\DDDh;\C) \to \CCh^1(\DDDh;\C) \to \, \cdots\,.
$$
\end{proposition}

\begin{proof}[Proof (cf. {\cite[Lemma III.4.2,
  p. 220]{hartshorne}}).]
This can be checked on the stalks and is therefore a purely
combinatorial problem. 
Let $x\in\Dh$ and assume w.l.o.g. that the irreducible components
$\Dh_1,\ldots,\Dh_s$ but no other component of $\Dh$ are passing through $x$.
For the stalks at $x$, we write
$$
\Ch^p := \CCh^p(\DDDh;\C)_x = 
\bigoplus_{1\le i_0 < \cdots < i_p \le s} 
 \C_{\Dh_{i_0 \, \cdots \, i_p,x}}.
$$
Observe that exactness at the zeroth and first
step is obvious. 
Using the
``coordinate functions'' for a set of indices $j_0,\ldots,j_p$
\begin{align*}
\Ch^p = 
\bigoplus_{1\le i_0 < \cdots < i_p \le s} \C_{\Dh_{i_0 \, \cdots \, i_p,x}}
& \surjection \C_{\Dh_{j_0 \, \cdots \, j_p,x}} \\
\oplus\,\alpha_{i_0\,\cdots\,i_p} & \mapsto \alpha_{j_0 \, \cdots \, j_p}
\end{align*}
we define a homotopy operator
$$
k^p : \Ch^p \to \Ch^{p-1}
$$
such that for any $\alpha\in\Ch^p$ 
the ``coordinates'' of its image $k^p(\alpha)$ satisfy
$$
k^p(\alpha)_{j_0 \, \cdots \, j_{p-1}} =
\begin{cases}
0 & \text{if $j_0=1$} \\
\alpha_{1\,j_0 \, \cdots \, j_{p-1}} & \text{else} .
\end{cases} 
$$
Now an elementary calculation shows that for $p\ge 1$
$$
d^{p-1}\circ k^p + k^{p+1} \circ d^p = {\rm id}_{\Ch^p}.
$$
Hence $k$ is a homotopy between the identity map and the zero map on
$\Ch\bul$, and we conclude that
$$
\hh^p\,\Ch\bul = 0 \quad \text{for} \quad p\ge 1 .
$$
\end{proof}

Thus we have a commutative diagram
$$
\begin{CD}
\C_{\Xh}[0] @>>> \iast \C_{\Dh}[0] \\
@VV{\wr}V @VV{\wr}V \\
\om{\Xh} @>>> \iast \wom{\Dh},
\end{CD}
$$
where the vertical maps are quasi-isomorphisms. As a consequence, we
also have a quasi-isomorphism between the respective mapping cones
\begin{equation}
\label{cplxwom}
\begin{CD}
0 @>>> \iast \C_{\Dh}[-1] @>>> M_{\C} @>>> \C_{\Xh}[0] @>>> 0 \\
@. @VV{\wr}V @VV{\wr}V @VV{\wr}V @. \\
0 @>>> \iast \wom{\Dh}[-1] @>>> \wom{\Xh,\Dh} @>>> \om{\Xh} @>>> 0 \, .
\end{CD}
\end{equation}
Taking hypercohomology proves the following proposition.

\begin{proposition}
\label{compancplx}
We have a natural isomorphism between the long exact sequence in
analytic \deRham cohomology (cf. Definition \ref{andr}) 
and the sequence in complex cohomology (cf. Definition \ref{cplxcoho})
\begin{equation*}
\begin{CD}
\cdots \to \, @. \hDR{p-1}(\Xh;\C) @>>>
\hDR{p-1}(\Dh;\C) @>>> \hDR{p}(\Xh,\Dh;\C) @. \, \to \cdots \\
@. @AA{\wr}A @AA{\wr}A @AA{\wr}A @. \\ 
\cdots \to \, @. \h^{p-1}(\Xh;\C) @>>> \h^{p-1}(\Dh;\C) @>>>
\h^p(\Xh,\Dh;\C) @. \, \to \cdots \, ,\\
\end{CD}
\end{equation*}
where $\Xh$ and $\Dh$ are as usual (see Subsection \ref{situation}).
\end{proposition}

\subsection{Singular Cohomology}
\index{Singular cohomology|textbf}

We extend the definition of singular cohomology given on page \pageref{singcoho}
\notation{H}{$\h\bul\sing(M;\C)$}{singular cohomology of a space $M$}
\begin{align*}
& \h\bul\sing(\Xh;\C) :=
  \hh\bul\Gamma\bigl(\Xh;\CC\bul\sing(\Xh;\C)\bigr) 
  \quad\text{and} \\
& \h\bul\sing(\Dh;\C) := \hh\bul\Gamma\bigl(\Dh;\CC\bul\sing(\Dh;\C)\bigr)
\end{align*}
by a relative version.

\begin{definition}[Relative singular cohomology]
\label{relsingcoho}
As usual, let $\Xh$ and $\Dh$ be as in Subsection \ref{situation}.
  We set
  $$
  \CC\bul\sing(\Xh,\Dh;\C):= \ker(\CC\bul\sing(\Xh;\C) \to
  \iast \CC\bul\sing(\Dh;\C))
  $$
\notation{C}{$\CC^p\sing(\Xh,\Dh;\C)$}{sheaf of singular $p$-cochains
  of a pair $(\Xh,\Dh)$}
  and define
  $$
  \h\bul\sing(\Xh,\Dh;\C) := \hh\bul\Gamma(\Xh;\CC\bul\sing(\Xh,\Dh;\C)).
  $$
\notation{H}{$\h\bul\sing(\Xh,\Dh;\C)$}{singular cohomology of a pair $(\Xh,\Dh)$}
\end{definition}

\begin{remark}
Similarly, one defines singular cohomology $\h\bul\sing(\Xh,\Dh;k)$
with coefficients in a field $k$ different from \C. We will need
$\h\bul\sing(\Xh,\Dh;\Q)$ occasionly. Since $?\tensor{\Q}\C$ is an
exact functor, we have
$$
\h\bul\sing(\Xh,\Dh;\Q) \tensor{\Q} \C = \h\bul\sing(\Xh,\Dh;\C).
$$
\end{remark}

\vspace{3mm}
The sheaf $\CC\bul\sing(\Xh,\Dh;\C)$   
is $\Gamma(\Xh;?)$-acyclic because of the surjectivity of
$$
{\rm C}\bul\sing(\Xh;\C) \surjection {\rm C}\bul\sing(\Dh;\C) = 
\Gamma(\Xh;\iast\CC\bul\sing(\Dh;\C)) .
$$
Thus the short exact sequence
$$
0 \to \CC\bul\sing(\Xh,\Dh;\C) \to
\CC\bul\sing(\Xh;\C) \to
\iast\CC\bul\sing(\Dh;\C) \to 0
$$
gives rise to a long exact sequence in singular cohomology
\begin{align*}
\cdots & \to \hspace{1mm} \h^{p-1}\sing(\Dh;\C) \to \\
\h^{p}\sing(\Xh,\Dh;\C) \to \, \h^{p}\sing(\Xh;\C) & \to \
\h^{p}\sing(\Dh;\C) \to \\
\h^{p+1}\sing(\Xh,\Dh;\C) \to \, \cdots \, . \hspace{1.6cm} &
\end{align*}

\subsection{Comparison of Complex Cohomology and Singular Cohomology}
\label{compcplxsing}

We have a commutative diagram
\begin{equation}
\label{cplxsing}
\begin{CD}
0 @>>> j_! \C_{\Uh}[0] @>>> \C_{\Xh}[0] @>>> \iast \C_{\Dh}[0] @>>> 0 \\
@. @VV{\alpha}V @VV{\beta}V @VV{\gamma}V @. \\
0 @>>> \CC\bul\sing(\Xh,\Dh;\C) @>>> \CC\bul\sing(\Xh;\C) @>>>
\iast\CC\bul\sing(\Dh;\C) @>>> 0 ,
\end{CD}
\end{equation}
where $\beta$ and $\gamma$ are quasi-isomorphisms, since the singular
cochain complex is a resolution of the constant sheaf (see Proposition \ref{singres}) and
$\iast$ is exact for closed immersions. Therefore $\alpha$ is also a
quasi-isomorphism.

Taking hypercohomology gives a natural isomorphism between the long
exact sequence in complex cohomology and the sequence in singular
cohomology. Thus we have finally proved the following theorem.

\subsection{Comparison Theorem}
\label{red1}

\begin{theorem}[Comparison theorem]
\label{comp}
Let $X$ be a smooth variety defined over\/ \C\ and $D$ a divisor
on $X$ with normal crossings 
(cf. Definition \ref{normalcrossings}). 
The associated complex analytic spaces 
(cf. Subsection \ref{defcomplexanspace})
are denoted by $\Xh$ and $\Dh$, respectively.

Then we have natural isomorphisms between
the long exact sequences in
\begin{itemize}
\item algebraic \deRham cohomology (cf. definitions \ref{algdr1}, \ref{algdr2}, \ref{algdr3}),
\item analytic \deRham cohomology (cf. Definition \ref{andr}),
\item complex cohomology (cf. Definition \ref{cplxcoho}), and
\item singular cohomology (cf. Definition \ref{relsingcoho})
\end{itemize}
as shown in the commutative diagram below
$$
\begin{CD}
\cdots \, \longto \, @.
\hDR{p}(X,D/\C) @>>> 
\underset{\makebox[0cm]{\text{\rm\footnotesize algebraic \deRham cohomology}}}
{\hDR{p}(X/\C)} @>>> 
\hDR{p}(D/\C) @. \, \longto \, \cdots \\
@. @VV{\wr}V @VV{\wr}V @VV{\wr}V @. \\
\cdots \, \longto \, @.
\hDR{p}(\Xh,\Dh;\C) @>>>
\underset{\makebox[0cm]{\text{\rm\footnotesize analytic \deRham cohomology}}} 
{\hDR{p}(\Xh;\C)} @>>>
\hDR{p}(\Dh;\C) @. \, \longto \, \cdots \\
@. @AA{\wr}A @AA{\wr}A @AA{\wr}A @. \\
\cdots \, \longto \, @.
\h^{p}(\Xh,\Dh;\C) @>>>
\underset{\makebox[0cm]{\text{\rm\footnotesize complex cohomology}}} 
{\h^{p}(\Xh;\C)} @>>>
\h^{p}(\Dh;\C) @. \, \longto \, \cdots \\
@.
@VV{\wr}V @VV{\wr}V @VV{\wr}V @. \\
\cdots \, \longto \, @.
\h^{p}\sing(\Xh,\Dh;\C) @>>>
\underset{\makebox[0cm]{\text{\rm\footnotesize singular cohomology}}}
{\h^{p}\sing(\Xh;\C)} @>>>
\h^{p}\sing(\Dh;\C) @. \, \longto \, \cdots \, .
\end{CD}
$$
\end{theorem}

\begin{proof}
Combine 
propositions \ref{compalgan}
and \ref{compancplx}
with the result from Subsection \ref{compcplxsing}. 
\end{proof}

\subsection{An Alternative Description of the Comparison
  Isomorphism}

For computational purposes, it will be useful to describe the
isomorphism from the Comparison Theorem \ref{comp}
$$
\hDR{\Bul}(\Xh;\C) \iso \h\bul\sing(\Xh;\C)
$$
more explicitely. 
We will formulate our result in Proposition \ref{compplus} and explain its
significance in Subsection \ref{alteriso}.
We start by defining yet another cohomology theory.

\subsubsection{Smooth Singular Cohomology}
\label{smoothsubsec}
\index{Smooth singular cohomology}

Let $M$ be a complex manifold.
We write $\Triangle_p$ for the standard $p$-simplex spanned by the
basis 
$\{
(0,\ldots,0,\underset{i}{1},0,\ldots,0)
\st 
i=1,\ldots,p+1\}$ 
of $\R^{p+1}$.
A singular $p$-simplex 
$\sigma: \Triangle_p \to M$
is said to be {\em smooth},
\index{Simplex!smooth}
if $\sigma$ extends to a smooth map defined on a
neighbourhood of the standard $p$-simplex $\Triangle_p$ in its
$p$-plane.
We denote by 
$\c^{\infty}_p(M;\C)$ 
\notation{C}{$\c^{\infty}_p(M;\C)$}{\C-vector space of smooth %
  singular $p$-chains of a complex manifold $M$}%
the \C-vector space generated by
all smooth $p$-simplices of $M$, thus getting a subcomplex
$\c^{\infty}_{\Bul}(M;\C)$ 
of the complex 
$\c^{\rm sing}_{\Bul}(M;\C)$%
\notation{C}{$\c^{\rm sing}_p(M;\C)$}{\C-vector space of singular %
$p$-chains of a space $M$} %
of all singular simplices on $M$
$$
\c_{\Bul}^{\infty}(M;\C)\subset\c_{\Bul}^{\rm sing}(M;\C).
$$
We can dualize this subcomplex
\notation{C}{$\c^p_{\infty}(M;\C)$}{\C-vector space of smooth singular $p$-cochains
  of a complex manifold $M$}
$$
\c_{\infty}\bul(M;\C) := \c^{\infty}_{\Bul}(M;\C)^{\vee},
$$
and define a complex of sheaves
\notation{C}{$\CC^p_{\infty} (M;\C)$}{sheaf of smooth singular %
$p$-cochains of a complex manifold $M$} %
$$
\CC\bul_{\infty} (M;\C) : V \mapsto \c\bul_{\infty}(V;\C)
\quad\text{for}\quad V\subseteq M\quad\text{open}.
$$
The sheaves $\CC^{p}_{\infty} (M;\Q)$ are flabby and define smooth
singular cohomology of $M$.

\begin{definition}[Smooth singular cohomology of a manifold]
\label{defsmooth1}
  Let $M$ be a complex ma\-ni\-fold. 
  We define {\em smooth singular cohomology} groups
  $$
  \h_{\infty}\bul(M;\C) := 
  \hh\bul\,\c_{\infty}\bul(M;\C)=
  \H\bul\bigl(M;\CC_{\infty}\bul(M;\C)\bigr).
  $$
\notation{H}{$\h_{\infty}\bul(M;\C)$}{smooth singular cohomology of a
complex manifold $M$}
\end{definition}

Since smooth singular cohomology coincides with singular cohomology
for smooth manifolds (cf. \cite[p.~291]{bredon}), we can show that
$\CC\bul_{\infty}(M;\C)$ is a resolution of the constant sheaf
$\C_M$ exactly as we did for $\CC\bul\sing(M;\C)$ in Lemma
\ref{singres}. Thus we have the following lemma. 

\begin{lemma}
\label{smoothres}
For $M$ a complex manifold, the following sequence is exact
$$
0 \longto \C_M \longto 
\CC^0_{\infty}(M;\C) \longto
\CC^1_{\infty}(M;\C) \longto \,\ldots\,.
$$
\end{lemma}

We are now going to extend the definition of smooth singular
cohomology to the case of a divisor and to the relative case.
Let $\Xh$ and $\Dh$ be as in Subsection \ref{situation}.

First, recall how $\wom{\Dh}$ was defined (cf. pages \pageref{defwom}
and \pageref{defwoman}): We composed the simplicial
scheme $D\bul$ with the functor $\iast \Omega\bul_{?^{\rm an}}$, then
obtained a double complex $\Omega^{\Bul,\Bul}_{\Dh{}\bul}$ by summing
up certain maps and finally denoted the total complex of
$\Omega^{\Bul,\Bul}_{\Dh{}\bul}$ by $\wom{\Dh}$. 

Similarly, we can use the functor $\iast \CC\bul_{\infty}(?^{\rm an};\C)$,
compose it to the simplicial scheme $D\bul$ and denote the
resulting double complex by 
$\CC^{\Bul,\Bul}_{\infty}(\Dh;\C)$.%
\notation{C}{$\CC^{\Bul,\Bul}_{\infty}(\Dh;\C)$}{double complex of %
sheaves of smooth singular cochains of $\Dh{}\bul$}%
\label{defccinf}%
For the cor\-res\-pon\-ding total
complex we write $\CC\bul_{\infty}(\Dh;\C)$. 
  \notation{C}{$\CC\bul_{\infty}(\Dh;\C)$}{complex of sheaves of %
    smooth singular cochains of a divisor $\Dh$} %

We define $\CC\bul_{\infty}(\Xh,\Dh;\C)$ similar to $\wom{\Xh,\Dh}$
(cf. pages \pageref{OmXD} and \pageref{defwoman}):
The natural restriction maps
$$
\CC\bul_{\infty}(\Xh;\C) \to  \iast \CC\bul_{\infty}(\Dh_j;\C)
$$
sum up to a natural map of double complexes
\begin{equation}
\label{cinfXD}
\CC\bul_{\infty}(\Xh;\C)[0] \to \iast \CC^{\Bul,\Bul}_{\infty}(\Dh;\C),
\end{equation}
where we view $\CC\bul_{\infty}(\Xh;\C)[0]$ as a double complex concentrated
in the zeroth column. 
Taking the total complex in (\ref{cinfXD}) yields a natural map
$$
\CC\bul_{\infty}(\Xh;\C) \longto \iast \CC\bul_{\infty}(\Dh;\C),
$$
for whose mapping cone (cf. the appendix) we write
$\CC\bul_{\infty}(\Xh,\Dh;\C)$. %
\notation{C}{$\CC\bul_{\infty}(\Xh,\Dh;\C)$}{complex of sheaves of %
smooth singular cochains of a pair $(\Xh,\Dh)$} %

Thus we are able to formulate:

\begin{definition}[Smooth singular cohomology for the case of a
    divisor and for the relative case]
\label{defsmooth2}
\notation{H}{$\h\bul_{\infty}(\Dh;\C)$}{smooth singular cohomology of %
  a divisor $\Dh$} %
\notation{H}{$\h\bul_{\infty}(\Xh,\Dh;\C)$}{smooth singular cohomology %
  of a pair $(\Xh,\Dh)$} %
  We define {\em smooth singular cohomology} groups
  \begin{align*}
  & \h\bul_{\infty}(\Dh;\C)  :=
    \hh\bul\,\c\bul_{\infty}(\Dh;\C)
  =\H\bul(\Dh;\CC\bul_{\infty}(\Dh;\C))
  \quad\text{and} \\
  & \h\bul_{\infty}(\Xh,\Dh;\C)  :=
  \hh\bul\,\c\bul_{\infty}(\Xh,\Dh;\C)
  =\H\bul(\Xh;\CC\bul_{\infty}(\Xh,\Dh;\C)),
  \end{align*}
where $\Xh$ and $\Dh$ are as in Subsection \ref{situation}.
\end{definition}

\begin{remark}
Similarly one defines smooth singular cohomology with coefficients in
\Q. All statements made below about smooth singular cohomology 
remain valid for rational coefficients
except those involving holomorphic differentials.
\end{remark}

As usual, we have a long exact sequence in smooth singular cohomology
$$
\cdots\,\to\,\h^p_{\infty}(\Xh,\Dh;\C) \to \h^p_{\infty}(\Xh;\C) \to 
\h^p_{\infty}(\Dh;\C)\,\to\,\cdots\,.
$$

\subsubsection{Comparison with Analytic \deRham Cohomology}
\label{compdrsmooth}

Our motivation for considering {\rm smooth} singular cohomology
instead of ordinary singular cohomology comes from the following fact:
We can integrate differential-$p$-forms 
over smooth
$p$-simplices $\gamma$ and thus have a natural morphism
\begin{align*}
\Omega^p_{\Xh} & \to \CC^p_{\infty}(\Xh;\C) \\
\omega & \mapsto 
\left( 
\begin{gathered}
\CC_p^{\infty}(\Xh;\C) \to \C \\
\gamma \mapsto \int_{\gamma} \omega 
\end{gathered}
\right) .
\end{align*}
As a consequence of Stoke's theorem
$$
\int_{\del \gamma} \omega = \int_{\gamma} d \omega,
$$
we get a well-defined map of complexes
$$
\om{\Xh} \to  \CC\bul_{\infty}(\Xh;\C).
$$
This map is functorial, 
hence gives a natural transformation of functors
$$
\Omega\bul_{?^{\rm an}} \longto \CC\bul_{\infty}(?^{\rm an};\C).
$$
With the definitions of the complexes 
$$
\wom{\Dh},
\CC\bul_{\infty}(\Dh;\C)
\quad\text{and}\quad
\wom{\Xh,\Dh},
\CC\bul_{\infty}(\Xh,\Dh;\C)
$$ 
(cf. pages \pageref{defwoman} and \pageref{defccinf}),
we get a commutative diagram
\begin{equation}
\label{womccinf}
\begin{CD}
0 @>>> \iast\wom{\Dh}[-1] @>>> \wom{\Xh,\Dh} @>>> \om{\Xh} @>>> 0 \\
@. @VVV @VVV @VVV @. \\
0 @>>> \iast\CC\bul_{\infty}(\Dh;\C)[-1] @>>>
\CC\bul_{\infty}(\Xh,\Dh;\C) @>>> \CC\bul_{\infty}(\Xh;\C) @>>> 0\,.
\end{CD}\, 
\end{equation}
Recall the short exact sequence (\ref{cplxmapcon}) of the mapping cone
$M_{\C}$ from Subsection \ref{cplxcoho}
$$
0 \to \iast \C_{\Dh}[-1] \to M_{\C} \to \C_{\Xh}[0] \to 0
$$
and Diagram (\ref{cplxwom}) relating this short exact sequence with
the first line of (\ref{womccinf}). 
Combining (\ref{cplxwom}) and (\ref{womccinf}) gives a big diagram
\newcommand{\myarrow}[1]{
\makebox[0cm]{\setlength{\unitlength}{1pt}
\begin{picture}(0,0)(0,75)
\curve(1,70, 20,30, 1,0)
\curve(1,0, 5,0) \curve(1, 0, 2, 3)
\put(23,20){\makebox(0,0)[t]{#1}}
\end{picture}}}
\begin{equation}
\label{drsmooth}
\begin{CD}
0 @>>> 
\iast \C_{\Dh}[-1] \myarrow{$\alpha$} @>>> 
M_{\C} \myarrow{$\beta$} @>>> 
\C_{\Xh}[0] \myarrow{$\gamma$} @>>> 0 \\
@. @VVV @VVV @VVV @. \\
0 @>>> \iast\wom{\Dh}[-1] @>>> \wom{\Xh,\Dh} @>>> \om{\Xh} @>>> 0 \\
@. @VVV @VVV @VVV @. \\
0 @>>> \iast\CC\bul_{\infty}(\Dh;\C)[-1] @>>>
\CC\bul_{\infty}(\Xh,\Dh;\C) @>>> \CC\bul_{\infty}(\Xh;\C) @>>> 0
\end{CD}
\end{equation}
with $\alpha$, $\beta$, and $\gamma$ being just the composition of the
two vertical maps. 
We want to show that $\alpha$, $\beta$, and $\gamma$ are quasi-isomorphisms.
We already know that $\gamma$ is a
quasi-isomorphism by Lemma \ref{smoothres}. 
If $\alpha$ is also a quasi-isomorphism, then it
will follow from the $5$-lemma that $\beta$ is a quasi-isomorphism as
well. We split $\alpha$ as a composition of a quasi-isomorphism 
(cf. Proposition \ref{closedres}) and a
map $\wt{\alpha}$ 
$$
\alpha: \C_{\Dh}[0] \stackrel{\sim}{\longto} \CCh\bul(\DDD;\C)
\stackrel{\wt{\alpha}}{\longto} \CC\bul_{\infty}(\Dh;\C).
$$
The map $\wt{\alpha}$ is a map of total complexes induced by the morphism of double complexes
$$
\CCh\bul(\DDD;\C)[0] \longto \CC^{\Bul,\Bul}_{\infty}(\Dh;\C),
$$
where we consider $\CCh\bul(\DDD;\C)[0]$ as a double complex
concentrated in the zeroth row.
By Lemma \ref{smoothres} the sequence
$$
0 \longto \C_{\Dh_I} \longto \CC_{\infty}^0(\Dh_I;\C) \longto \CC_{\infty}^1(\Dh_I;\C) \longto
\, \ldots 
$$
is exact. Hence the sequence 
$$
0 \longto \CCh^q(\DDD;\C) \longto
\CC_{\infty}^{0,q}(\Dh;\C) \longto
\CC_{\infty}^{1,q}(\Dh;\C) \longto\,\ldots
$$
is also exact, since 
$\oplus$ and $\iast$ are
exact functors.
Applying Lemma \ref{lemma3} shows that $\wt{\alpha}$ is a
quasi-isomorphism.
Thus all vertical maps in diagram (\ref{drsmooth}) are indeed
quasi-isomorphisms. Taking hypercohomology provides us with
isomorphisms between long exact cohomology sequences generalizing the
Comparison Theorem \ref{comp}
\begin{equation}
\label{compdrsmoothdia}
\begin{CD}
\cdots \, @>>>
\h^{p}(\Xh;M_{\C}) @>>>
{\h^{p}(\Xh;\C)} @>>>
\h^{p}(\Dh;\C) @>>> \, \cdots \\
@.
@VV{\wr}V @VV{\wr}V @VV{\wr}V @. \\
\cdots \, @>>>
\hDR{p}(\Xh,\Dh;\C) @>>>
{\hDR{p}(\Xh;\C)} @>>>
\hDR{p}(\Dh;\C) @>>> \, \cdots \\
@. @VV{\wr}V @VV{\wr}V @VV{\wr}V @. \\
\cdots \, @>>>
\h^{p}_{\infty}(\Xh,\Dh;\C) @>>>
{\h^{p}_{\infty}(\Xh;\C)} @>>>
\h^{p}_{\infty}(\Dh;\C) @>>> \, \cdots \, .
\end{CD}
\end{equation}

\subsubsection{Comparison with Singular Cohomology}
\label{compsingsmooth}

We will compare smooth singular cohomology with ordinary singular
cohomology. 
First we need morphisms between the defining complexes.
We already have a natural map 
$$
\CC\bul\sing(\Xh;\C) \stackrel{\sim}{\longto}
\CC\bul_{\infty}(\Xh;\C)
$$
induced by the duals of the standard inclusions
$\c_{\Bul}^{\infty}(V;\C) \subset \c_{\Bul}^{\rm sing}(V;\C)$ 
for $V\subseteq\Xh$ open.

Next we discuss the case of a divisor:
By restricting singular cochains from $\Dh$ to an irreducible 
component $\Dh_j$ of $\Dh$, we
get a map
$$
\CC\bul\sing(\Dh;\C) \to \iast\CC\bul\sing(\Dh_j;\C),
$$
where $i:\Dh_j\inclusion\Dh$ is the natural inclusion.
By summing up the composition maps
$$
\CC\bul\sing(\Dh;\C) \to \iast\CC\bul\sing(\Dh_j;\C) \to
\iast\CC\bul_{\infty}(\Dh_j;\C),
$$
we obtain a morphism of complexes
$$
\CC\bul\sing(\Dh;\C) \to \CC^{\Bul,0}_{\infty}(\Dh;\C),
$$
which gives a morphism of double complexes
$$
\CC\bul\sing(\Dh;\C)[0] \to \CC^{\Bul,\Bul}_{\infty}(\Dh;\C),
$$
where we consider
$\CC\bul\sing(\Dh;\C)[0]$ as a double complex concentrated in the
zeroth column. 

Now we turn to the relative case of a pair $(\Xh;\Dh)$:
let $M\bul\sing$ denote the mapping cone of 
  \notation{M}{$M\bul\sing$}{mapping cone of 
  $\CC\bul\sing(\Xh;\C) \longto \iast \CC\bul\sing(\Dh;\C)$}
$$
\CC\bul\sing(\Xh;\C) \longto \iast \CC\bul\sing(\Dh;\C).
$$
By the functoriality of the mapping cone construction, we get a
commutative diagram
\newcommand{\primed}[1]{\textstyle #1 \hspace{0.2mm}'}
\newcommand{\dprimed}[1]{\textstyle #1 \hspace{0.2mm}''}
$$
\begin{CD}
0\,\longto\,\iast \CC\bul\sing(\Dh;\C)[-1] @>>> M\bul\sing 
@>>> \CC\bul\sing(\Xh;\C) \,\longto\,0 \\
\hspace{1cm} 
@V{\primed{\alpha}}VV 
@V{\primed{\beta}}VV 
\hspace{-1cm} 
@V{\primed{\gamma}}VV \\
0\,\longto\, \iast \CC\bul_{\infty}(\Dh;\C)[-1] @>>> \CC\bul_{\infty}(\Xh,\Dh;\C)
@>>> \hspace{2mm} \CC\bul_{\infty}(\Xh;\C) \hspace{0.2cm} \,\longto\,0 \, .
\end{CD}
$$
It is not hard to show directly that $\alpha'$, $\beta'$, $\gamma'$ 
are quasi-isomorphisms. 
However, we can use a shorter indirect argument:
We have the following commutative diagram
(cf. diagrams (\ref{cplxsing}), (\ref{drsmooth}) and the one above)
\renewcommand{\myarrow}[1]{
\makebox[0cm]{\setlength{\unitlength}{1pt}
\begin{picture}(0,0)(0,75)
\curve(1,70, 20,30, 1,0)
\curve(1,0, 5,0) \curve(1, 0, 2, 3)
\put(23,20){\makebox(0,0)[t]{$#1$}}
\end{picture}}}
$$
\begin{CD}
0\,\longto \hspace{7mm} 
\iast \C_{\Dh}[-1] \myarrow{\dprimed{\alpha}} \hspace{5mm} @>>> 
M_{\C} \myarrow{\dprimed{\beta}} @>>> \hspace{7mm} 
\C_{\Xh}[0] \myarrow{\dprimed{\gamma}} \hspace{5mm} \longto\,0 \\
\hspace{1cm} 
@V{\textstyle \alpha}V{\wr}V 
@V{\textstyle \beta}V{\wr}V 
\hspace{-1cm} 
@V{\textstyle \gamma}V{\wr}V \\
0\,\longto \hspace{2mm} \iast \CC\bul_{\infty}(\Dh;\C)[-1] @>>> 
\CC\bul_{\infty}(\Xh,\Dh;\C) @>>> 
\hspace{2mm} \CC\bul_{\infty}(\Xh;\C) \hspace{0.2cm} \longto\,0 \, \\
\hspace{1cm} 
@A{\primed{\alpha}}AA 
@A{\primed{\beta}}AA 
\hspace{-1cm} 
@A{\primed{\gamma}}AA \\
0\,\longto\,\iast \CC\bul\sing(\Dh;\C)[-1] @>>> 
M\bul\sing  @>>> 
\CC\bul\sing(\Xh;\C) \,\longto\,0
\end{CD}
$$
In Subsection \ref{compdrsmooth}, we have already seen that $\alpha$,
$\beta$, $\gamma$ are quasi-isomorphisms. Furthermore $\alpha''$ and
$\gamma''$ are quasi-isomorphisms (cf. Subsection \ref{compcplxsing}),
hence so is $\beta''$ (by the $5$-lemma). Therefore $\alpha'$,
$\beta'$, $\gamma'$ have to be quasi-isomorphisms as well.

Taking hypercohomology in this diagram 
gives us
\begin{equation}
\label{compsingsmoothdia1}
\begin{CD}
\cdots \, @>>>
\h^{p}(\Xh;M_{\C}) \myarrow{\dprimed{\beta}_{\ast}} @>>>
{\h^{p}(\Xh;\C)} \myarrow{\dprimed{\gamma}_{\ast}} @>>>
\h^{p}(\Dh;\C) \myarrow{\dprimed{\alpha}_{\ast}} @>>> \, \cdots \\
@. @VV{\wr}V @VV{\wr}V @VV{\wr}V @. \\
\cdots \, @>>>
\h^{p}_{\infty}(\Xh,\Dh;\C) @>>>
{\h^{p}_{\infty}(\Xh;\C)} @>>>
\h^{p}_{\infty}(\Dh;\C) @>>> \, \cdots \, \\
@. @VV{\wr}V @VV{\wr}V @VV{\wr}V @. \\
\cdots \, @>>>
\h^{p}(\Xh;M\bul\sing) @>>>
{\h^{p}\sing(\Xh;\C)} @>>>
\h^{p}\sing(\Dh;\C) @>>> \, \cdots \, .
\end{CD}
\end{equation}
In order to get rid of the mapping
cones $M_{\C}$ and $M\bul\sing$,
we apply Lemma \ref{homalg} two each of the two short exact sequences
of complexes below (where $\GG\bul_{?}$ denotes the Godement
resolution --- see page \pageref{godement})
$$
\begin{matrix}
0 
& \!\!\!\!\to\!\!\!\!
& \Gamma(\Xh;\GG\bul_{j_!\C_{\Uh}}) 
& \!\!\!\!\to\!\!\!\! 
& \Gamma(\Xh;\GG\bul_{\C_{\Xh}}) 
& \!\!\!\!\to\!\!\!\! 
& \Gamma(\Xh;\GG\bul_{\iast \C_{\Dh}}) 
& \!\!\!\!\to\!\!\!\! 
& 0 \\
& & \downarrow & & \downarrow & & \downarrow & & \\
0 
& \!\!\!\!\to\!\!\!\! 
& \Gamma(\Xh;\GG\bul_{\CC\bul\sing(\Xh,\Dh;\C)})
& \!\!\!\!\to\!\!\!\! 
& \Gamma(\Xh;\GG\bul_{\CC\bul\sing(\Xh;\C)}) 
& \!\!\!\!\to\!\!\!\! 
& \Gamma(\Xh;\GG\bul_{\iast\CC\bul\sing(\Dh;\C)}) 
& \!\!\!\!\to\!\!\!\! 
& 0
\end{matrix}
$$
and obtain a commutative diagram (cf. definitions \ref{cplxcoho} and \ref{relsingcoho})
\begin{diagram}[PostScript=dvips,landscape,size=2em,LaTeXeqno]
\cdots
& \rTo
& \h^{p}(\Xh,\Dh;\C) 
& 
& \rTo
& 
& \h^{p}(\Xh;\C)
& 
& \rTo
&
& \h^{p}(\Dh;\C)
& \rTo \cdots
& {}
&
&
\\ 
&
& \vLine
& \rdTo
& 
& 
& \vLine
& \rdTo^{\eps_p}
&
&
& \vLine 
& \rdTo^{\eps_p}
&
&
&
\\ 
& 
& \cdots
& \rTo
& \h^{p}(\Xh;M_{\C}) 
& \rTo
& \HonV 
& 
& \h^{p}(\Xh;\C) 
& \rTo
& \HonV
&
& \h^{p}(\Dh;\C)
& \rTo
& \cdots
\\ 
&
& \dTo
& 
& \dTo_{\dprimed{\beta}_{\ast}}
& 
& \dTo
& 
& \dTo_{\dprimed{\gamma}_{\ast}}
&
& \dTo
&
& \dTo_{\dprimed{\alpha}_{\ast}}
&
&
\\ 
\cdots
& \rTo
& \h^{p}\sing(\Xh,\Dh;\C) 
& \hLine 
& \VonH 
& \rTo
& \h^{p}\sing(\Xh;\C)
& \hLine
& \VonH  
& \rTo
& \h^{p}\sing(\Dh;\C)
& \rTo \cdots
& \VonH 
&
\\ 
&
&
& \rdTo
& 
& 
& 
& \rdTo^{\eps_p}
& 
&
&
& \rdTo^{\eps_p}
&
&
&
\\ 
&
& \cdots
& \rTo 
& \h^{p}(\Xh;M\bul\sing) 
& 
& \rTo
& 
& \h^{p}\sing(\Xh;\C)
&
& \rTo
&
& \h^{p}\sing(\Dh;\C) 
& \rTo
& \cdots
\label{compsingsmoothdia2}
\end{diagram}

Thus we have the following supplement 
to the Comparison Theorem \ref{comp}.

\subsubsection{Supplement to the Comparison Theorem \ref{comp}}

\begin{proposition}
\label{compplus}
With the notation of Theorem \ref{comp} and definitions
\ref{defsmooth1}, \ref{defsmooth2}, we have the following
commutative diagram, all whose vertical maps are isomorphisms
\renewcommand{\myarrow}{
\makebox[0cm]{\setlength{\unitlength}{1pt}
\begin{picture}(0,0)(0,115)
\curve(1,110, 20,60, 1,0)
\curve(1,0, 4,1) \curve(1, 0, 1, 4)
\end{picture}}}
$$
\begin{CD}
\cdots \,\longto\,
\h^{p}(\Xh,\Dh;\C) \myarrow @>>>
{\h^{p}(\Xh;\C)} \myarrow @>>>
\h^{p}(\Dh;\C) \myarrow \,\longto\, \cdots \\
\hspace{1.2cm} @VV{\wr}V @VV{\wr}V \hspace{-1.2cm} @VV{\wr}V \\
\cdots \,\longto\,
\hDR{p}(\Xh,\Dh;\C) @>>>
{\hDR{p}(\Xh;\C)} @>>>
\hDR{p}(\Dh;\C) \,\longto\, \cdots \\
\hspace{1.2cm} @VV{\wr}V @VV{\wr}V \hspace{-1.2cm} @VV{\wr}V\\
\cdots \,\longto\,
\h^{p}_{\infty}(\Xh,\Dh;\C) @>>>
{\h^{p}_{\infty}(\Xh;\C)} @>>>
\h^{p}_{\infty}(\Dh;\C) \,\longto\, \cdots \\
\hspace{1.2cm} @AA{\wr}A @AA{\wr}A \hspace{-1.2cm} @AA{\wr}A\\
\cdots \,\longto\,
\h^{p}\sing(\Xh,\Dh;\C) @>>>
{\h^{p}\sing(\Xh;\C)} @>>>
\h^{p}\sing(\Dh;\C) \,\longto\, \cdots \, .
\end{CD}
$$
\end{proposition}

\begin{proof}
Combine diagrams (\ref{compdrsmoothdia}),
(\ref{compsingsmoothdia1}) and (\ref{compsingsmoothdia2}).
\end{proof}

\subsubsection{Application of the Supplement to the Comparison
  Theorem}
\label{alteriso}

With Proposition \ref{compplus}, we can compute the isomorphism from
the Comparison Theorem \ref{comp}
$$
\begin{CD}
\substack{\rm analytic\\ \rm\deRham \rm\ cohomology}
@<{\sim}<<
\substack{\rm complex\\ \rm cohomology}
@>{\sim}>>
\substack{\rm singular \\ \rm cohomology}
\end{CD}
$$
alternatively as
$$
\begin{CD}
\substack{\rm analytic\\ \rm\deRham \rm\ cohomology}
@>{\sim}>{\alpha}>
\substack{\rm smooth\ singular\\ \rm cohomology}
@<{\sim}<{\beta}<
\substack{\rm singular \\ \rm cohomology}.
\end{CD}
$$
Why should this be useful? We show that we can easily find  the image of a
\deRham\ cohomology class under the second isomorphism
$\beta^{-1}\circ\alpha$ 
by computing a few integrals.

First, we define {\em smooth singular homology}
\index{Smooth singular homology}
$\h^{\infty}_{\Bul}$
as the dual of $\h\bul_{\infty}$ 
(cf. definitions \ref{defsmooth1}, \ref{defsmooth2}), 
or equivalently:

\begin{definition}[Smooth singular homology]
\label{smoothsinghom}
  \notation{H}{$\h^{\infty}_{\Bul}(M;\C)$}{smooth singular homology of
  a complex manifold $M$}
  \notation{H}{$\h^{\infty}_{\Bul}(\Dh;\C)$}{smooth singular homology of
  a divisor $\Dh$}
  \notation{H}{$\h^{\infty}_{\Bul}(\Xh,\Dh;\C)$}{smooth singular homology of
  a pair $(\Xh,\Dh)$}
For $\Xh$, $\Dh$ as usual (see Subsection \ref{situation}), we define
  {\em smooth singular homology groups}
\begin{align*}
& \h^{\infty}_{\Bul}(\Xh;\C) := \hh\bul\, \c^{\infty}_{\Bul}(\Xh;\C) \\
& \h^{\infty}_{\Bul}(\Dh;\C) := \hh\bul\, \tot \oplus_{|I|=\,\Bul+1}\,
  \iast \c^{\infty}_{\Bul}(\Dh_I;\C) \\
& \h^{\infty}_{\Bul}(\Xh,\Dh;\C) := 
\hh\bul \,{\rm mapping\ cone} \left( 
\tot \oplus_{|I|=\,\Bul+1}\, \iast \c^{\infty}_{\Bul}(\Dh_I;\C) \to
\c^{\infty}_{\Bul}(\Xh;\C)
\right).
\end{align*}
\end{definition}

This gives immediately an isomorphism between singular homology and
its smooth version.

Now if $[\omega]$ is a \deRham cohomology class in the image of 
$$
\hh\bul\Gamma(\Xh;\om{\Xh}) \longto \hDR{\Bul}(\Xh;\C),
$$
represented by a differential form $\omega\in\Gamma(\Xh;\om{\Xh})$, we can find
the image $(\beta^{-1}\circ\alpha)[\omega]$ of $[\omega]$ in
$\h\bul\sing(\Xh;\C)$ as follows:
Let $\{\gamma_1,\ldots,\gamma_t\}$ be a basis of
$\h^{\infty}_{\Bul}(\Xh;\C)$ and
$\{\gamma^{\ast}_1,\ldots,\gamma^{\ast}_t\}$ its dual. Then 
$\{\beta^{-1}\gamma^{\ast}_1,\ldots,\beta^{-1}\gamma^{\ast}_t\}$ is a
basis of $\h\bul\sing(\Xh;\C)$ and we can express
$(\beta^{-1}\circ\alpha)[\omega]$ in terms of this basis
\begin{align*}
(\beta^{-1}\circ\alpha)[\omega] & = 
\sum_{i=1}^t \,\langle\,
\beta\gamma_i\,,\,(\beta^{-1}\circ\alpha)[\omega]\,\rangle\,\cdot\,
\beta^{-1}\gamma^{\ast}_i \\
& = \sum_{i=1}^t \,
\langle\, \gamma_i\,,\,\alpha[\omega]\,\rangle
\,\cdot\, \beta^{-1}\gamma^{\ast}_i \\
& = \sum_{i=1}^t \,\left(\, \int_{\gamma_i} \omega \,\right) \,\cdot\,
\beta^{-1}\gamma^{\ast}_i ,
\end{align*}
where $\langle\,\cdot\,,\,\cdot\,\rangle$
\notation{0}{$\langle\,\cdot\,,\,\cdot\,\rangle$}{homology-cohomology-pairing}
denotes the homology-cohomology-pairing for singular and smooth
singular (co)homology.

Thus, if $[\omega]\in\hDR{\Bul}(\Xh;\C)$ is a cohomology class
represented by a differential form $\omega\in\Gamma(\Xh;\om{\Xh})$, 
we can find its image $(\beta^{-1}\circ\alpha)[\omega]$ in $\h\bul\sing(\Xh;\C)$
by computing the integrals $\int_{\gamma_i}\omega$ for $i=1,\ldots,t$.
For an application, see Example \ref{noiso}.\\[0.5mm]
The case of a divisor and the relative case can be dealt with
analogously.\label{pper}
The only difference is that a cohomology class $[\omega]$ in the image of
$$
\hh\bul\Gamma(\Dh;\wom{\Dh}) \longto \hDR{\Bul}(\Dh;\C)
$$
or
$$
\hh\bul\Gamma(\Xh;\wom{\Xh;\Dh}) \longto \hDR{\Bul}(\Xh,\Dh;\C)
$$
is represented by a formal sum $\oplus_I\omega_I$ of differential forms $\omega_I$ living on either $\Xh$
or one of the various $\Dh_I$. 
Si\-mi\-lar\-ly, a homology class
$\gamma\in\h^{\infty}_{\Bul}(\Dh;\C)$ or
$\gamma\in\h^{\infty}_{\Bul}(\Xh,\Dh;\C)$ is represented by a formal
sum 
$\oplus_I\Gamma_I$
of 
smooth simplicial chains $\Gamma_I$ living on $\Xh$ and the $\Dh_I$'s.
In computing the pairing 
$\langle\gamma,[\omega]\rangle$, 
a summand
$\langle\gamma_J,[\omega_I]\rangle$ 
can only give a non-zero contribution
if the domains of definition coincide $\Dh_I=\Dh_J$. 
(Here $\omega_I$ lives on $\Dh_I$ with the convention
that $\omega_{\emptyset}$ corresponds to $\Xh$.)

\subsection{A Motivation for the Theory of Periods}

Assume we are dealing with a smooth variety $X_0$ defined over \Q\
and a normal-crossings-divisor $D_0$ on $X_0$. We denote the
base change to \C\ by $X$ and $D$, respectively.

It seems a natural question to ask whether the isomorphism
\begin{align*}
\hDR{\Bul}(X_0,D_0/\Q) \tensor{\Q} \C = 
& \hDR{\Bul}(X,D/\C) \\
& \hspace{1cm} \downiso \\
\h\bul\sing(\Xh,\Dh;\Q) \tensor{\Q} \C = & \h\bul\sing(\Xh,\Dh;\C)
\end{align*}
is induced by a natural isomorphism of the form
$$
\hDR{\Bul}(X_0,D_0/\Q) \to \h\bul\sing(\Xh,\Dh;\Q).
$$
The following example shows that the answer to this question is {\bf
no} and this negative result may be regarded as the starting point
of the theory of periods.
(This example uses absolute (co)homology, but taking $D_0:=\emptyset$
we can reformulate it for relative (co)homology.)

\begin{example}
\label{noiso}
Let
\notation{A}{$\A^n_k$}{affine $n$-space over a field $k$}%
$$
X_0 := \Spec \Q[t,t^{-1}] =\A^1_{\Q} \setminus \{0\}
$$
be the affine line with point $0$ deleted. 

Then the singular homology group $\h_1^{\rm sing}(\Xh;\Q)$ of $\Xh =
 \C^{\times}$ is generated by the unit circle
$\sigma:=S^1$. Hence
$\h^1\sing(\Xh;\Q)$ is generated by its dual $\sigma^{\ast}$.

For the algebraic \deRham cohomology group $\hDR{1}(\Xo/\Q)$ of $\Xo$ we get
\begin{align*}
\hDR{1}(\Xo/\Q) & = \H^{1}(\Xo;\om{\Xo/\Q}) \\
& \stackrel{(\ast)}{=} \hh^{1}\Gamma(\Xo;\om{\Xo/\Q}) \\
& = \coker( d: \Q[t,t^{-1}] \to \Q[t,t^{-1}]dt) \\
& = \Q \frac{dt}{t} .
\end{align*}
Here we used at $(\ast)$ that the sheaves $\Omega^p_{X_0/\Q}$ are
quasi-coherent, hence acyclic for the global section functor
$\Gamma(X_0;?)$ by 
\cite[Thm. III.3.5, p.~215]{hartshorne}, 
since $X_0$ is affine.

Under the isomorphism $\hDR{1}(X/\C) \stackrel{\sim}{\longrightarrow}
\h^{1}\sing(X;\C)$ the generator $\frac{dt}{t}$ is mapped to $2\pi i\,
\sigma^{\ast}$ because of $\int_{\sigma} \frac{dt}{t} = 2\pi i$ . But
$2\pi i$ is not a rational number, hence the isomorphism
$$
\hDR{1}(X/\C) {\to} \h^{1}\sing(X;\C)
$$
is {\bf not} induced by a map
$$
\hDR{1}(\Xo/\Q) {\to} \h^{1}\sing(X;\Q).
$$ 
The complex number $2\pi i$
is our first example of a period.
\end{example}

\section{Definitions of Periods}
\label{per123}

A period is to be thought of as an integral that occurs in a geometric context.
In their papers \cite{kontsevich} and \cite{kontsevich_zagier},
Kontsevich and Zagier list various ways of how to define a period.

It is stated in their papers without reference, that all these variants give the same 
definition.
We give a partial proof of this statement 
in the Period Theorem \ref{perioperio}.

\subsection{First Definition of a Period: Pairing Periods}

Let $X_0$ be a smooth variety defined over \Q\ and $D_0$ a divisor 
with normal crossings on $X_0$. We denote by $\Xh$ and $\Dh$ the
complex analytic spaces (cf. Subsection \ref{defcomplexanspace}) associated to the base change to $\C$ of $X_0$
and $D_0$.

From the discussion of \deRham cohomology in the previous sections
(especially Pro\-po\-si\-tion \ref{basechange} and Theorem \ref{comp}), 
we see that inside the \C-vector space $\hDR{p}(\Xh,\Dh;\C)$ there sits
a \Q-lattice $\hDR{p}(X_0,D_0/\Q)$ of full rank. 
Since 
\begin{align*}
\hDR{\Bul}(\Xh,\Dh;\C) & \stackrel{(\ref{comp})}{\iso} \h\bul\sing(\Xh,\Dh;\C) 
 = \h\bul\sing(\Xh,\Dh;\Q) \tensor{\Q} \C,
\end{align*}
the classical perfect pairing
$$
\langle \cdot,\cdot \rangle : 
\h_{\Bul}^{\rm sing}(\Xh,\Dh;\Q) \times \h\bul\sing(\Xh,\Dh;\Q) 
 \to \Q
$$
gives us a new non-degenerate pairing
\begin{equation}
\label{pairing}
\langle \cdot,\cdot \rangle : 
\h_{\Bul}^{\rm sing}(\Xh,\Dh;\Q) \times  \hDR{\Bul}(X_0,D_0/ \Q) 
\to \C,
\end{equation}
which is natural in $(X_0,D_0)$. 

\begin{definition}[Pairing period, short: p-period]
\label{per1}
\index{Period!pairing}\index{p-period|see{Period, pairing}}
The complex numbers 
$$
\langle \gamma,\omega_0 \rangle \in \C
\quad \text{for} \quad 
\omega_0 \in \hDR{\Bul}(X_0,D_0/ \Q), \,
\gamma \in \h_{\Bul}^{\rm sing}(\Xh,\Dh;\Q),
$$
which appear in the image of the natural pairing (\ref{pairing}), are
called {\em p-periods}.
The set of all p-periods for all pairs $(X_0,D_0)$ will be denoted by 
$\P_p$.
\notation{P}{$\P_p$}{set of all pairing periods}

We call a p-period $\langle\gamma,\omega_0\rangle$ {\em special} 
\index{Period!pairing!special}\index{p-period@${\rm p}'$-period|see{Period, pairing, special}}
or a 
{\em ${\rm p}'$-period}, if\/ $\omega_0$ is contained in the image of
$$
\hh\bul\Gamma(X_0;\wom{X_0,D_0/\Q})\to\hDR{\Bul}(X_0,D_0/\Q)
$$ 
and write $\P^{\,\prime}_p$ for the set of all special pairing periods
for all $X_0,D_0,\omega_0,\gamma$.
\notation{P}{$\P^{\,\prime}_p$}{set of all special pairing periods}

Running through bases
\begin{center}
$\{\omega_1,\ldots,\omega_t\}$
of\/ 
$\hDR{\Bul}(X_0,D_0/ \Q)$ and 
$\{\gamma_1,\ldots,\gamma_t\}$ 
of\/
$\h_{\Bul}^{\rm sing}(\Xh,\Dh;\Q)$
\end{center}
yields the so-called {\em period matrix} $P:=\bigl(\langle \gamma_i,\omega_j\rangle\bigr)_{i,j=1}^t$ of $(X_0,D_0)$.
\index{Period!matrix}
\end{definition}

\begin{remark}
\label{det}
As a consequence of the Comparison Theorem \ref{comp}, the period
matrix must be a square matrix
(cf. \cite[p. 63]{kontsevich}). In loc.~cit. also a statement about its
determinant is made, being a square root of a rational number times a
power of $2\pi i$
$$
\text{determinant} = \pm \sqrt{a} \cdot (2\pi i)^n \quad
\text{for} \quad a\in\Q^{\times},\, n\in\N_0.
$$
An indication of a proof is not given there.
\end{remark}

\subsection{Second Definition of a Period: Abstract Periods}
\label{per2sec}

There exists another definition of a period, which also emphasizes
their geometric origin but does not involve algebraic \deRham
cohomology.

For this definition the following data is needed,
\begin{itemize}
\item
$X_0$ a smooth algebraic variety defined
over {\Q} of dimension $d$,
\item
$D_0$ a divisor on $X_0$ with normal crossings,
\item
$\omega_0 \in \Gamma(X_0;\Omega^d_{X_0/\Q})$ an algebraic differential form of top
degree,
\item
$\gamma \in \h_d^{\rm sing}(\Xh,\Dh;\Q)$ a homology class of singular chains
on the complex manifold $\Xh$ with boundary on the divisor $\Dh$.
\end{itemize}

As usual, we denote by $\Xh$ and $\Dh$ the
complex analytic spaces (cf. Subsection \ref{defcomplexanspace}) associated to the base change to $\C$ of $X_0$
and $D_0$.
Thus, we are dealing with objects $\Xh$, $\Dh$ and $\gamma$ of
real dimension $2d$, $2d-2$ and $d$, respectively.

The differential form $\omega_0$ on $X_0$ gives rise to a differential
form 
$$
\omega := \pi^{\ast} \omega_0
$$
on the base change $X=X_0\times_{\Q}\C$, where $\pi: X \rightarrow X_0$ is the natural projection.
Now $\omega_0$ is closed for dimension reasons,
hence $\omega$ is closed as well. 

We choose a representative $\Gamma\in\alpha^{-1}\gamma$ of the
preimage $\alpha^{-1}\gamma$ of $\gamma$ under the isomorphism
(cf. Proposition \ref{compplus})
$$
\alpha: \h_d^{\infty}(\Xh,\Dh;\Q) \stackrel{\sim}{\longto} \h_d^{\rm
  sing}(\Xh,\Dh;\Q)
$$
and define
$$
\int_{\gamma}\omega:=\int_{\Gamma_d}\omega,
$$
where the chain $\Gamma_d$ 
  \notation{G}{$\Gamma_p$}{top-dimensional part of a 
  $p$-chain $\Gamma\in\c^{\infty}_p(\Xh,\Dh;\Q)$}
consists of the $d$-simplices of $\Gamma$
(that is we ignore those $d-q-1$-simplices of $\Gamma$ that live on one of the 
$\Dh_I$ for $|I|=q+1,\, q\ge 0$ --- see page \pageref{shorthand} for
the definition of $D_I$).

\begin{proof}[Observe that this integral $\int_{\gamma}\omega$ is well-defined]
\label{defev1}
The restriction of $\omega$ to some irreducible component $\Dh_j$ of
$\Dh$ is a holomorphic $d$-form on a complex manifold of dimension
$d-1$, hence zero. Therefore the integral $\int_{\triangle} \omega$
evaluates to zero for smooth singular simplices $\triangle$ that are
supported on $\Dh$.
Now if $\Gamma'$, $\Gamma''$ are two representatives for $\alpha^{-1}\gamma$, we have
$$
\Gamma'_d - \Gamma''_d \sim \del (\Gamma_{d+1})
$$
modulo simplices living on some $\Dh_I$ for a smooth singular chain 
$\Gamma$ of dimension $d+1$
$$
\Gamma \in {\rm C}_{d+1}^{\infty} (\Xh,\Dh;\Q).
$$
Using Stoke's theorem we get 
$$
\int_{\Gamma'_d} \omega - \int_{\Gamma''_d} \omega = 
\int_{\del (\Gamma_{d+1})} \omega =
\int_{\Gamma_{d+1}} d \omega = 0 ,
$$
since $\omega$ is closed.
\end{proof}

\begin{definition}[Abstract period, short: a-period]
\label{per2}
\index{Period!abstract}\index{A-period@a-period|see{Period, abstract}}
We will call the complex number
$\int_{\gamma}\omega$ the {\em a-period} of the quadruple
$(X_0,D_0,\omega_0,\gamma)$ and denote the set of all a-periods 
for all $(X_0,D_0,\omega_0,\gamma)$ by $\P_a$. 
\notation{P}{$\P_a$}{set of all abstract periods}
\end{definition}

\begin{remark}
This definition was motivated by Kontsevich's discussion of effective
periods \cite[def.~20, p.~62]{kontsevich},
which we will partly quote in subsection \ref{per4}. 
Modulo Conjecture \ref{ev}, abstract and effective periods are
essentially the same.
\end{remark}

We can ask whether a different definition of
the tuples $(X_0,D_0,\omega_0,\gamma)$ yields more period values. 
A partial answer is given by the following remarks.

\begin{remark}{\bf (cf. \cite[p. 62]{kontsevich})}
\label{periorem1}
For example, we could have considered algebraic varieties defined over
\QQ. However, doing so does not give us more period
values.
\end{remark}

{\em Proof.}

\begin{itemize}
\item Any variety $X_0$ defined over $\ol{\Q}$ is already defined over
  a finite extension $\Q^{\,\prime}$ of \Q, i.e. there exists a variety $X_0'$
  such that
  $$
  X_0 = X_0' \times_{\Q^{\,\prime}} \ol{\Q} .
  $$
  Via
  $$
  X_0' \longrightarrow \Spec \Q^{\,\prime} \longrightarrow \Spec \Q
  $$
  we may consider $X_0'$ as a variety defined over \Q.
\item An analogous argument applies to a normal-crossings-divisor
  $D_0$ on $X_0$. By possibly replacing $\Q^{\,\prime}$ by a finite
  extension field, we can therefore assume w.l.o.g. that $D_0$ is the base
  extension
  of a normal-crossings-divisor $D_0^{\,\prime}$ on $X_0'$
  $$
  D_0 = D_0^{\,\prime} \times_{\Q} \ol{\Q} .
  $$
\item Also any differential $d$-form $\omega_0 \in \Gamma(X_0;\Omega^d_{X_0/\ol{\Q}})$ can
  be considered as a differential form $\omega_0'$ on $X_0'$ by
  eventually replacing  $\Q^{\,\prime}$ by a finite extension.
\end{itemize}
Finally, we discuss singular homology classes.
\begin{itemize}
\item 
First observe that 
\begin{align*}
\vert \Spec \C \tensor{\Q} \Q^{\,\prime} \vert & = \Hom_{\C} ( \C \tensor{\Q}
\Q^{\,\prime},\C) \\
& = \Hom_{\Q}(\Q^{\,\prime},\C)
\end{align*}
is a finite discrete set of points. Thus $X'$ is a finite union of
disjoint copies of $X$
\begin{align*}
X' &= X_0' \times_{\Q} \C \\
& = X_0' \times_{\Q^{\,\prime}} \Q^{\,\prime} \times_{\Q} \C \times_{\C} \C \\
& = (X_0' \times_{\Q^{\,\prime}} \C) \times_{\C} 
\bigl( \Q^{\,\prime} \times_{\Q} \C) \bigr) \\
&= X \times \Hom_{\Q}(\Q^{\,\prime},\C) \\
&= \coprod_{\sigma\,:\, \Q^{\,\prime} \rightarrow \C} X .
\end{align*}
Similarly, $D^{\,\prime}=\coprod_{\sigma} D$. Therefore we have 
$$
\h^{\rm sing}_d(X^{\prime\,\rm an},D^{\,\prime\,\rm an}; \Q) 
  = \underset{\sigma}{\bigoplus}\, \h^{\rm sing}_d(\Xh,\Dh;\Q).
$$
If $\gamma \in \h^{\rm sing}_d(\Xh,\Dh;\Q)$ is some singular homology
class, we can pick any $\sigma':\Q^{\,\prime} \rightarrow \C$ and find
$$
\gamma':=\{ \sigma' \} \times\gamma \in 
\h^{\rm sing}_d(X^{\prime\,\rm an},D^{\,\prime\,\rm an};\Q)
=\oplus_{\sigma} \,\h^{\rm sing}_d(\Xh,\Dh;\Q).
$$ 
\end{itemize}
Now
$$
\int_{\gamma'} \omega' = \int_{\gamma} \omega
$$
which finishes the proof.{\hfill $\square$}

\begin{remark}{\bf (cf. \cite[p. 62]{kontsevich})}
\label{periosmooth}
The requirement of $X$ to be smooth was also made only for convenience
and is not a real restriction.
\end{remark}

We postpone the proof of Remark \ref{periosmooth} until the Period Theorem
\ref{perioperio} is proved.

\subsection{Third Definition of a Period: Na{\"\i}ve Periods}

For the last definition of a period, we need the notion of
semi-algebraic sets.

Let \K\ be any field contained in \R.

\begin{definition}[\K-semi-algebraic sets,
  \label{semialgset}
  {\cite[Def. 1.1, p. 166]{hironaka_74}}]
  \index{Semi-algebraic!set}
  A subset of $\R^n$ is said to be {\em \K-semi-algebraic}, if it is of the
  form
  $$
  \{ \underline{x} \in  \R^n \vert f(\underline{x}) \ge 0\}
  $$
  for some polynomial $f \in \K [ x_1, \ldots, x_n ]$ or can be
  obtained from sets of this form in a finite number of steps, where
  each step consists of one of the following basic operations:
  \renewcommand{\labelenumi}{(\roman{enumi})}
  \begin{enumerate}
  \item complementary set,
  \item finite intersection,
  \item finite union.
  \end{enumerate}
  \renewcommand{\labelenumi}{(\arabic{enumi})}
\end{definition}

Denote the integral closure of \Q\ in \R\ by \QR.
Note that \QR\ is a field.
\notation{Q}{\QR}{integral closure of \Q\ in {\R}}
Now we can define
\begin{definition}[Na{\"\i}ve periods, short: n-periods]
\label{per3}
\index{Period!na{\"\i}ve}\index{n-period|see{Period, na{\"\i}ve}}
Let 
\begin{itemize}
\item
  $G\subseteq\R^n$ be an oriented compact \QR-semi-algebraic set 
  which is equidimensional of dimension $d$,
  and
\item
  $\omega_0$ a rational differential $d$-form on $\R^n$ having coefficients
  in \QQ,
  which does not have poles on $G$. 
\end{itemize}
Then we call the complex number
$\int_{G}\omega_0$ a {\em n-period} and denote the set of all
n-periods for all $G$ and $\omega_0$ by
$\P_n$.
\notation{P}{$\P_n$}{set of all na{\"\i}ve periods}
\end{definition}

This set $\P_n$ enjoys additional structure.

\begin{proposition}
\label{palgebra}
The set $\P_n$ is a $\ol{\Q}$-algebra.
\end{proposition}
\begin{proof} 
  {\em Additive structure:} 
  Let 
  $\int_{G_1} \omega_1$ and 
  $\int_{G_2} \omega_2 \in \P_n$ 
  be periods with domains of integration 
  $G_1 \subseteq \R^{n_1}$ and
  $G_2 \subseteq \R^{n_2}$. Using the inclusions
  \begin{align*}
  i_1 : \R^{n_1} & \iso \R^{n_1}\times \{1/2\}\times \{\ul{0}\} \subset
    \R^{n_1} \times \R \times \R^{n_2}\quad\text{and} \\
  i_2 : \R^{n_2} & \iso \{\ul{0}\} \times \{-1/2\}\times \R^{n_2} \subset
    \R^{n_1} \times \R \times \R^{n_2},
  \end{align*}
  we can write $i_1(G_1)\union i_2(G_2)$ for the disjoint union of
  $G_1$ and $G_2$. 
  With the projections $p_j: \R^{n_1}\times\R\times\R^{n_2} \to
  \R^{n_j}$ for $j=1,2$, we can lift $\omega_j$ on $\R^{n_j}$ to
  $p_j^{\ast}\omega_j$ on $\R^{n_1}\times\R\times\R^{n_2}$. For
  $q_1,q_2\in\ol{Q}$ we get
  $$
  q_1\! \int_{G_1} \omega_1 + q_2\! \int_{G_2} \omega_2 = \\
  \int_{i_1(G_1)\union i_2(G_2)} q_1\cdot(1/2+t)\cdot
  p_1^{\ast}\omega_1+ q_2\cdot(1/2-t)\cdot p_2^{\ast}\omega_2 \in
  \P_n, 
  $$
  where $t$ is the coordinate of the ``middle'' factor \R\ of 
  $\R^{n_1}\times \R \times \R^{n_2}$. This shows that $\P_n$ is a
  \QQ-vector space. 

  {\em Multiplicative structure:}
  In order to show that $\P_n$ is closed under
  multiplication, we write
  $$
  p_i : \R^{n_1} \times \R^{n_2} \longrightarrow \R^{n_i},\quad i=1,2
  $$
  for the natural projections and obtain
  $$
  \left( \int_{G_1} \omega_1 \right) \cdot \left( \int_{G_2} \omega_2 \right) =
  \int_{G_1 \times G_2} p_1^{\ast} \omega_1 \wedge p_2^{\ast} \omega_2 \in \P_n
  $$
  by the Fubini formula. 
\end{proof}

\begin{remark}
\label{kzdef}
The Definition \ref{per3} was inspired by the one given in
\cite[p.~772]{kontsevich_zagier}
\begin{quote}
{\bf Definition.} {\em A {\em [na{\"\i}ve] period} is a complex number 
whose real and imaginary
part are values of absolutely convergent integrals of rational
functions with rational coefficients, over domains in $\R^n$ given by
polynomial inequalities with rational coefficients.}
\end{quote}
We will not work with this definition and use the modified version
\ref{per3} instead, since this gives us more flexibility in proofs.

From the statements made in
\cite[p.~773]{kontsevich_zagier} it would follow that both
definitions of a na{\"\i}ve period agree.
\end{remark}

Examples of na{\"\i}ve periods are
\label{expper}
\begin{itemize}
\item 
  $\displaystyle \int_1^2 \frac{dt}{t}= \ln 2$,
\item 
  $\displaystyle \int_{\{x^2+y^2\, \le \,1\}} dx\,dy = \pi$ and
\item 
  $\displaystyle 
  \int_G \frac{dt}{s} = 
  \int_1^2 \frac{dt}{\sqrt{t^3+1}} = 
  \text{elliptic integrals}$, \\
  \hbox{\hspace{3.5cm} for}
  $ G:=\{(t,s)\in\R^2\st 1\le t\le2,\,0\le s,\,s^2=t^3+1\}$.
\end{itemize}

As a problematic example, we consider the following identity.

\begin{proposition}[{cf. \cite[p.~62]{kontsevich}}]
\label{zetaprop}
We have
\index{zeta-value@$\zeta$-value}
\begin{equation}
\label{zeta}
\int_{0\,\le\,t_1\,\le\,t_2\,\le\,1} \frac{dt_1 \wedge dt_2}{(1-t_1)\,
  t_2} = \zeta(2).
\end{equation}
\end{proposition}

\begin{proof}
This equality follows by a simple power series manipulation:
For $0\le t_2<1$, we have
$$
\int_0^{t_2} \frac{dt_1}{1-t_1} = - \log (1-t_2) = 
\sum_{n=1}^{\infty} \frac{t_2^n}{n}.
$$
Let $\epsilon>0$. The power series $\sum_{n=1}^{\infty}
\frac{t_2^{n-1}}{n}$ converges uniformly for $0 \le t_2 \le 1-\epsilon$ and
we get
$$
\int_{0\,\le\,t_1\,\le\,t_2\,\le\,1-\epsilon} 
\frac{dt_1\,dt_2}{(1-t_1)\,t_2} = 
\int_0^{1-\epsilon} \sum_{n=1}^{\infty} \frac{t_2^{n-1}}{n}\,dt_2 =
\sum_{n=1}^{\infty} \frac{(1-\epsilon)^n}{n^2} .
$$
Applying Abel's Theorem
\cite[vol. 2, XII, 438, $6^{\circ}$,
p. 411]{fichtenholz}
at $(\ast)$, 
using $\sum_{n=1}^{\infty} \frac{1}{n^3} < \infty$
gives us
$$
\int_{0\,\le\,t_1\,\le\,t_2\,\le\,1} \frac{dt_1\,dt_2}{(1-t_1)\,t_2} = 
\lim_{\epsilon\rightarrow 0} \sum_{n=1}^{\infty} \frac{(1-\epsilon)^n}{n^2}
\stackrel{(\ast)}{=} 
\sum_{n=1}^{\infty} \frac{1}{n^2} = \zeta(2) .
$$
\end{proof}

Equation (\ref{zeta}) is not a valid representation of $\zeta(2)$ as
an integral for a na{\"\i}ve period in our sense, because the pole locus
$\{t_1=1\}\union\{t_2=0\}$ of $\frac{dt_1\,\wedge\,dt_2}{(1-t_1)\,t_2}$
is not disjoint with the domain of integration $\{0\le t_1 \le t_2 \le
1\}$. 
But (\ref{zeta}) gives a valid period integral according to the
original definition of \cite[p.~772]{kontsevich_zagier} --- see Remark
\ref{kzdef}.
We will show in Example \ref{exp4} how to circumvent this difficulty.

\medskip
We will prove in the Period Theorem \ref{perioperio}, that 
special pairing, abstract and na{\"\i}ve periods
are essentially the same, i.e.
$$
\P^{\,\prime}_p = \P_a = \P_n,
$$
but have to make a digression on the triangulation of varieties first.

Note that in Subsection \ref{per4} we will give a fourth definition of
a period.

\section{Triangulation of Algebraic Varieties}

If \Xo\ is a variety defined over \Q\, we may ask whether any singular
homology class
$\gamma \in \h_{\bullet}^{\text{sing}} (\Xh;\Q)$ can be represented by
an object described by polynomials.
This is indeed the case: For a precise statement we
need several definitions. 
The result will be formulated in Proposition \ref{semialg}.

\subsection{Semi-algebraic Sets}

We already defined \QR-semi-algebraic sets in Definition \ref{semialgset}.

\begin{definition}[\QR-semi-algebraic map 
  {\cite[p. 168]{hironaka_74}}]
  \index{Semi-algebraic!map}
  A continuous map $f$ between \QR-semi-algebraic sets 
  $A\subseteq \R^n$ and $B\subseteq \R^m$ is said to be
  {\em \QR-semi-algebraic} if its graph 
  $$ 
  \Gamma_f := \bigl\{\bigl(a,f(a)\bigr) \st a\in A\bigr\} \subseteq A \times B
  \subseteq \R^{n+m}
  $$ 
  is \QR-semi-algebraic.
\end{definition}

\begin{example}
\label{semialgmapexp}
  Any polynomial map
  \begin{align*}
  f:A & \longrightarrow B \\
  (a_1,\ldots,a_n) & \mapsto
  (f_1(a_1,\ldots,a_n),\ldots,f_m(a_1,\ldots,a_n))
  \end{align*}
  between \QR-semi-algebraic sets $A\subseteq \R^n$ and $B\subseteq \R^m$
  with
  $ f_i \in \QR[x_1,\ldots,x_n] $ for $i=1,\ldots,m$
  is \QR-semi-algebraic, since it is continuous and its graph
  $\Gamma_f \subseteq \R^{n+m}$ is cut out from $A\times B$ by the
  polynomials
  \begin{equation}
  \label{polmapgraph}
  y_i - f_i(x_1,\ldots,x_n) \in \QR[x_1,\ldots,x_n,y_1,\ldots,y_m]
  \quad\text{for}\quad i=1,\ldots,m .
  \end{equation}
  We can even allow $f$ to be a rational map with rational component
  functions
  $$
  f_i \in \QR(x_1,\ldots,x_n),\quad i=1,\ldots,m
  $$
  as long as none of the denominators of the $f_i$ vanish at a point
  of $A$. The argument remains the same except that the expression
  (\ref{polmapgraph}) has to be multiplied by the denominator of $f_i$.
\end{example}  

\begin{fact}[{\cite[Prop. II, p. 167]{hironaka_74}},
  {\cite[Thm. 3, p. 370]{seidenberg}}]
  \label{seidenberg_tarski}
  By a result of Seidenberg-Tarski, the image
  (respectively preimage) of a \QR-semi-algebraic set under a
  \QR-semi-algebraic map is again \QR-semi-algebraic.
\end{fact}

As the name suggests any algebraic set should be in particular
\QR-semi-algebraic.

\begin{lemma}
  \label{algsemialg}
  Let $X_0$ be an algebraic variety defined over \QR. Then we can
  regard the complex analytic space \Xh 
  (cf. Subsection \ref{defcomplexanspace})
  associated to the base change $X=X_0\times_{\QR}\C$ 
  as a bounded \QR-semi-algebraic subset 
  \begin{equation}
  \label{semialgincl}
  \Xh \subseteq \R^N 
  \end{equation}
  for some $N$. Moreover, if 
  $
  f_0: X_0 \rightarrow Y_0
  $
  is a morphism  of varieties defined over \QR, we can
  consider $f_{\rm an}: \Xh \rightarrow \Yh$ as a \QR-semi-algebraic map.
\end{lemma}

\pagebreak

\begin{proof} 
{\em First step $X_0=\Q P^n$:}
Consider
\begin{itemize}
\item
  $\C \Ph^n$ with homogenous coordinates $x_0,\ldots,x_n$, which we split
  as $x_m=a_m + ib_m$ with $a_m, b_m \in \R$ in real and imaginary part, and
\item
  $\R^N$, $N=2(n+1)^2$, with coordinates
  $\{y_{kl},z_{kl}\}_{k,l=0,\ldots,n}$.
\end{itemize}
We define a map
\begin{align*}
\psi : 
  \underset{[x_0:\ldots:x_n]}{\C \Ph^n}
& \underset{\strut}{\longrightarrow}
  \underset{(y_{00},z_{00},\ldots,y_{nn},z_{nn})}{\R^N}\\[3mm]
[x_0 : \ldots : x_n] 
& \mapsto
\biggl( \ldots,
  \myunderbrace{\ts y_{kl}}{
    \frac{\Re x_k \ol{x}_l}{\sum_{m=0}^n |x_m|^2}},
  \myunderbrace{\ts z_{kl}}{
    \frac{\Im x_k \ol{x}_l}{\sum_{m=0}^n |x_m|^2}}, \ldots
  \biggr) \\[3mm]
[a_0 + i b_0 : \ldots : a_n + i b_n] 
& \mapsto 
\biggl( \ldots, 
  \myunderbrace{\ts y_{kl}}{
    \frac{a_k a_l + b_k b_l}{\sum_{m=0}^n a_m^2 + b_m^2}},
  \myunderbrace{\ts z_{kl}}{
    \frac{b_k a_l - a_k b_l}{\sum_{m=0}^n a_m^2 + b_m^2}}, \ldots
  \biggr).
\end{align*}
Rewriting the last line (with the convention $0\cdot
\cos(\myatop{\text{indeterminate}}{\text{angle}})=0$) as
\begin{equation}
\label{psimodarg}
[r_0 e^{i \phi_0}:\ldots:r_n e^{i \phi_n}] \mapsto 
\left( \ldots, \frac{r_k r_l \cos(\phi_k - \phi_l)}{\sum_{m=0}^n r_m^2}, 
  \frac{r_k r_l \sin(\phi_k - \phi_l)}{\sum_{m=0}^n
    r_m^2}, \ldots \right)
\end{equation}
shows that $\psi$ is injective: Assume
$$
\psi\bigl([r_0 e^{i \phi_0}:\ldots:r_n e^{i \phi_n}]\bigl) =
(y_{00},z_{00},\ldots,y_{nn},z_{nn})
$$
where $r_k \ne 0$, or equivalently $y_{kk} \ne 0$, for a fixed $k$.
We find
\begin{align*}
& \frac{r_l}{r_k} = \frac{\sqrt{ y_{kl}^2 + z_{kl}^2}}{y_{kk}}, \quad
\text{and} \\
& \phi_k - \phi_l = 
\begin{cases}
\arctan ( z_{kl} / y_{kl} ) & \text{if } y_{kl} \ne 0, \\
\pi / 2 & \text{if } y_{kl}=0, z_{kl}>0, \\
\text{indeterminate} & \text{if } y_{kl}=z_{kl}=0, \\
- \pi / 2 & \text{if } y_{kl}=0, z_{kl}<0; 
\end{cases}
\end{align*}
that is the preimage of $(y_{00},z_{00},\ldots,y_{nn},z_{nn})$ is
uniquely determined.

Therefore we can consider $\C \Ph^n$ via $\psi$ as a subset of
$\R^N$. 
It is bounded since it is contained in the unit sphere
$S^{N-1}\subset \R^N$.
We claim that $\psi(\C \Ph^n)$ is also \QR-semi-algebraic. 
The composition of the projection
\begin{align*}
\pi : \R^{2(n+1)}\setminus\{(0,\ldots,0)\} & \longrightarrow \C \Ph^n \\
(a_0,b_0,\ldots,a_n,b_n) & \mapsto [a_0 + i b_0 : \ldots : a_n + i
b_n]
\end{align*} 
with the map $\psi$ is a polynomial map, hence
\QR-semi-algebraic by Example \ref{semialgmapexp}. Thus
$$
\im \psi \circ \pi = \im \psi \subseteq \R^N
$$
is \QR-semi-algebraic by Fact \ref{seidenberg_tarski}.

{\em Second step (zero set of a polynomial):}
We use the notation
\begin{align*}
& V(g):=\{ \ul{x} \in \C \Ph^n \st g(\ul{x}) =0 \}\quad
\text{for\quad $g\in\C[x_0,\ldots,x_n]$ homogenous},\quad\text{and} \\
& W(h):=\{\ul{t}\in \R^N \st h(\ul{t})=0 \}\quad
\text{for\quad $h\in\C[y_{00},\ldots,z_{nn}]$}.
\end{align*}
Let $\Xh=V(g)$ for some homogenous $g\in \QR[x_0,\ldots,x_n]$. 
Then $\psi(\Xh)\subseteq \R^N$ is a \QR-semi-algebraic subset, as a little
calculation shows. Setting for $k=0,\ldots,n$
\begin{align*}
g_k & :=\text{``}g(\ul{x}\, \ol{x}_k)\text{''} \\
& = g(x_0 \ol{x}_k, \ldots, x_n \ol{x}_k) \\
& = g\bigl((a_0 a_k + b_0 b_k) +  i(b_0 a_k - a_0 b_k),\ldots,
      (a_n a_k + b_n b_k) +  i(b_n a_k - a_n b_k)\bigr), \\
\intertext{where $x_j=a_j + i b_j$ for $j=0,\ldots,n$, and}
h_k & :=g(y_{0k}+i z_{0k},\ldots,y_{nk}+i z_{nk}),
\end{align*}
we obtain
\begin{align*}
\psi(\Xh) & = \psi(V(g))\\
& = \bigcap_{k=0}^n \psi(V(g_k)) \\
& = \bigcap_{k=0}^n \psi(\C \Ph^n) \cap W(h_k) \\
& = \bigcap_{k=0}^n \psi(\C \Ph^n) \cap W(\Re h_k) \cap W(\Im h_k) .
\end{align*}

{\em Final step:}\/ 
We can choose an embedding
$$
X_0 \subseteq \QR P^n,
$$
thus getting
$$
\Xh \subseteq \C \Ph^n .
$$
Since $X_0$ is a locally closed subvariety of $\QR P^n$, $\Xh$ can be
expressed in terms of subvarieties of the form $V(g)$ with
$g\in \QR[x_0,\ldots,x_n]$, 
using only the following basic operations
\begin{enumroman}
\item complementary set,
\item finite intersection,
\item finite union.
\end{enumroman}
Now \QR-semi-algebraic sets are stable under these operations as well and
the first assertion is proved. {\hfill $\blacksquare$}

{\em Second assertion:}\/
The first part of the lemma provides us with \QR-semi-algebraic
inclusions
\begin{align*}
\psi : 
  \Xh \subseteq 
  \underset{\ul{x}=[x_0:\ldots:x_n]}{\C \Ph^n} & \subseteq 
  \underset{(y_{00},z_{00},\ldots,y_{nn},z_{nn})}{\R^N}, \\[1mm]
\phi: 
  \Yh \subseteq 
  \underset{\ul{u}=[u_0:\ldots:u_m]}{\C \Ph^m} & \subseteq 
  \underset{(v_{00},w_{00},\ldots,v_{mm},w_{mm})}{\R^M},
\end{align*}
and a choice of coordinates as indicated. We use the notation
\begin{align*}
V(g) &:=
  \{ (\ul{x},\ul{u}) \in \C \Ph^n \times \C \Ph^m \st
  g(\ul{x},\ul{u})=0 \},\\
& \quad \text{for 
  $g\in \C[x_0,\ldots,x_n,u_0,\ldots,u_m]$ homogenous in both
  $\ul{x}$ and $\ul{u}$},\quad\text{and} \\
W(h) &:=\{\ul{t}\in \R^{N+M} \st h(\ul{t})=0 \},\quad
\text{for $h\in \C[y_{00},\ldots,z_{nn},v_{00},\ldots,w_{mm}]$}.
\end{align*}
Let $\{U_i\}$ be a finite open affine covering of $X_0$ such that
$f_0(U_i)$ satisfies
\begin{itemize}
\item
$f_0(U_i)$ does not meet the hyperplane 
$\{u_j=0\}\subset\QR P^m$
for some $j$%
, and
\item
$f_0(U_i)$ is contained in an open affine subset $V_i$ of $Y_0$.
\end{itemize}
This is always possible, since we can start with the open covering
$Y_0\cap \{u_j\ne 0\}$ of $Y_0$, take a subordinated open affine covering $\{V_{i'}\}$,
and then choose a finite open affine covering $\{U_i\}$ subordinated
to $\{f^{-1}(V_{i'})\}$.
Now each of the maps
$$
f_i := {f_{\rm an}}_{|U_i} : \Uh_i \longrightarrow \Yh
$$
has image contained in $\Vh_i$ and does not meet the hyperplane
$\{\ul{u}\in\C \Ph^m\st u_j=0\}$ for an appropriate $j$
$$
f_i : \Uh_i \longrightarrow \Vh_i .
$$
Being associated to an algebraic map between affine varieties,
this map is rational
$$
f_i : \ul{x} \mapsto
\left[
  \frac{g'_0(\ul{x})}{g''_0(\ul{x})}:\ldots:
  \frac{g'_{j-1}(\ul{x})}{g''_{j-1}(\ul{x})}:
  \underset{j}{1}:
  \frac{g'_{j+1}(\ul{x})}{g''_{j+1}(\ul{x})}:\ldots:
  \frac{g'_m(\ul{x})}{g''_m(\ul{x})}
\right],
$$
with
$g'_k,g''_k \in \QR[x_0,\ldots,x_n]$, $k=0,\ldots,\widehat{j},\ldots,m$.
Since the graph $\Gamma_{f_{\rm an}}$ of ${f_{\rm an}}$ is the finite union of the graphs
$\Gamma_{f_i}$ of the $f_i$, it is sufficient to prove that $(\psi
\times \phi)(\Gamma _{f_i})$ is a \QR-semi-algebraic subset of
$\R^{N+M}$. Now
$$
\Gamma_{f_i} = 
\left(\Uh_i \times \Vh_i \right) \cap 
\bigcap_{\myatop{k=0}{k\ne j}}^n 
  V\left(\frac{y_k}{y_j} - \frac{g'_k(\ul{x})}{g''_k(\ul{x})}\right) =
\left(\Uh_i \times \Vh_i \right) \cap 
\bigcap_{\myatop{k=0}{k\ne j}}^n 
V(y_k g''_k(\ul{x}) - y_j g'_k(\ul{x})),
$$
so all we have to deal with is 
$$
V(y_k g''_k(\ul{x}) - y_j g'_k(\ul{x})).
$$
Again a little calculation is necessary. Setting
\begin{align*}
g_{pq} &:= \text{``}u_k \ol{u}_q g''_k(\ul{x}\,\ol{x}_p) -
  u_j \ol{u}_q g'_k(\ul{x}\,\ol{x}_p)\text{''} \\
& = u_k \ol{u}_q g''_k(x_0 \ol{x}_p,\ldots,x_n \ol{x}_p) -
  u_j \ol{u}_q g'_k(x_0 \ol{x}_p,\ldots,x_n \ol{x}_p) \\
& = \bigl((c_k c_q + d_k d_q) +i(d_k c_q - c_k d_q)\bigr) \\
&\hspace{1.5cm}  g''_k\bigl((a_0 a_p+b_0 b_p)+ i(b_0 a_p- a_0 b_p),\ldots, 
    (a_n a_p+b_n b_p)+ i(b_n a_p- a_n b_p)\bigr) \\
&\quad  -  \bigl((c_j c_q + d_j d_q) +i(d_j c_q - c_j d_q)\bigr) \\
&\hspace{1.5cm}  g'_k\bigl((a_0 a_p+b_0 b_p)+ i(b_0 a_p- a_0 b_p),\ldots, 
    (a_n a_p+b_n b_p)+ i(b_n a_p- a_n b_p)\bigr) ,
\end{align*}
where 
$x_l=a_l+ib_l$ for $l=0,\ldots,n$,
$u_l=c_l+id_l$ for $l=0,\ldots,m$, and
$$
h_{pq}:=(v_{kq}+iw_{kq})g''_k(y_{0p}+iz_{0p},\ldots,y_{np}+iz_{np}) -
  (v_{jq}+iw_{jq})g'_k(y_{0p}+iz_{0p},\ldots,y_{np}+iz_{np}),
$$
we obtain
\begin{align*}
(\psi\times\phi)\,\Bigl(V\bigl(y_k g''_k(\ul{x}) - & y_j g'_k(\ul{x})\bigr)\Bigr) =  \\
= & \bigcap_{p=0}^n \bigcap_{q=0}^m (\psi\times\phi)(V(g_{pq})) \\
= & \bigcap_{p=0}^n \bigcap_{q=0}^m 
  (\psi\times\phi)(\Uh_i\times\Vh_j) \cap W(h_{pq}) \\
= & \bigcap_{p=0}^n \bigcap_{q=0}^m 
  (\psi\times\phi)(\Uh_i\times\Vh_j) \cap W(\Re h_{pq}) \cap W(\Im h_{pq}). 
\end{align*}
\end{proof}

\subsection{Semi-algebraic Singular Chains}

We need further prerequisites 
in order to state the announced Proposition \ref{semialg}.

\begin{definition}[Open simplex, 
  {\cite[p. 168]{hironaka_74}}]
  By an {\em open simplex} $\osimplex$
\index{Simplex!open}
\notation{D}{$\osimplex$}{interior of a simplex or a point}
  we mean the interior of a simplex
  ($=$ the convex hull of\/ $r+1$ points in $\R^n$ which span a
  $r$-dimensional subspace). For convenience, a point is considered  as
  an open simplex as well. 
 
  The notation $\Triangle_d$ 
  \notation{D}{$\Triangle_p$}{standard $p$-simplex}
  will be reserved for the {\em closed standard simplex}
  spanned by the standard basis\/
  $\{(0,\ldots,0,\underset{i}{1},0,\ldots,0)\st i=1,\ldots,d+1\}$ 
  of\/ $\R^{d+1}$.
\end{definition}

Consider the following data ($\ast$):
\begin{itemize}

\item \Xo\ a variety defined over \QR,

\item \Do\ a divisor on \Xo\ with normal crossings,

\item and finally $\gamma \in \h_p^{\rm sing} (\Xh,\Dh;\Q)$, $\,p\in\N_0$.

\end{itemize}

As usual, we have denoted by \Xh\ (resp. \Dh) the complex analytic space
associated to the base change $X=X_0\times_{\QR}\C$ (resp. $D=D_0\times_{\QR}\C$). 

By Lemma \ref{algsemialg},
we may consider both \Xh\ and \Dh\ as bounded \QR-semi-algebraic
subsets of $\R^N$.

We are now able to formulate our proposition.

\begin{proposition}
  \label{semialg}
  With data {\rm ($\ast$)} as above,
  we can find a representative of $\gamma$ that is a rational linear 
  combination of singular
  simplices each of which is \QR-semi-algebraic.
\end{proposition}

The proof of this proposition relies on the following proposition due to Lojasiewicz which has
been written down by Hironaka.

\begin{proposition}[Triangulation of \QR-semi-algebraic sets, {\cite[p.~170]{hironaka_74}}]
  \label{triangulation}
  For $\{ X_{i} \}$ a finite system of bounded \QR-semi-algebraic
  sets in $\R^n$, there exists a simplicial decomposition
  $$
  \R^n = \coprod_{j} \osimplex_{j} 
  $$
  by open simplices $\osimplex_j$ and a \QR-semi-algebraic automorphism 
  $$
  \kappa : \R^n \to \R^n
  $$
  such that each
  $X_{i}$ is a finite union of some of the
   $\kappa ( \osimplex
  _{j} )$.
\end{proposition}

\begin{note} Although Hironaka considers \R-semi-algebraic sets,
  it has been checked by the author, that we
  can safely replace \R\ by \QR\ in his article (including the fact he
  cites from \cite{seidenberg}). The only problem that could possibly arise
  concerns a ``good direction lemma'':

  \begin{lemma}[Good direction lemma for \R,
  {\cite[p. 172]{hironaka_74}}, {\cite[Thm. 5.I, p. 242]{koopman_brown}}]
    Let $Z$ be a \R-semi-algebraic subset of $\R^n$, which is nowhere
    dense.
    A direction $v\in |\R P^{n-1}|$ is called {\em good}, if any line
    $l$ in $\R^n$ parallel to $v$ meets $Z$ in a discrete (maybe
    empty) set of points; otherwise $v$ is called {\em bad}.
    Then the set $B(Z)$ of bad directions
    is a Baire category set in $|\R P^{n-1}|$.
  \end{lemma}
  
  This gives immediately good directions 
  $v \in |\R P^{n-1}| \setminus B(Z)$, but not
  necessarily  $v\in |\QR P^{n-1}| \setminus B(Z)$. 
  However, in Remark 2.1 of
  \cite{hironaka_74}, which follows directly after the lemma, the
  following statement is made: If $Z$ is compact, then $B(Z)$ is
  closed in $|\R P^{n-1}|$. 
  In particular $|\QR P^{n-1}| \setminus B(Z)$ 
  will be non-empty. 
  Since we only consider {\bf bounded} \QR-semi-algebraic
  sets $Z'$, we may take $Z:=\ol{Z'}$ (which is compact by
  Heine-Borel), and thus find a good direction 
  $v\in |\QR P^{n-1}| \setminus B(Z')$ using $B(Z') \subseteq B(Z)$. 

  \begin{lemma}[Good direction lemma for \QR]
    Let $Z'$ be a bounded \QR-semi-algebraic subset of $\R^n$, which is
    nowhere dense. 
    Then the set $|\QR P^{n-1}| \setminus B(Z)$ of good directions is
    non-empty.
  \end{lemma}
\end{note}

\newcommand{\trial}{\ensuremath{\osimplex_{j}}}
\newcommand{\vartrial}{\ensuremath{\varosimplex_{j}}}
\newcommand{\oltrial}{\ensuremath{\csimplex_{j}}}

{\em Proof of Proposition \ref{semialg}.}
Applying Proposition \ref{triangulation} to the two-element system of
\QR-semi-algebraic sets
$\Xh, \Dh \subseteq \R^N$, we obtain a \QR-semi-algebraic
decomposition
$$
\R^N = \coprod_{j} \trial
$$
of $\R^N$ by open simplices $\trial$ and a \QR-semi-algebraic automorphism
$$
\kappa: \R^N \to \R^N .
$$
We write $\oltrial$ for the closure of $\trial$. The sets
$$
K:= \{ \trial \st \kappa(\trial) \subseteq \Xh \} \text{\quad and \quad}
L:= \{ \trial \st \kappa(\trial) \subseteq \Dh \}
$$
can be thought of as finite simplicial complexes, but built out of
open simplices instead of closed ones. We define their {\em geometric
realizations}
$$
\vert K \vert := \Union_{\vartrial \in K} \trial \text{\quad and \quad}\\
\vert L \vert := \Union_{\vartrial \in L} \trial.
$$
Then Proposition \ref{triangulation} states that $\kappa$ maps the pair
of topological spaces $(\vert K \vert, \vert L \vert)$
homeomorphically to $(\Xh,\Dh)$. 

{\em Easy case:} If \Xo\ is complete, so is $X$ 
(by \cite[Cor. II.4.8(c), p. 102]{hartshorne}), hence \Xh\ and \Dh\ will be compact
\cite[B.1, p. 439]{hartshorne}. 
In this situation
$$
\ol{K} := \{ \oltrial \st \kappa(\oltrial) \subseteq \Xh \} 
\text{\quad and \quad}
\ol{L} := \{ \oltrial \st \kappa(\oltrial) \subseteq \Dh \}
$$
are (ordinary) simplicial complexes, whose geometric realizations
coincide with those of $K$ and $L$, respectively. 
In particular 
\begin{equation}
\label{simplsingiso}
\begin{aligned}
\h_{\bullet}^{\rm simpl}(\ol{K},\ol{L};\Q) & \iso
\h_{\bullet}^{\rm sing}(\left| {\ol{K}} \right|,\left| \ol{L} \right| ;\Q) \\
& \iso \h_{\bullet}^{\rm sing}(\vert K \vert,\vert L \vert;\Q) \\
& \iso \h_{\bullet}^{\rm sing}(\Xh,\Dh;\Q).
\end{aligned}
\end{equation}
Here $\h_{\Bul}^{\rm simpl}(\ol{K},\ol{L};\Q)$ denotes simplicial
homology of course.

We write $\gamma_{\rm simpl}\in \h_p^{\rm simpl}(\ol{K},\ol{L};\Q)$ and
$\gamma_{\rm sing}\in \h_p^{\rm sing}(\left| \ol{K} \right|,\left| \ol{L}
\right|;\Q)$
for the image of $\gamma$ under this isomorphism. 
Any representative
$\Gamma_{\rm simpl}$ of $\gamma_{\rm simpl}$ is a rational linear
combination
$$
\textstyle
\Gamma_{\rm simpl} = \sum_{j} a_{j} \, \csimplex_{j}, 
\quad
a_{j} \in \Q
$$
of oriented closed
simplices
$\oltrial \in \ol{K}$.
We can choose orientation-preserving affine-linear maps of the
standard simplex $\Triangle_p$ to $\triangle_j$
$$
\sigma_{j} : \Triangle_p \longrightarrow \oltrial
\quad\text{for}\quad
\oltrial \in \Gamma_{\rm simpl}.
$$
These maps yield a representative 
$$
\textstyle
\Gamma_{\rm sing} := \sum_{j} a_{j} \, \sigma_{j} 
$$ 
of
$\gamma_{\rm sing}$. Composing with $\kappa$ yields
$\Gamma:=\kappa_{\ast}\Gamma_{\rm sing} \in \gamma$, where $\Gamma$ has the
desired properties. 

In the {\bf general case}, we perform a barycentric subdivision \BB\
on $K$ twice (once is not enough) and define $|K|$ and $|L|$ not as the ``closure''  of $K$ and
$L$, but as follows (see Figure \ref{triagprooffig1})
\begin{equation}
\label{Kbardef}
\begin{aligned}
& \ol{K} := \{ \csimplex \st \osimplex \in \BB^2(K) \text{ and }
\csimplex \subseteq |K|\}, \\
& \ol{L} := \{ \csimplex \st \osimplex \in \BB^2(K) \text{ and }
\csimplex \subseteq |L|\}.
\end{aligned}
\end{equation}

\ignore{
\begin{figure}[htb]
\rule{14.3cm}{0.1mm}
\vspace*{2mm}

\begin{tabular}{ccc}
\epsfig{file=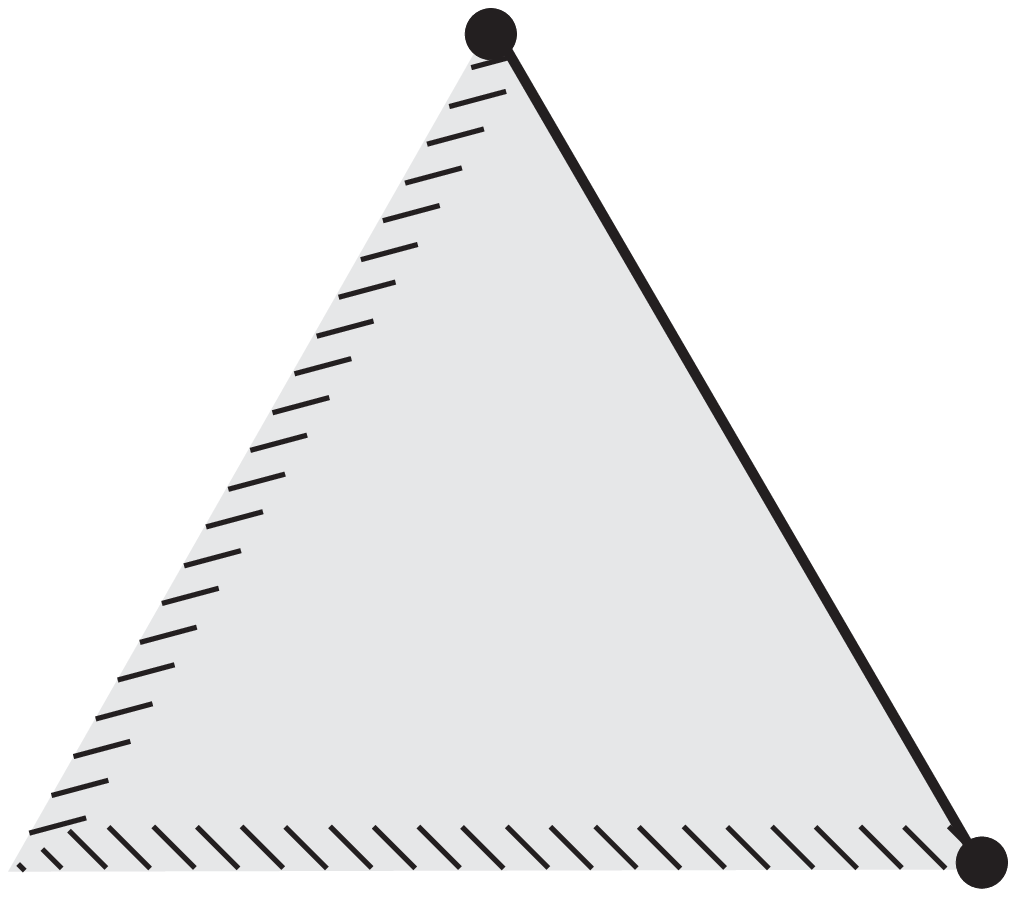,height=3.5 cm} &
\epsfig{file=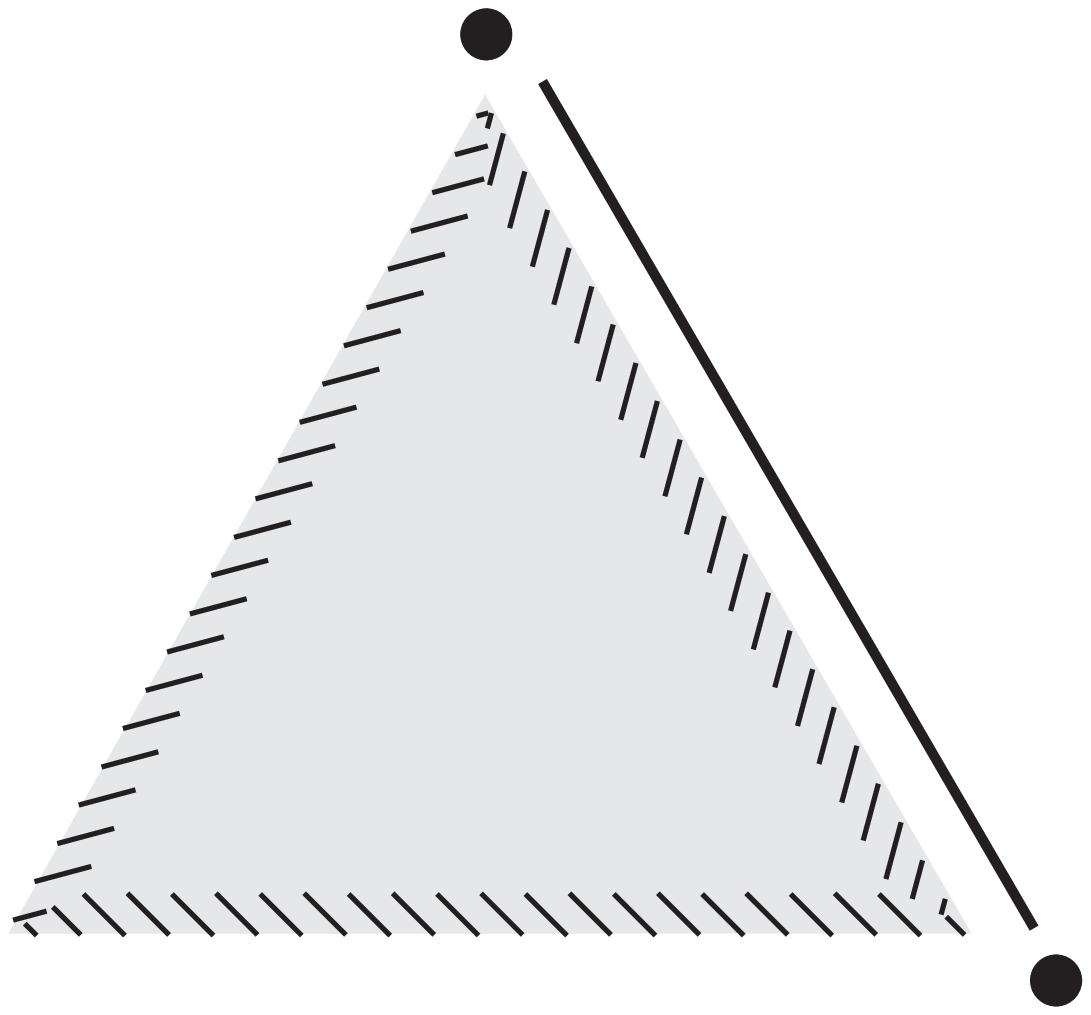,height=3.5 cm} &
\epsfig{file=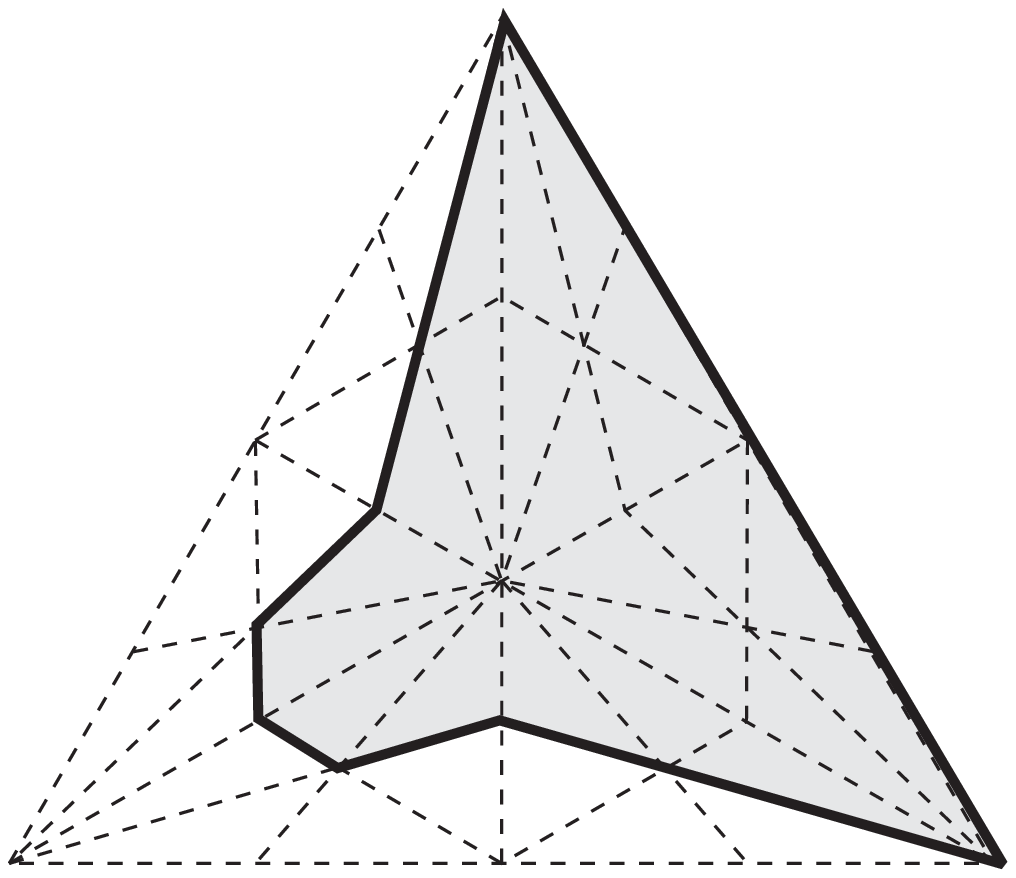,height=3.5 cm} \\

$\kappa^{-1}(\Xh) \cap \csimplex_{j}$ &
$K \cap \csimplex_{j}$ &
$\left| \overline{K} \right| \cap \csimplex_{j}$ \\[0.5mm]

\parbox{4.5cm}{
Intersection of $\kappa^{-1}(\Xh)$ with a closed
$2$-simplex $\csimplex_j$, where we assume that
part of the boundary $\del \csimplex_j$ does not 
belong to $\kappa^{-1}(\Xh)$} 
&
\parbox{4.3cm}{
Open simplices of $K$ contained in $\csimplex_j$}
& 
\parbox{4.3cm}{
Intersection of $\left|\ol{K}\right|$ with $\csimplex_j$
(the dashed lines show the barycentric subdivision)}

\end{tabular}
\caption{Definition of $\ol{K}$}
\label{triagprooffig1}
\rule{14.3cm}{0.1mm}
\end{figure}
}

The point is that the pair of topological spaces $(\left| \ol{K} \right|,
\left| \ol{L} \right|)$ is a strong deformation retract of $(\vert K
\vert, \vert L \vert)$. Assuming this we see that in the general case
with $\ol{K}$, $\ol{L}$ defined as in (\ref{Kbardef}), the isomorphism 
(\ref{simplsingiso}) still holds and we can proceed as in the easy case
to prove the proposition. 

We 
define the retraction map
$$
\rho : (|K|\times [0,1],|L|\times [0,1]) \rightarrow
(\left|\ol{K}\right|,\left|\ol{L}\right|)
$$
as follows: Let $\trial \in K$ be an open simplex which is not
contained in the boundary of any other simplex of $K$ and set
$$
inner:=\oltrial \cap \ol{K}, \qquad outer:=\oltrial \setminus \ol{K}.
$$

\begin{wrapfigure}{r}{4.5cm}
\epsfig{file=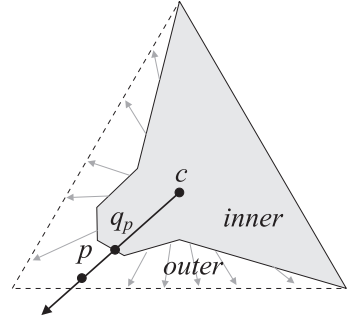,width=4.5 cm}
\caption{Definition of $q_p$}
\label{triagprooffig2}
\end{wrapfigure}
Note that $inner$ is closed. 
For any point $p\in outer$ the ray 
$\overrightarrow{c\,p}$ from the center $c$ of \trial\ through $p$
``leaves'' the set $inner$ at a point $q_p$, i.e. $\overrightarrow{c\, p} \cap
inner$ equals the line segment $c\, q_p$; see Figure \ref{triagprooffig2}. 
The map
\begin{align*}
\rho_{j} : \oltrial \times [0,1] & \rightarrow \oltrial\\
(p,t) & \mapsto
  {\begin{cases}
    p \text{\quad if $p\in inner$}, \\
    q_p +t\cdot (p-q_p) \text{\quad if\quad $p\in outer$}
  \end{cases}}
\end{align*}
retracts $\oltrial$ onto $inner$.

Now these maps $\rho_{j}$ glue together
to give the desired homotopy $\rho$.
\hfill $\Box$

\vspace{3mm}
We want to state one of the intermediate results of this proof
explicitly.

\begin{corollary}
Let $X_0$ and $D_0$ be as above. Then the pair of topological spaces 
$(\Xh,\Dh)$ is homotopic equivalent to a pair of (realizations of) simplicial complexes 
$(\vert X^{\rm simpl}\vert,\vert D^{\rm simpl}\vert)$.
\end{corollary}

\section{Periods Revisited}

\subsection{Comparison of Definitions of Periods}

We are now ready to prove the theorem announced in the introduction.

\begin{maintheorem}[Period theorem]
\label{perioperio}
It holds
$$
\P^{\,\prime}_p=\P_a=\P_n,
$$
that is the following three definitions of a period given in Section 
\ref{per123} coincide
\begin{itemize}
\item special pairing periods (cf. Definition \ref{per1}),
\item abstract periods (cf. Definition \ref{per2}), and
\item na{\"\i}ve periods (cf. Definition \ref{per3}).
\end{itemize}
\end{maintheorem}

Note that we will give a fourth definition of a period 
in Subsection \ref{per4}.

\begin{proof}
{$\P_a\subseteq\P^{\,\prime}_p$:}
Let $\int_{\gamma}\omega\in\P_a$ be an abstract period coming from
$(X_0,D_0,\omega_0,\gamma)$. 
We want to understand $\omega_0$ as an element of
$\Gamma(X_0;\wt{\Omega}^d_{X_0,D_0/\Q})=\Gamma(X_0;\Omega^d_{X_0/\Q})\oplus\Gamma(D_0;\wt{\Omega}^{d-1}_{D_0/\Q})$ 
via
$\Gamma(X_0;\Omega^d_{X_0/\Q})\injection\Gamma(X_0;\wt{\Omega}^d_{X_0,D_0/\Q})$,
where $d:={\rm dim\,}X_0$.
Then $\omega_0\in\Gamma(X_0;\wom{X_0,D_0/\Q})$ is a
cocycle, since $\omega_0$ is closed (for dimension reasons) and its
restriction to each of the irreducible components $D_j$ of $D_0$ is zero (also for dimension reasons).
Thus we can consider $[\omega_0]$ via
$$
\hh\bul\Gamma(X_0;\wom{X_0,D_0/\Q})\to\hDR{\Bul}(X_0,D_0/\Q)
$$
as an algebraic \deRham cohomology class and find 
(cf. Subsection \ref{alteriso})
$$
\int_{\gamma}\omega=\langle\gamma,[\omega_0]\rangle\in\P^{\,\prime}_p.
$$
{\hfill $\blacksquare$}

{$\P^{\,\prime}_p\subseteq\P_n$:}
As explained in Subsection \ref{alteriso}, a special pairing period
$$
\langle\gamma,\omega_0\rangle\in\P^{\,\prime}_p
$$
is a sum of integrals 
$$
\textstyle
\sum_I\langle\gamma_I,\omega_I\rangle=\sum_I\int_{\gamma_I}\omega_I,
$$
where $\omega_I\in\hDR{\Bul}(D_I/\Q)$ and
$\gamma_I\in\h\bul_{\infty}(D_I;\Q)$.
Since $\P_n$ is closed under addition by Proposition \ref{palgebra},
it suffices to prove 
$\int_{\gamma_I}\omega_I\in\P_n$ for all index sets $I$.

We choose an embedding
$$
D_I \subseteq \underset{(x_0:\ldots:x_n)}{\Q P^n}
$$ 
and equip $\Q P^n$ with coordinates as indicated. 
Lemma \ref{algsemialg} provides us with a map
$$
\psi : \C \Ph^n \inclusion \R^N
$$
such that $\Dh_I$ and $\C \Ph^n$ become \QR-semi-algebraic subsets of
$\R^N$. Then, by Proposition \ref{semialg}, $\psi_{\ast}\gamma_I$ has a
representative which is a rational linear combination of singular
simplices $\Gamma_i$, each of which is \QR-semi-algebraic.

We set $G_i:=\im \Gamma_i$ and are left to prove
$$
\int_{\psi^{-1}(G_i)}\omega_I\in\P_n.
$$
This will be clear as soon as we find a rational differential form
$\omega'_I$ on $\R^N$ such that $\psi^{\ast}\omega'_I=\omega_I$, since then
$$
\int_{\psi^{-1}(G_i)}\omega_I = 
\int_{\psi^{-1}(G_i)}\psi^{\ast}\omega'_I =
\int_{G_i} \omega'_I \in \P_n.
$$
After eventually applying a barycentric subdivision to $\Gamma_i$, we
may assume w.l.o.g. that there exists a hyperplane in $\C \Ph^n$, 
say $\{x_0=0\}$,
which does not meet $\psi^{-1}(G_i)$.
Furthermore, we may assume that $\psi^{-1}(G_i)$ lies entirely in
$\Uh$ for $U_0$ an
open affine subset of $D_I \cap \{x_0\ne 0\}$. 
(As usual, $\Uh$ denotes the complex analytic space associated to the
base change to \C\ of $U$.)
The restriction of
$\omega_I$ to the open affine subset can be represented 
in the form
(cf. \cite[II.8.4A, II.8.2.1, II.8.2A]{hartshorne})
$$
\sum_{|J|=d} f_J (x_0,\ldots,x_n)\,
d\!\left(\frac{x_{j_1}}{x_0}\right)\wedge\ldots\wedge d\!\left(\frac{x_{j_d}}{x_0}\right)
$$
with $f_J (x_1,\cdots,x_n)\in \Q(x_0,\cdots,x_n)$
being homogenous of degree zero. This expression defines a rational
differential form on all of $\Q P^n$ with coefficients in \Q\ and it
does not have poles on $\psi^{-1}(G_i)$.

We construct the rational differential form $\omega'_I$ on $\R^N$ with
coefficients in $\Q(i)$ as follows
$$
\omega'_I:=\sum_{|J|=d}
f_J\!\left(1,\frac{y_{10}+i z_{10}}{y_{00}+iz_{00}},\cdots,
           \frac{y_{n0}+i z_{n0}}{y_{00}+iz_{00}}\right)
\,
d\!\left(\frac{y_{{j_1}0}+i z_{{j_1}0}}{y_{00}+i z_{00}}\right)
\,\wedge\,\ldots\,\wedge\,
d\!\left(\frac{y_{{j_d}0}+i z_{{j_d}0}}{y_{00}+i z_{00}}\right),
$$
where we have used the notation from the proof of Lemma
\ref{algsemialg}.
Using the explicit form of $\psi$ given in this proof, we obtain
\begin{align*}
\psi^{\ast}
f_J\!\left(1,\frac{y_{10}+i z_{10}}{y_{00}+iz_{00}},\cdots,
           \frac{y_{n0}+i z_{n0}}{y_{00}+iz_{00}}\right) 
& =
f_J\!\left(\frac{x_0 \ol{x}_0}{| x_0 |^2},
         \frac{x_1 \ol{x}_0}{| x_0 |^2},\ldots,
         \frac{x_n \ol{x}_0}{| x_0 |^2}\right) \\
& = f_J(x_0,x_1,\ldots,x_n)
\end{align*}
and
$$
\psi^{\ast}
d\!\left(\frac{y_{j0}+i z_{j0}}{y_{00}+i z_{00}}\right) =
d\!\left(\frac{x_j \ol{x}_0}{| x_0 |^2}\right) =
d\!\left(\frac{x_j}{x_0}\right).
$$
This shows that $\psi^{\ast}\omega'_I=\omega_I$ and we are done.
{\hfill $\blacksquare$}

{$\P_n\subseteq\P_a$:}
In this part of the proof, we will use objects over various base fields: We will
use subscripts to indicate which base field is used:
\newcommand{\trick}{\phantom{\huge A}}
\begin{center}
\begin{tabular}{|c|c|}
\hline
subscript & base field \\
\hline \hline 
$0$ & \trick \QR \trick \\
\hline 
$1$ & \trick $\ol{\Q}$ \trick \\
\hline
\R  & \trick \R \trick \\
\hline
none & \trick \C \trick \\
\hline
\end{tabular}
\end{center}
Furthermore, we fix an embedding $\ol{\Q}\subset\C$. 

Let $\int_G \omega_{\R}\in\P_n$ be a na{\"\i}ve period with
\begin{itemize}
\item
  $G\subset\R^n$ an oriented \QR-semi-algebraic set, 
  equidimensional of dimension $d$, and
\item
  $\omega_{\R}$ a rational differential $d$-form on $\R^n$ with
  coefficients in $\ol{\Q}$, which does not have poles on $G$.
\end{itemize}
The \QR-semi-algebraic set $G\subset\R^n$ is given by polynomial inequalities and
equalities. By omitting the inequalities but keeping the equalities in
the definition of $G$, we see that $G$ is supported on (the set of
\R-valued points of) a variety
$Y_{\R}\subseteq\A_{\R}^n$ of same dimension $d$. This variety $Y_{\R}$ is already
defined over \QR\
$$
Y_{\R}=Y_0 \times_{\QR} \R
$$
for a variety $Y_0\subseteq\A_{\QR}^n$ over \QR.
Similarly the boundary $\del G$ of $G$ is supported on a variety
$E_{\R}$, likewise defined over \QR
$$
E_{\R} = E_0 \times_{\QR} \R.
$$
Note that $E_0$ is a divisor on $Y_0$.
By eventually enlarging $E_0$, we
may assume w.l.o.g. that $E_0$ contains the singular locus of
$Y_0$. In order to obtain an abstract period, we need smooth varieties.
The resolution of
singularities according to Hironaka \cite{hironaka_64} provides us with a
Cartesian square
\begin{equation}
\label{cartsq}
\begin{aligned}
\widetilde{E}_0 & \subseteq \widetilde{Y}_0 \\
\downarrow & \text{\hspace{6mm}} \downarrow \pi_0 \\
E_0 & \subseteq Y_0
\end{aligned}
\end{equation}
where 
\begin{itemize}
\item
$\wt{Y}_0$ is smooth and quasi-projective,
\item
$\pi_0$ is proper, surjective and birational, and 
\item
$\widetilde{E}_0$ is a divisor with normal crossings. 
\end{itemize}
In fact, $\pi_0$ is an isomorphism away from $\wt{E}_0$
since the singular locus of $Y_0$ is contained in $E_0$
\begin{equation}
\label{p0}
\pi_{0|\wt{U}_0} : \wt{U}_0 \stackrel{\sim}{\longto} U_0
\end{equation}
with $\wt{U}_0:=\wt{Y}_0\setminus \wt{E}_0$ and 
$U_0:=Y_0\setminus E_0$.

We apply the functor {\bf an} (associated complex analytic space ---
cf. Subsection \ref{defcomplexanspace}) to
the base change to \C\ of the map $\pi_0:\wt{Y}_0\to Y_0$ 
and obtain a projection
$$
\pi_{\rm an}: \wt{Y}^{\rm an} \to Y^{\rm an}.
$$
We want to show that the ``strict transform'' of $G$
$$
\wt{G}:=\ol{\pi_{\rm an}^{-1}(G\setminus \Eh)}\subseteq\wt{Y}^{\rm an}
$$
can be triangulated. Since $\C \Ph^n$ is the projective closure of
$\C^n$, we have $\C^n \subset \C \Ph^n$ and thus get an embedding
$$
\Yh \subseteq \C^n \subset \C \Ph^n.
$$
We also choose an embedding
$$
\wt{Y}^{\rm an} \subseteq \C \Ph^m
$$
for some $m\in\N$. 
Using Lemma \ref{algsemialg}, we may consider both $\Yh$ and
$\wt{Y}^{\rm an}$ as \QR-semi-algebraic sets via some maps
\begin{align*}
& \psi: \Yh \subset \C \Ph^n \inclusion \R^N,\quad\text{and} \\
& \wt{\psi}: \wt{Y}^{\rm an} \subseteq \C \Ph^m \inclusion \R^M .
\end{align*}
In this setting, the induced
projection
$$
\pi_{\rm an} : \wt{Y}^{\rm an} \longto \Yh
$$
becomes a \QR-semi-algebraic map. The composition of $\psi$ with the
inclusion $G\subseteq\Yh$ is a \QR-semi-algebraic map; hence
$G\subset\R^N$ is \QR-semi-algebraic by Fact \ref{seidenberg_tarski}.
Since $\Eh$ is also \QR-semi-algebraic via $\psi$, we find that
$G\setminus \Eh$ is \QR-semi-algebraic. Again by Fact
\ref{seidenberg_tarski},
$\pi_{\rm an}^{-1}(G\setminus\Eh)\subset\R^M$ is 
\QR-semi-algebraic. Thus $\wt{G}\subset\R^M$, being the closure of a
\QR-semi-algebraic set, is \QR-semi-algebraic.
From Proposition \ref{semialg}, we see that $\wt{G}$ can be
triangulated
\begin{equation}
\label{unionsimpl}
\wt{G} = \union_j \triangle_j,
\end{equation}
where the $\triangle_j$ are (homeomorphic images of) $d$-dimensional
simplices. 

Our next aim is to define an algebraic differential form
$\wt{\omega}_1$ replacing $\omega_{\R}$. We first make a base change
in (\ref{cartsq}) from \QR\ to $\ol{\Q}$ and obtain
\begin{equation*}
\begin{aligned}
\widetilde{E}_1 & \subseteq \widetilde{Y}_1 \\
\downarrow & \text{\hspace{6mm}} \downarrow \pi_1 \\
E_1 & \subseteq Y_1 \, .
\end{aligned} 
\end{equation*}
The differential $d$-form $\omega_{\R}$ can be written as
\begin{equation}
\label{omegaR}
\omega_{\R}=\sum_{|J|=d} f_J(x_1,\ldots,x_n)\,
  dx_{j_1}\wedge\ldots\wedge dx_{j_d},
\end{equation}
where $x_1,\ldots,x_n$ are coordinates of $\R^n$ and
$f_J\in\ol{\Q}(x_1,\ldots,x_n)$. We
can use equation (\ref{omegaR}) to define a differential form $\omega_1$
on $\A_{\ol{\Q}}^n$
$$
\omega_{\R}=\sum_{|J|=d} f_J(x_1,\ldots,x_n)\,
dx_{j_1} \wedge\ldots\wedge dx_{j_d},
$$
where now $x_1,\ldots,x_n$ denote coordinates of $\A_{\ol{\Q}}^n$.
The pole locus of $\omega_1$ gives us a variety
$Z_1\subset\A_{\ol{\Q}}^n$.
We set
\begin{gather*}
X_1:=Y_1\setminus Z_1, \quad D_1:=E_1\setminus Z_1,\quad\text{and} \\
\wt{X}_1:=\pi_1^{-1}(X_1), \quad \wt{D}_1:=\pi_1^{-1}(D_1).
\end{gather*}
The restriction ${\omega_1}_{|X_1}$ of $\omega_1$ to $X_1$ is a
(regular) algebraic differential form on $X_1$; the pullback
$$
\wt{\omega}_1:=\pi_1^{\ast}({\omega_1}_{|X_1})
$$
is an algebraic differential form on 
$\wt{X}_1$.

We consider the complex analytic spaces $\wt{X}^{\rm an}$,
$\wt{D}^{\rm an}$, $\Zh$ associated to the base change to \C\ of $\wt{X}_1$,
$\wt{D}_1$, $Z_1$. 
Since $\omega_{\R}$ has no poles on $G$, we have 
$G \cap \Zh = \emptyset$; hence 
$\wt{G} \cap \pi_{\rm an}^{-1}(\Zh) = \emptyset$.
This shows $\wt{G}\subseteq\wt{X}=\wt{Y}\setminus\pi_{\rm an}^{-1}(\Zh)$.

Since $G$ is oriented, so is $\pi_{\rm an}^{-1}(G\setminus\Eh)$, because
$\pi_{\rm an}$ is an isomorphism away from $\Eh$. Every $d$-simplex
$\triangle_j$ in (\ref{unionsimpl}) intersects 
$\pi_{\rm an}^{-1}(G\setminus\Eh)$ 
in a dense open subset, hence inherits an orientation.
As in the proof of Proposition \ref{semialg}, we choose
orientation-preserving homeomorphisms from the standard $d$-simplex
$\Triangle_d$ to $\triangle_j$
$$
\sigma_j: \Triangle_d \longto \triangle_j.
$$
These maps sum up to a singular chain
$$
\wt{\Gamma} = \oplus_j\,\sigma_j \in \c^{\rm sing}_d(\wt{X}^{\rm an};\Q).
$$
It might happen that the boundary of the singular chain 
$\wt{\Gamma}$
is not supported on $\del \wt{G}$. Nevertheless, it will always be supported
on $\wt{D}^{\rm an}$:
The set $\pi_{\rm an}^{-1}(G\setminus\Eh)$ is oriented and therefore
the boundary components of $\del \triangle_j$ that do not belong to
$\del \wt{G}$ cancel if they have non-zero intersection with 
$\pi_{\rm an}^{-1}(G\setminus\Eh)$.
Thus $\wt{\Gamma}$ gives rise to a singular
homology class 
$$
\wt{\gamma} \in 
\h^{\rm sing}_d(\wt{X}^{\rm an},\wt{D}^{\rm an};\Q).
$$
We denote the base change to \C\ of $\omega_1$ and $\wt{\omega}_1$
by $\omega$ and $\wt{\omega}$, respectively.
Now
\begin{multline*}
\int_G \omega_{\R} =
\int_G \omega =
\int_{G\cap \Uh}\omega \\
\stackrel{(\ref{p0})}{=} \int_{\pi^{-1}(G\cap U^{\rm an})} \pi^{\ast}\omega = 
\int_{\wt{G}\cap\wt{U}^{\rm an}}\wt{\omega} \\
= \int_{\wt{G}}\wt{\omega} = 
\int_{\wt{\Gamma}}\wt{\omega} =
\int_{\wt{\gamma}}\wt{\omega} \in \P_a
\end{multline*}
is an abstract period for the quadruple
$(\wt{X}_1,\wt{D}_1,\wt{\omega}_1,\wt{\gamma})$ by Remark 
\ref{periorem1}.
\end{proof}

Now Remark \ref{periosmooth} is an easy corollary of the Period Theorem
\ref{perioperio}.

\begin{proof}[Proof of Remark \ref{periosmooth}]
Let $X_0$ be a possibly non-regular variety defined over \Q\ of
dimension $d$ and $D_0$ a
Cartier divisor on $X_0$. 
Furthermore, let 
$\omega_0 \in \Omega^d_{X_0}$ and 
$\gamma\in \h^{\rm sing}_d(X^{\rm an},D^{\rm an};\Q)$.

Exactly the same argumentation 
used for $\int_{D_I}\omega_I$ in proving $\P^{\,\prime}_p\subseteq\P_n$
shows that $\int_{\gamma} \omega \in \P_n$. 
Because of $\P_n=\P_a$, there must exists a
``smooth'' quadruple
$(\wt{X}_0,\wt{D}_0,\omega_0,\gamma)$ as described in the definition
of an abstract period (cf. Definition \ref{per2}), which gives the same period value.
\end{proof}

\begin{remark}
It might be expected that also $\P_p=\P^{\,\prime}_p$ holds. The
following ideas make this conjecture probable. Let
$\langle\gamma,\omega_0\rangle$ be a pairing period for $(X_0,D_0)$.
\begin{itemize}
\item
If $X_0$ is affine, then the complex $\wom{X_0,D_0/\Q}$ consists of quasi-coherent
sheaves by \cite[Thm. III.3.5, p. 215]{hartshorne}.
Hence we have a surjection
$$
\h\bul\Gamma(X_0;\wom{X_0,D_0/\Q}) \surjection \hDR{\Bul}(X_0,D_0/\Q),
$$
that is, only special
periods appear. 
\item
If $X_0$ is not affine, one might try to break up $\gamma$ into
subchains $\gamma_j$ with \QR-semi-algebraic representatives
$\Gamma_j$ that are already supported on $\Uh_j$ for open affine subsets
$i: U_j \inclusion X_0$. Finding a divisor $D_j$ on $U_j$ such that 
$\del \Gamma_j\subset \Dh_j$ would give us a singular chain 
$$
[i^{-1}\Gamma_j]\in\h^{\rm sing}_{\Bul}(\Uh_j,\Dh_j;\Q)
$$
and so would allow us to reduce the general case to the affine one
$$
\langle\gamma_j,\omega_0\rangle=\langle[i^{-1}\Gamma_j],i^{\ast}\omega_0\rangle.
$$
\item
If $\gamma$ and $\omega_0$ are top-dimensional of real dimension $2d$,
the procedure described above will definitely not work: A boundary
$\del \Gamma_j$ has real codimension one, but $D^{\rm an}$ has real
codimension two for $D_0$ a divisor on $X_0$. On the other hand, the
top-dimensional case is simpler, since the relative algebraic \deRham
cohomology group involved is in fact an absolute one
$$
\hDR{2d}(X_0,D_0/\Q) = \hDR{2d}(X_0/\Q),
$$
because of the long exact sequence in algebraic \deRham cohomology
and the vanishing of $\hDR{p}(D_0/\Q)$ for $p\ge {\rm
  dim}_{\R}(\Dh)=2d-2$. 

Now, generalizing a result of Huisman \cite[Thm~5.1, p.~7]{huisman} to
varieties over $\Q(i)$, we find that every algebraic \deRham
cohomology class in $\hDR{\Bul}\bigl(X_0\times_{\Q} \Q(i) / \Q(i)\bigr)$ can be
realized by a rational differential form on 
${\displaystyle\pi}_{\Q(i)/\Q}X_0$ with coefficients in \Q,
where ${\displaystyle\pi}_{\Q(i)/\Q}X_0$ is the restriction of scalars according to
Weil with respect to the field extension $\Q(i)/\Q$ of $X_0$.
\end{itemize}
\end{remark}

\medskip
The set $\P:=\P_a=\P^{\,\prime}_p=\P_n$
\notation{P}{$\P$}{equal to $\P_a$, $\P^{\,\prime}_p$ and $\P_n$}
is not only a
$\ol{\Q}$-algebra, but, at least conjecturally, comes with a triple
coproduct imposing on $\Spec\P[\frac{1}{2\pi i}]$ the structure of a
torsor, 
this notion being defined as follows.

\subsection{Torsors}

Heuristically, a torsor is a group that has forgotten its identity,
more formally:

\begin{definition}[{Torsor, \cite[p.~61]{kontsevich}}] 
A torsor
\index{Torsor} 
is a non-empty set $X$ together with
a map
$$ 
(\cdot,\cdot,\cdot) : X \times X \times X \rightarrow X 
$$
satisfying
\renewcommand{\labelenumi}{(\roman{enumi})}
\begin{enumerate}
\item $(a,a,c)=c$
\item $(a,b,b)=a$
\item $((a,b,c),d,e)=(a,(b,c,d),e)=(a,b,(c,d,e))$.
\end{enumerate}
\renewcommand{\labelenumi}{\arabic{enumi}}
\end{definition}

\begin{example}
A group $G$ becomes a torsor by setting
\begin{equation*}
\begin{aligned}
(\cdot,\cdot,\cdot) : G \times G \times G & \to G \\
(a,b,c) & \mapsto a \circ b^{-1} \circ c .
\end{aligned}
\end{equation*}

Conversely, a torsor $X$ becomes a group by choosing a distinguished
element 
$e \in X$ (`the identity') and setting
\begin{gather*}
a \circ b := (a,e,b) \\
a^{-1} := (e,a,e) .
\end{gather*}
\end{example}

\begin{example}
\label{te2}
The fibre $X$ of a principal $G$-bundle is a torsor: For any $a,b,c
\in X$ we can work out the `difference' of $a$ and $b$, which is an
element $g \in G$, and let $g$ act on $c$ in order to define
$(a,b,c)$.
\end{example}

\begin{example}
\label{te3}
Let $\cal C$ be a category and $E,F$ two isomorphic objects of $\cal
C$. The set of isomorphisms $\Iso(E,F)$ comes naturally with a right
action of the automorphism group $\Aut(E)$ and a left action of
$\Aut(F)$; both actions being simply transitive. This can be encoded
in a single map, thus turning $\Iso(E,F)$ into a torsor:
\begin{equation*}
\begin{aligned}
\Iso(E,F) \times \Iso(E,F) \times \Iso(E,F) & \rightarrow \Iso(E,F) \\
(a,b,c) & \mapsto a \circ b^{-1} \circ c.
\end{aligned}
\end{equation*} 
\end{example}

These examples are generic in the sense that we could have defined a
torsor alternatively as a principal homogenous space over a group
(Example \ref{te2}), or as a category containing just two isomorphic objects
(Example \ref{te3}).

\subsection{Fourth Definition of a Period: Effective Periods}
\label{per4}
\index{Period!effective}

We pick up our discussion of periods and show how periods could give
us a torsor.

\medskip
The usual tools for proving identities between integrals are the
{\em change of variables} and the {\em Stoke's formula}. These rules are
formalized by the following definition.

\begin{definition}[Effective Periods, {\cite[def.~20, p.~62]{kontsevich}}] 
\label{effperiodef}
The \Q -vector space $\PP_+$ 
\notation{P}{$\PP_+$}{\Q-vector space of effective periods}
of {\em effective periods} has as generators the quadruples
$(X_0,D_0,\omega_0,\gamma)$, as considered in the definition of an
abstract period (cf. Definition \ref{per2}),  modulo the following relations
\begin{itemize}
\item (linearity in both $\omega_0$ and $\gamma$)
\begin{gather*}
(X_0,D_0,q_1 \omega_{0,1}+q_2 \omega_{0,2},\gamma) \defeq
q_1 (X_0,D_0,\omega_{0,1},\gamma)+ q_2 (X_0,D_0,\omega_{0,2},\gamma) \text{ for } q_1,q_2 \in \Q,
\\
(X_0,D_0,\omega_0,q_1 \gamma_1+q_2 \gamma_2) \defeq
q_1 (X_0,D_0,\omega_0,\gamma_1)+ q_2 (X_0,D_0,\omega_0,\gamma_2) \text{ for } q_1,q_2 \in \Q.
\end{gather*}
\item (change of variables)

If $f: (X_1,D_1) \rightarrow (X_2,D_2)$ is a morphism of pairs defined
over \Q , $\gamma_1 \in \h_d^{\rm sing}(\Xh_1,\Dh_1;\Q)$ and
$\omega_2 \in \Omega^d_{X_2}$ then
$$
(X_1,D_1,f^{\ast}\omega_2,\gamma_1) \defeq (X_2,D_2,\omega_2,f_{\ast}\gamma_1).
$$
\item (Stoke's formula)

Write $D_0 = \bigcup_{i=1}^n D_i$ as union of its prime divisors,
then
\begin{equation}
\label{periostokes}
(X_0,D_0,d \omega_0, \gamma) \defeq \text{``}(D_0,\emptyset,\omega_0,\del
\gamma) \text{''} = \sum_{i=1}^n (D_i,\bigcup_{j \neq i} D_j \cap D_i,
\omega_0, \del \gamma \cap \Dh_i),
\end{equation}
where $\del : \h^{\rm sing}_d(\Xh,\Dh;\Q) \rightarrow
\h^{\rm sing}_{d-1}(\Dh;\Q)$ is the boundary operator.
\end{itemize}
\end{definition}

Here the term $\del \gamma \cap \Dh_i$ in the Stoke's formula is defined in an ad-hoc fashion:
We consider \Dh\ as a \QR-semi-algebraic set (by Lemma
\ref{algsemialg})
$$
\Dh \subset \R^N
$$
and pick a representative $\Gamma \in \del\gamma$ which is 
a rational combination $\Gamma = \sum_{k} a_{k}
\Gamma_{k}$ of \QR-semi-algebraic simplices $\Gamma_{k}$ (by
Proposition \ref{semialg}). Then the images of the $\Gamma_{k}$ are
\QR-semi-algebraic sets
$$
G_{k} := \im \Gamma_{k} \subset \R^N
$$
by Fact \ref{seidenberg_tarski}. Applying Proposition \ref{triangulation} to
the system $\{G_{k},\Dh_i\}$ yields triangulations of the
$G_{k}$'s, i.e. we can write 
$$
G_{k} = \bigcup_{l}\, \triangle_{k,l}
$$
where the
$\triangle_{k,l}$ are (homeomorphic images of) oriented
closed simplices of dimension $d-1$ whose interiors are disjoint.
They enjoy the additional property that each
$\triangle_{k,l}$ lies entirely in one of the $\Dh_i$
(because the same $\triangle_{k,l}$ triangulate also the
$\Dh_i$'s).
The induced triangulation on
$$
G_{k}^i := G_{k} \cap \Dh_i
$$
gives us a singular chain $\Gamma^i_{k}$ on $\Dh_i$ and hence a
singular homology class
$$
\gamma_{k}^i \in \h^{\rm sing}_{d-1}(\Dh_i,\union_{i \ne j} \Dh_i
\cap \Dh_j;\Q)
$$
which we call $\del \gamma \cap \Dh_i$. 

We will denote the equivalence class represented by the quadruple
$(X_0,D_0,\omega_0,\gamma)$ by $[ X_0,D_0,\omega_0,\gamma ]$.

\begin{proposition}
We have an evaluation morphism
\begin{equation*}
\begin{aligned}
{\rm ev}: \PP_+ & \rightarrow & \C \\
[ X_0,D_0,\omega_0,\gamma ] & \mapsto & \int_{\gamma} \omega .
\end{aligned}
\end{equation*}
\notation{E}{$\rm ev$}{evaluation morphism $\PP_+ \rightarrow \C$}
\end{proposition}

\begin{proof}
On page \pageref{defev1}, we have already seen that this map is
well-defined on the level of representatives
\newcommand{\ev}{\wt{\rm ev}}
$$
\wt{\rm ev}: (X_0,D_0,\omega_0,\gamma) \mapsto \int_{\gamma} \omega,
$$
so it remains to show that the map $\rm ev$ factors through the
relations imposed in Definition \ref{effperiodef}.
\begin{itemize}
\item (linearity)
\begin{multline*}
\ev(X_0,D_0,q_1 \omega_{0,1} + q_2 \omega_{0,2},\gamma) = 
\int_{\gamma} q_1 \omega_1 + q_2 \omega_2 = \\
= q_1 \int_{\gamma} \omega_1 + q_2 \int_{\gamma} \omega_2 =
q_1 \ev(X_0,D_0,\omega_{0,1},\gamma) + q_2 \ev(X_0,D_0,\omega_{0,2},\gamma),
\end{multline*}
\begin{multline*}
\ev(X_0,D_0,\omega_0,q_1 \gamma_1 + q_2 \gamma_2) = 
\int_{q_1 \gamma_1 + q_2 \gamma_2} \omega = \\
= q_1 \int_{\gamma_1} \omega + q_2 \int_{\gamma_2} \omega =
q_1 \ev(X_0,D_0,\omega_0,\gamma_1) + q_2
\ev(X_0,D_0,\omega_0,\gamma_2) ,
\end{multline*}
\item (change of variables)
$$
\ev(X_1,D_1,f^{\ast} \omega_2,\gamma_1) =
\int_{\gamma_1} f^{\ast} \omega_2 
=\int_{f_{\ast} \gamma_1} \omega_2 =
\ev(X_2,D_2,\omega_2,f_{\ast} \gamma_1) ,
$$
\item (Stoke's formula)
\begin{multline*}
\ev(X_0,D_0,d \omega_0,\gamma) = \int_{\gamma} d\omega =
\int_{\del\gamma} \omega = \\
= \sum_i \int_{\del\gamma \cap \Dh_i} \omega
= \sum_i \ev(D_i,\union_{i \ne j} D_i \cap D_j, \omega_0, \del\gamma
\cap D_i) 
\end{multline*}
\end{itemize}
\end{proof}

We state some easy properties of $\PP_+$.

\begin{proposition}[Basic properties of $\PP_+$,
  {\cite[p. 63]{kontsevich}}]
\end{proposition}
{\em
\begin{enumerate}
\item The set $\PP_+$ is a \Q -algebra with product structure as below.
\item $\QQ \subset \im({\rm ev})$.
\end{enumerate}}

\begin{proof}[Proof (cf. {\cite[p. 63]{kontsevich}})]
{\it First fact:\/} Let $[X_i,D_i,\omega_i,\gamma_i] \in \PP_+$,
$i=1,2$ and define their product to be
$$
[X_1 \times_{\Q} X_2, D_1 \times_{\Q} X_2 \cup X_1 \times_{\Q} D_2,
p_1^{\ast}\omega_1 \wedge p_2^{\ast}\omega_2, \gamma_1 \tensor{} \gamma_2] \in \PP_+,
$$
where $p_i: X_1\times_{\Q} X_2 \to X_i$, $i=1,2$, are the natural projections.
By the Fubini formula
$$
\int_{\gamma_1 \tensor{} \gamma_2} p_1^{\ast}\omega_1 \wedge p_2^{\ast}\omega_2 =
\left( \int_{\gamma_1} \omega_1 \right) \cdot \left( \int_{\gamma_2}
\omega_2 \right),
$$
i.e. the map $\rm ev$ becomes a {\Q}-algebra-morphism.

{\it Second fact:\/} Let $\alpha \in \QQ$, $p \in \Q[ t ]$ its
minimal polynomial and $X_0:=\Spec \Q [ t ]/(p)=\Spec \Q (\alpha)$.
We may regard $\alpha$ as a closed point of
$X$ via
\begin{align*}
\Spec \C & \to  X \\
\C & \leftarrow  \C[t] / (p) \\
\alpha & \mapsfrom  t .  
\end{align*}
Any closed point $x$ of $X$ gives an element of zeroth homology $\{x\} \in
\h^{\rm sing}_0(\Xh;\Q)$. Now $t \in \Omega^0_{X_0}$ and we get 
$$
{\rm ev} [ X,\emptyset,t,\{\alpha\}] = \alpha.
$$
\end{proof}

\label{evsurj}
Note that the \Q-algebra morphism ${\rm ev}: \PP_+ \to \P$ is
surjective as a consequence of the definition of abstract periods
(cf. Definition \ref{per2}).

\begin{conjecture}[{\cite[p.~62]{kontsevich}}]
\label{ev}
It is conjectured that the map $\rm ev:\PP_+\to\P$ is injective, i.e. that
all identities between integrals can be proved using standard rules
only.
\end{conjecture} 
This would give us a \Q-algebra-isomorphism $\PP_+ \iso \P$.

\begin{center}
\em
For the remainder of this section, we assume Conjecture \ref{ev}.
\end{center}

We adjoin formally $(2\pi i)^{-1}$ to $\PP_+$
$$
\PP:=\PP_+[{\textstyle\frac{1}{(2\pi i)}}]\stackrel{?}{=}
\P[{\textstyle\frac{1}{(2\pi i)}}]
$$
  \notation{P}{$\PP$}{$\PP_+$ adjoint $\frac{1}{2\pi i}$}
and define a coproduct on $\PP$.

\begin{definition}[The triple coproduct $\Delta$ on $\PP$,
    {\cite[p.~63]{kontsevich}}]
\label{coprod}
\index{Triple coproduct|textbf}
\end{definition}
\vspace{-2.5mm}
{\em
Let $\alpha:=[X_0,D_0,\omega_0,\gamma]$ be an effective period
represented by an quadruple
$(X_0, D_0, \omega_0, \gamma)$.
We write $P$
for the period matrix of $(X_0,D_0)$ 
as defined in Definition \ref{per1}
and assume 
${\rm ev}(\alpha)=P_{i,j}$, the entry in $P$ at position $(i,j)$.
We define
$$
\Delta(P_{i,j}):=
\sum_{k,l} P_{i,k} \tensor{} (P^{-1})_{k,l} \tensor{} P_{l,j}
\in\PP\tensor{\Q}\PP\tensor{\Q}\PP.
$$
Here we used Remark \ref{det} about the determinant of $P$: 
${\rm det}\,P$  being a
period in $\P$ times a power of $(2\pi i)$.
Thus $P^{-1}$ has coefficients in 
$\PP=\PP_+[\frac{1}{2\pi i}] \stackrel{?}{=} \P[\frac{1}{2\pi i}]$.

We adopt the convention
$$
\Delta\left(\frac{P_{i,j}}{(2\pi i)^n}\right):=
\sum_{k,l} \frac{P_{i,k}}{(2\pi i)^n} \tensor{} (2\pi i)^n
(P^{-1})_{k,l} \tensor{} \frac{P_{l,j}}{(2\pi i)^n}
\quad\text{for}\quad n\in\N_0.
$$

Extending the map $\Delta$ in a $\Q$-linear way to all of $\PP$
yields the {\em triple coproduct}
$$
\Delta: \PP \to \PP\tensor{\Q}\PP\tensor{\Q}\PP.
$$
\notation{D}{$\Delta$}{triple coproduct}
}

\begin{remark}
It is far from being obvious that $\Delta$ is well-defined. In his
paper \cite{kontsevich}, Kontsevich states that the proof of the
correctness of Definition \ref{coprod} relies on an
unpublished result by M.~Nori \cite[Thm. 6, p.~63]{kontsevich}.
\end{remark}

According to \cite{kontsevich}, 
also the following proposition would follow from Nori's theorem.
\begin{proposition}
Modulo Conjecture \ref{ev} and Nori's theorem, the triple coproduct on $\PP$
induces the structure of a torsor on $\Spec\PP$.
\end{proposition}

We calculate a triple coproduct in Example \ref{exp1}.

\section{Examples}

\subsection{First Example: $\A^1_{\C}\setminus\{0\}$}
\label{exp1}

{\em First part (cf. {\cite[p.~63]{kontsevich}}):}
Let
$$
\Xo := \A^1_{\Q} \setminus \{ 0 \} = \Spec \Q[t,t^{-1}]
$$
be the affine line with the point $0$ deleted and 
$$
D_0 := \{ 1, \alpha \} \quad \text{with} \quad \alpha \ne 0,1 
$$
a divisor on \Xo.
\begin{figure}[h]
\begin{center}
\epsfig{file=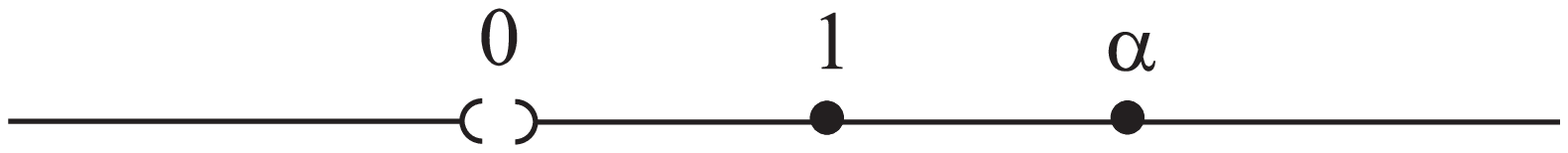,width=9.5cm}
\caption{The algebraic pair $(X_0,D_0)$}
\label{exp1a}
\end{center}
\vspace{-5mm}
\end{figure}

Then the {\em singular homology} of the pair
$(\Xh,\Dh)=(\C^{\times},\{1,\alpha\})$ 
(cf. Definition \ref{relsingcoho})
is generated by a
small loop $\sigma$ turning counter-clockwise around $0$ once and the
interval $[1,\alpha]$. 

In order to compute the 
{\em algebraic \deRham cohomology} of $(\Xo,\Do)$ (cf. Definition
\ref{algdr3}), 
we first note that 
$$
\hDR{\Bul}(\Xo,\Do/\Q) = \hh\bul\Gamma(\Xo;\wom{\Xo,\Do/\Q}),
$$
since $\Xo$ is affine and 
the sheaves $\wt{\Omega}^p_{\Xo,\Do/\Q}$ are quasi-coherent, 
hence acyclic for the global section functor 
by \cite[Thm. III.3.5, p.~215]{hartshorne}. 

We spell out the complex $\Gamma(X_0;\wom{\Xo,\Do/\Q})$ in detail (cf. page \pageref{defwomrel})
\begin{align*}
0 \hspace{1cm} & \\
\up[\phantom{d}] \hspace{1cm} & \\
\Gamma(\Xo;\wt{\Omega}^1_{\Xo,\Do/\Q}) & = 
\Gamma\bigl(\Xo;\Omega^1_{\Xo/\Q}
\oplus \bigoplus_j\, \iast \OO_{D_j}\bigr) = 
\Q[t,t^{-1}]dt \oplus \underset{1}{\Q} \oplus \underset{\alpha}{\Q} \\[-2mm]
\up[d] \hspace{1cm} & \\
\Gamma(\Xo;\OO_{\Xo}) & = \Q[t,t^{-1}] 
\end{align*}
and observe that the  map ``evaluation-at-a-point''
\begin{align*}
\Q[t,t^{-1}] & \surjection \underset{1}{\Q} \oplus \underset{\alpha}{\Q} \\
f(t) & \mapsto \bigl(f(1),f(\alpha)\bigr)
\end{align*}
is surjective with kernel
$$
(t-1)(t-\alpha)\Q[t,t^{-1}] =
\linspan_{\Q}\{t^{n+2}-(\alpha+1)t^{n+1}+\alpha t^n \st n\in\Z\}.
$$
Differentiation maps this kernel to 
$$
\linspan_{\Q} \{ (n+2)t^{n+1} - (n+1)(\alpha+1)t^n - 
n\alpha t^{n-1} \st n\in \Z\}dt.
$$
Therefore we get
\begin{align*}
\hDR{1}(\Xo,\Do;\Q) 
& = \Gamma(X_0;\wt{\Omega}_{X_0,D_0/\Q})\,/\,\Gamma(X_0;\OO_{X_0}) \\
& = \Q[t,t^{-1}]dt\oplus
  \underset{1}{\Q}\oplus\underset{\alpha}{\Q}\,/\,
  d(\Q[t,t^{-1}]) \\[-1mm]
& = \Q[t,t^{-1}]dt/
\linspan_{\Q} \{ (n+2)t^{n+1} - (n+1)(\alpha+1)t^n - 
n\alpha t^{n-1}\}dt.
\end{align*}
By the last line, we see that the class of $t^ndt$ in
$\hDR{1}(X_0,D_0;\Q)$ for $n\ne-1$ is linearly dependent of 
\begin{itemize}
\item
$t^{n-1}dt$ and $t^{n-2}dt$, and
\item
$t^{n+1}dt$ and $t^{n+2}dt$,
\end{itemize}
hence linearly dependent of $\frac{dt}{t}$ and $dt$ by an induction
argument.
Hence $\hDR{1}(X_0,D_0;\Q)$ is spanned by 
$$
\frac{dt}{t} \quad \text{and} \quad \frac{1}{\alpha-1}dt .
$$
We obtain the following {\em period matrix} $P$ for $(\Xo,\Do)$
(cf. Definition \ref{per1})

\begin{equation}
\label{permat}
\begin{tabular}{c|cc}
& $\frac{1}{\alpha-1}dt$ & $\frac{dt}{t}$ \\
\hline 
$[1,\alpha]$ & $1$ & $\ln\alpha$ \\[1mm]
$\sigma$ & $0$ & $2\pi i$ \, .
\end{tabular}
\end{equation}
Let us the compute the {\em triple coproduct} of $\ln\alpha$
\index{Triple coproduct} (cf. Definition \ref{coprod})
from this example. We have
$$
P^{-1}=
\left(
  \begin{matrix}
1 & \frac{-\ln\alpha}{2\pi i} \\[1mm]
0 & \frac{1}{2\pi i}
  \end{matrix}
\right)
$$
and thus get for the triple coproduct (cf. \cite[p.~63]{kontsevich})
\begin{equation}
\label{coprodln}
\textstyle
\Delta(\ln\alpha)=
\ln\alpha\tensor{}\frac{1}{2\pi i}\tensor{}2\pi i
-1\tensor{}\frac{\ln\alpha}{2\pi i}\tensor{}2\pi i
+1\tensor{}1\tensor{}\ln\alpha.
\end{equation}

\bigskip
{\em Second part:} Now consider the degenerate configuration $\alpha=1$,
i.e. $D=2\cdot\{1\}$ with $D_1=D_2=\{1\}$. Although $D_0$ is not a
divisor with normal crossings anymore, the machinery developed in
Section \ref{alg} still works.
\begin{figure}[h]
\begin{center}
\epsfig{file=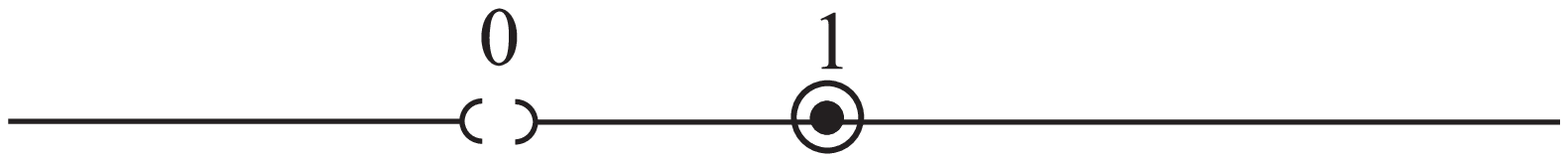,width=9.5cm}
\caption{The algebraic pair $(X_0,D_0)$}
\label{exp1b}
\end{center}
\vspace{-5mm}
\end{figure}

Repeating the calculations
done for the general case gives
$$
\hDR{1}(X_0,D_0/\Q) = \Q\, \frac{dt}{t} \oplus
\Q\,(1_{D_1} 
+ 0_{D_2}),
$$
where $c_{D_i}$ is the constant function equal to $c$ on the
irreducible component $D_i$ of $D$,
$i=1,2$. Computing the smooth singular cohomology group 
$\h^1_{\infty}(\Xh,\Dh;\Q)$ 
(cf. Definition \ref{defsmooth2})
in a similar fashion gives
$$
\h^1_{\infty}(\Xh,\Dh;\Q) = \Q\,\sigma \oplus \Q\,(D_1 - D_2),
$$
and we find that the period matrix of $(X_0,D_0)$ is precisely the
limit of (\ref{permat}) for $\alpha \to 1$
\begin{center}
\begin{tabular}{c|cc}
& $1_{D_1}+0_{D_2}$ & $\frac{dt}{t}$ \\
\hline
$D_1-D_2$ & $1$ & $0$ \\
$\sigma$ & $0$ & $2\pi i$
\end{tabular}
\end{center}

\bigskip
{\em Third part:}
Finally, we let 
$$
D_0:=\{1,\alpha,\beta\}
\quad\text{with}\quad
\alpha\ne 0,1
\quad\text{and}\quad
\beta\ne 0,1,\alpha,
$$
but keep $X_0:=\A^1_{\Q}\setminus\{0\}=\Spec \Q[t,t^{-1}]$. 
\begin{figure}[h]
\begin{center}
\epsfig{file=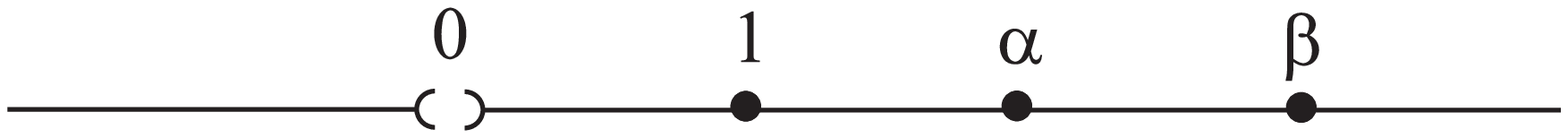,width=9.5cm}
\caption{The algebraic pair $(X_0,D_0)$}
\label{exp1c}
\end{center}
\vspace{-5mm}
\end{figure}

Then
$\h^{\rm sing}_1(\Xh,\Dh;\Q)$
is generated by 
the loop $\sigma$ from the first part
and the intervals $[1,\alpha]$ and $[\alpha,\beta]$. 
Hence the differential forms $\frac{dt}{t}$, $dt$ and
$2t\,dt$ give a basis of 
$\hDR{1}(X_0,D_0;\Q)$: 
If they were linearly dependent, the 
period matrix $P$ would not be of full rank

\begin{center}
\begin{tabular}{c|ccc}
& $\frac{dt}{t}$ & $dt$ & $2t\,dt$ \\
\hline 
$\sigma$ & $2\pi i$ & $0$ & $0$ \\[1mm]
$[1,\alpha]$ & $\ln\alpha$ & $\alpha-1$ & $\alpha^2-1$ \\[1mm]
$[\alpha,\beta]$ & $\ln\left(\frac{\beta}{\alpha}\right)$ &
$\beta-\alpha$ & \,\,$\beta^2-\alpha^2$\,.
\end{tabular}
\end{center}

Observe that 
$\det P=2\pi i (\alpha-1)(\beta-\alpha)(\beta-1)\ne 0$.

We can compute the triple coproduct of $\ln\alpha$ again.
We have
$$
P^{-1}=\left(\begin{matrix}
\frac{1}{2\pi i} & 0 & 0 \\[1mm]
\frac{\ln\beta(\alpha^2-1)-\ln\alpha(\beta^2-1)}{2\pi
  i(\beta-\alpha)(\alpha-1)(\beta-1)} &
\frac{\alpha+\beta}{(\alpha-1)(\beta-1)} &
\frac{\alpha+1}{(\alpha-\beta)(\beta-1)} \\[1mm]
\frac{-\ln\beta(\alpha-1)+\ln\alpha(\beta-1)}{2\pi
  i(\beta-\alpha)(\alpha-1)(\beta-1)} &
\frac{-1}{(\alpha-1)(\beta-1)} & 
\frac{-1}{(\alpha-\beta)(\beta-1)} 
\end{matrix}\right)\, ,
$$
and therefore get for the triple coproduct
\begin{align*}
\Delta(\ln\alpha)=&
\ln\alpha\tensor{}\frac{1}{2\pi i}\tensor{}2\pi i \displaybreak[0] \\
&\hspace{0.5cm}+ (\alpha-1)\tensor{}
  \frac{-\ln\beta(\alpha^2-1)+\ln\alpha(\beta^2-1)}
    {2\pi i(\beta-\alpha)(\alpha-1)(\beta-1)}\tensor{}
  2\pi i \displaybreak[0] \\
&\hspace{1cm}+ (\alpha-1)\tensor{}
  \frac{\alpha+\beta}{(\alpha-1)(\beta-1)}\tensor{}
  \ln\alpha \displaybreak[0] \\
&\hspace{1.5cm}+ (\alpha-1)\tensor{}\frac{\alpha+1}{(\alpha-\beta)(\beta-1)}\tensor{}
  \ln\left(\frac{\beta}{\alpha}\right) \displaybreak[0] \\
&\hspace{2cm}+ (\alpha^2-1)\tensor{}
  \frac{\ln\beta(\alpha-1)-\ln\alpha(\beta-1)}
    {2\pi i(\beta-\alpha)(\alpha-1)(\beta-1)}\tensor{}
  2\pi i \displaybreak[0] \\
&\hspace{2.5cm}+ (\alpha^2-1)\tensor{}
  \frac{-1}{(\alpha-1)(\beta-1)}\tensor{}
  \ln\alpha \displaybreak[0] \\
&\hspace{3cm}+ (\alpha^2-1)\tensor{}
  \frac{-1}{(\alpha-\beta)(\beta-1)}\tensor{}
  \ln\left(\frac{\beta}{\alpha}\right) \displaybreak[0] \\
=&\ln\alpha\tensor{}\frac{1}{2\pi i}\tensor{}2\pi i
-1\tensor{}\frac{\ln\alpha}{2\pi i}\tensor{}2\pi i
+1\tensor{}1\tensor{}\ln\alpha,
\end{align*}
i.e. we obtain the same result as in (\ref{coprodln}).

\subsection{Second Example: Quadratic Forms}
\index{Quadratic form}

Let
$$
\begin{matrix}
Q(\ul{x}) : & \Q^3 & \longrightarrow & \Q\\
& \ul{x}=(x_0,x_1,x_2) & \mapsto & \ul{x}\, A\, \ul{x}^T
\end{matrix}
$$
be a quadratic form with $A\in \Q^{3\times 3}$ being a regular,
symmetric matrix.

The zero-locus of $Q(\ul{x})$
$$
\ol{X}_0 := \{ \ul{x} \in \Q P^2 \st Q(\ul{x}) = 0 \}
$$ 
is a {\em quadric} or non-degenerate {\em conic}. We are interested
in its affine piece
$$
X_0:= \ol{X}_0 \cap \{x_0 \ne 0\} \subset \Q^2 \subset \Q P^2.
$$

We show that we can assume $Q(\ul{x})$ to be of a particular nice form.
A vector $v\in \Q^3$ is called {\em $Q$-anisotropic}, if
$Q(v)\ne 0$. Since ${\rm char}\,\Q \ne 2$, there exist such vectors;
just suppose the contrary:
$$
\begin{aligned}
Q(1,0,0) = 0 & \quad\text{gives} & A_{11}=0 , \\
Q(0,1,0) = 0 & \quad\text{gives} & A_{22}=0, \\
Q(1,1,0) = 0 & \quad\text{gives} & 2\cdot A_{12}=0
\end{aligned}
$$
and $A$ would be degenerate. In particular
$$
Q(1,\lambda,0) = Q(1,0,0) + 2 \lambda Q(1,1,0) + 
\lambda^2 Q(0,1,0)
$$
will be different form zero for almost all $\lambda \in 
\Q$. Hence, we can assume that $(1,0,0)$ is anisotropic after applying
a coordinate transformation of the form 
$$
x'_0:= x_0, \quad x'_1:= -\lambda x_0 + x_1, \quad x'_2:= x_2 .
$$ 
After another affine change of coordinates, we can also assume that
$A$ is a diagonal matrix by \cite[Ch. IV, \S\ 1, 2, Thm. 1, p. 30]{serre}.
An inspection of the proof of this result reveals that we can choose
this coordinate transformation such that the $x_0$-coordinate is
left unaltered.
(Just take for $e_1$ the anisotropic vector $(1,0,0)$ in the proof.)
Such a transformation does not change the isomorphy-type of $X_0$, and
we can take $X_0$ to be cut out by an equation of the form 
$$
a x^2 + b y^2 = 1 \quad \text{for} \quad a,b\in\Q^{\times}
$$
with affine coordinates $x:=\frac{x_1}{x_0}$ and $y:=\frac{x_2}{x_0}$.

Since $X_0$ is affine, and hence the sheaves $\Omega^p_{X_0/\Q}$ are 
$\Gamma(X_0;?)$-acyclic by \cite[Thm. III.3.5, p.~215]{hartshorne}, 
we can compute its {\em algebraic 
\deRham cohomology} (cf. Definition \ref{algdr1}) by 
$$
\hDR{\Bul}(\Xo/\Q) = \hh\bul\Gamma(\Xo;\om{\Xo/\Q});
$$
so we write down the complex $\Gamma(X_0;\om{\Xo/\Q})$ in detail
\begin{align*}
0 \hspace{1cm} & \\
\uparrow \hspace{1cm} & \\
\Gamma(\Xo;\Omega^1_{\Xo/\Q}) & = 
\Q[x,y]/(ax^2+by^2-1) \{dx,dy\}
\, / \, (axdx+bydy) \\
{\scriptstyle d}\uparrow \hspace{1cm} & \\
\Gamma(\Xo;\OO_{\Xo}) & = \Q[x,y]/(ax^2+by^2-1) .
\end{align*}
Obviously, $\hDR{1}(X_0/\Q)$ can be presented with generators $x^n
y^m dx$ and $x^n y^m dy$ for $m,n\in \N_0$ modulo numerous relations.
We easily get
\begin{itemize}
\item
$y^m \,dy = d\, \frac{y^{m+1}}{m+1} \sim 0$
\item
$x^n \,dx = d\, \frac{x^{n+1}}{n+1} \sim 0$
\item[{\footnotesize $n\ge 1$}$\quad\bullet$]
$
x^n y^m \,dy = \frac{-n}{m+1} x^{n-1}y^{m+1} \,dx + d\, \frac{x^n
  y^{m+1}}{m+1} \\[1mm]
\phantom{x^n y^m \,dy} \sim \frac{-n}{m+1} x^{n-1}y^{m+1} \,dx
\quad\text{for}\quad n\ge 1
$
\item
$x^n y^{2m} \,dx = x^n \big( \frac{1-ax^2}{b} \big)^m \,dx$
\item
$x^n y^{2m+1} \,dx = x^n \big( \frac{1-ax^2}{b} \big)^m y \,dx$
\item
$ 
xy \,dx = \frac{-x^2}{2} \,dy + d\, \frac{x^2 y}{2} \\[1mm]
\phantom{xy \,dx} \sim \frac{by^2-1}{2a} \,dy \\[1mm]
\phantom{xy \,dx} = \frac{b}{2a}y^2\,dy-\frac{1}{2a}\,dy \sim 0 
$
\item[{\footnotesize$n\ge 2$}\quad$\bullet$]
$
x^ny\,dx = 
\frac{-b}{a}x^{n-1} y^2\,dy + x^n y\,dx + \frac{b}{a}x^{n-1} y^2\, dy \\[1mm]
\phantom{x^ny\,dx} =\frac{-b}{a}x^{n-1}y^2\,dy +\frac{x^{n-1}y}{2a}\,
d(a x^2 + b y^2 -1 ) \\[1mm]
\phantom{x^ny\,dx} = \frac{-b}{a} x^{n-1} y^2 \,dy + d\, \big( \frac{ (x^{n-1} y) ( ax^2 +
  b y^2 -1)}{2a} \big) \\[1mm]
\phantom{x^ny\,dx} \sim \frac{-b}{a} x^{n-1}y^2 \,dy \\[1mm]
\phantom{x^ny\,dx} = \big( x^{n+1} - \frac{x^{n-1}}{a} \big) \,dy \\[1mm]
\phantom{x^ny\,dx} = \big(- (n+1) x^ny + \frac{n-1}{a} x^{n-2} y \big) \,dx 
+ d\, \big(x^{n+1}y - \frac{x^{n-1}}{a}y \big)
$
\item[$\Rightarrow$]
$x^n y \,dx \sim \frac{n-1}{(n+2)a} x^{n-2}y \,dx$\quad for\quad $n\ge 2$.
\end{itemize}
Thus we see that all generators are linearly dependent of $y\, dx$
\begin{align*}
\hDR{1}(X_0/\Q) & = \hh^1\Gamma(X_0;\om{X_0/\Q}) \\
& = \Q \,\,y\,dx.
\end{align*}
What about $X$, the base change to $\C$ of $X_0$? 
We use the symbol $\sqrt{\phantom{a}}$ for the principal branch of the
square root. 
Over \C, the change of coordinates
$$
u:=\sqrt{a} x - i\sqrt{b}y, \quad v:=\sqrt{a} x + i \sqrt{b} y
$$
gives
\begin{align*}
X =\, & \Spec \C [x,y] / (ax^2 + by^2 -1) \\
=\, & \Spec \C [u,v] / ( uv-1) \\
=\, & \Spec \C [u,u^{-1}] \\
=\, & \A^1_{\C } \setminus \{0\}.
\end{align*}
Hence the {\em first singular homology group} $\h_{\Bul}^{\rm sing}(\Xh;\Q)$
of $\Xh$ is generated by
$$
\sigma : [0,1] \to \Xh, s \mapsto u=e^{2\pi i s},
$$ 
i.e. a circle with radius $1$ turning counter-clockwise around 
$u=0$ once.

The {\em period matrix} (cf. Definition \ref{per1}) consists of a single entry
\begin{align*}
\int_{\sigma} y\, dx 
= & \,
\int_{\sigma} \frac{v-u}{2i\sqrt{b}}\,\, d\,
\frac{u+v}{2\sqrt{a}} \\
\stackrel{\makebox[0mm]{\scriptsize Stokes}}{=} & \,
\int_{\sigma}  \frac{v\, du - u\,
  dv}{4i\sqrt{ab}} \\
= & \, \frac{1}{2i\sqrt{ab}} \int_{\sigma} \frac{du}{u}\\
= & \,\frac{\pi}{\sqrt{ab}} \\
= & \,\frac{\pi}{\sqrt{\text{\scriptsize discriminant}}}.
\end{align*}
The denominator squared is nothing but the {\em discriminant}  of the
quadratic form $Q$ (cf. \cite[Ch. IV, \S 1, 1, p. 27]{serre})
\index{Quadratic form!discriminant of}
$$
{\rm disc}\, Q := {\rm det}\, A \in \Q^{\times}/{\Q^{\times}}^2.
$$
This is an important invariant, that distinguishes some, but not all
isomorphy classes of quadratic forms.
Since ${\rm disc}\,Q$ is well-defined modulo
$({\Q^{\times}})^2$, it makes sense to write
$$
\hDR{1}(X_0/\Q) = \Q\, \frac{\pi}{\sqrt{{\rm disc}\,Q}} 
\subset \h^1\sing(\Xh;\Q)\tensor{\Q}\C.
$$

\subsection{Third Example: Elliptic Curves}
\index{Elliptic curves}

In this example, we give a brief summary of Chapter VI in Silverman's
classical textbook \cite{silverman} on elliptic curves. His Chapter VI
deals with elliptic curves over \C, but some results generalize easily
to the case of elliptic curves over \Q.

An {\em elliptic curve} $E_0$ is an one-dimensional, nonsingular, 
complete scheme defined over a field $k$, which has genus $1$, where the
{\em geometric genus} $p_g$ is defined as
$$
p_g(E_0) := {\rm dim}_k\, \Gamma(E_0;\Omega^1_{E_0/k}) .
$$
For simplicity, we assume $k=\Q$. It can be shown, using the
Riemann-Roch~theorem, that such an elliptic curve $E_0$ can be given
as the zero locus in $\Q P^2$ of a {\em Weierstra\ss-equation}
(cf. \cite[Ch. III, \S 3, p. 63]{silverman})
\begin{equation}
\label{weierstrass}
Y^2Z=4X^3-60G_4XZ^2-140G_6Z^3
\end{equation}
with coefficients $G_4, G_6\in\Q$ and projective coordinates $X$, $Y$
and $Z$.

\begin{figure}[h]
\begin{center}
\epsfig{file=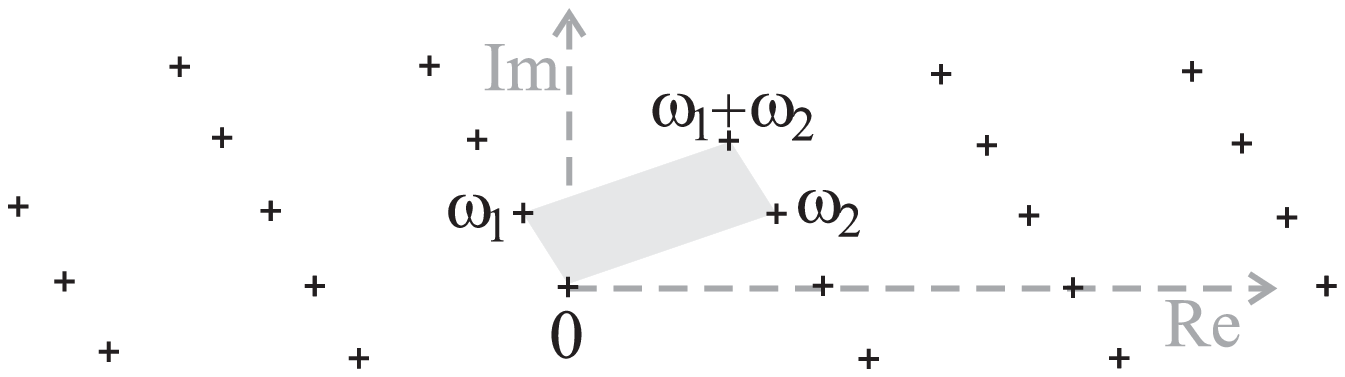,width=7.5cm}
\caption{The lattice $\Lambda_{\omega_1,\,\omega_2}\subset\C$}
\label{lattice}
\end{center}
\end{figure}

The base change $E:=E_0 \times_{\Q} \C$ of $E_0$ gives us a {\rm
complex torus} $\Eh$ (see \cite[Ch. VI, \S 5, Thm. 5.1,
p. 161]{silverman}), i.e. an isomorphism
\begin{equation}
\label{eanclt}
\Eh \iso \C/\Lambda_{\omega_1,\,\omega_2}
\end{equation}
in the complex analytic category 
(cf. Subsection \ref{defcomplexanspace}), 
with 
\begin{gather*}
\Lambda_{\omega_1,\,\omega_2} 
:= \omega_1 \Z \oplus \omega_2 \Z \\
\text{for } \omega_1,\omega_2 \in \C \text{ linearly
  independent over \R},
\end{gather*}
being a 
lattice of full rank (cf. Figure \ref{lattice}).

\begin{figure}[b]
\begin{center}
\epsfig{file=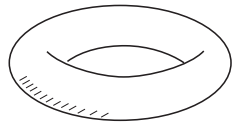,width=3cm}
\caption{The standard torus $S^1\times S^1$}
\label{torus}
\end{center}
\end{figure}

Thus
all elliptic curves over \C\ are diffeomorphic to the standard torus $S^1\times
S^1$ (cf. Figure \ref{torus}), but carry different complex structures as the parameter
$\tau:=\omega_2/\omega_1$ varies.
We can describe the isomorphism (\ref{eanclt}) quite explicitly using
periods.
 Let $\alpha$ and $\beta$ be a basis of 
$$
\h_1^{\rm sing}(\Eh;\Z) = 
\h_1^{\rm sing}(S^1\times{}S^1;\Z) =
\Z\,\alpha \, \oplus \, \Z\,\beta.
$$ 
The \Q-vector space $\Gamma(E_0;\Omega^1_{E_0/\Q})$ is
spanned by the so-called {\em invariant differential}
\cite[p.~48]{silverman}
$$
\omega=\frac{d(X/Z)}{Y/Z}.
$$
One
can show (cf. \cite[Ch. VI, \S 5, Prop. 5.2,
  p. 161]{silverman}), that the (abstract) periods (cf. Definition \ref{per2})
$$
\omega_1:=\int_{\alpha} \omega\quad\text{and}\quad
\omega_2:=\int_{\beta}\omega
$$
are \R-linearly independent and that the
map
\begin{equation}
\label{elliso}
\begin{aligned}
\Eh \to & \quad \C/\Lambda_{\omega_1,\,\omega_2} \\
P \mapsto & \int_O^P \omega \,\text{ modulo }
\Lambda_{\omega_1,\,\omega_2}
\end{aligned}
\end{equation}
is an isomorphism (here $O=[0:1:0]$ denotes a base-point on $E$). 

Note that under the natural projection $\pi: \C \to
\C/\Gamma_{\omega_1,\,\omega_2}$ any meromorphic function $f$ on the
torus 
$\C/\Gamma_{\omega_1,\,\omega_2}$
lifts to a doubly-periodic function $\pi^{\ast}f$ on the complex plane
$\C$ with periods $\omega_1$ and $\omega_2$
$$
f(x+n\omega_1+m\omega_2) = f(x)
\quad\text{for all}\quad n,m\in\Z\quad\text{and}\quad x\in\C.
$$
This example is possibly the origin of the term ``period'' for the
elements of the sets $\P_p$, $\P_a$, and $\P_n$.

The inverse map $\C/\Lambda_{\omega_1,\,\omega_2} 
\to \Eh$ 
for the isomorphism (\ref{elliso}) can be described in terms of the {\em
Weierstra\ss-{}$\wp$-function} of the lattice
$\Lambda:=\Lambda_{\omega_1,\,\omega_2}$ 
$$
\wp(z)=\wp(z,\Lambda):=\frac{1}{z^2} + 
\sum_{\substack{
\omega\in\Lambda \\
\omega\ne 0}}
\frac{1}{(z-\omega)^2} - \frac{1}{\omega^2}
$$
and takes the form (cf. \cite[Prop.~3.6, p.~158]{silverman})
\begin{align*}
\C/\Lambda_{\omega_1,\,\omega_2} & \to \Eh \subset \C \Ph^2 \\
z & \mapsto [\wp(z):\wp'(z):1] .\
\end{align*}
The defining coefficients $G_4, G_6$ of $E_0$ can be recovered from
$\Lambda_{\omega_1,\,\omega_2}$ by 
computing {\em Eisenstein series} (cf. \cite[Ch. VI, \S 3, Thm. 3.5,
  p. 157]{silverman})
$$
G_{2k} := 
\sum_{\substack{
\omega\in\Lambda \\
\omega\ne 0}}
\omega^{-2k}\quad\text{for}\quad k=2,3.
$$

Thus the periods $\omega_1$ and $\omega_2$ determine the elliptic
curve $E_0$ uniquely. However, they are not invariants of $E_0$, since
they depend on the chosen Weierstra\ss-equation (\ref{weierstrass})
of $E_0$. It is proved in \cite[p.~50]{silverman}, that a change of
coordinates which preserves the shape of (\ref{weierstrass}), must
be of the form
$$
X'=u^2 X, \quad Y'=u^3 Y, \quad Z'=Z\quad \text{for}\quad u\in\Q^{\times}.
$$
In the new parameterization $X', Y', Z'$, we have
\begin{align*}
& G'_4=u^4 G_4 ,\quad G'_6=u^6 G_6, \\
& \omega' = u^{-1} \omega \\
& \omega'_1 = u^{-1} \omega_1 \quad\text{and}\quad \omega'_2=u^{-1} \omega_2.
\end{align*}
Hence $\tau=\omega_2/\omega_1$ is a true invariant of $E_0$, which
allows us to distinguish non-isomorphic elliptic curves $E$ over
\C. If we consider both $\tau$ and the image of $\omega_1$ in
$\C^{\times}/\Q^{\times}$, we can even distinguish elliptic curves
over \Q.

We may ask whether the isomorphy-type of $E_0$ is already encoded in its
cohomology. The answer is yes if we consider also the
pure Hodge structure living on cohomology.

By the Hodge theorem \cite[p.~116]{griffiths_harris}, the (hypercohomology) spectral sequence
$$
\h^p(\Eh;\Omega^q_{\Eh}) \Rightarrow
\hDR{n}(\Eh;\C)
$$
splits, yielding a decomposition of 
the first classical \deRham cohomology group 
(cf. Subsection \ref{drclassic})
$$
\hDR{1}(\Eh;\C)=
\Gamma(\Eh;\Omega^1_{\Eh})\oplus
\h^1(\Eh;\OO_{\Eh}).
$$
This Hodge decomposition already lives on the algebraic \deRham
cohomology group $\hDR{1}(E_0/\Q)$ (cf. Definition \ref{algdr1})
\begin{equation}
\label{hodgedecomp}
\hDR{1}(E_0/\Q)=
\Gamma(E_0;\Omega^1_{E_0/\Q})\oplus
\h^1(E_0;\OO_{E_0}).
\end{equation}
In order to see this, observe that the morphism of spectral sequences
\begin{gather*}
\h^p(E_0;\Omega^q_{E_0/\Q})\tensor{\Q}\C \Rightarrow 
\hDR{n}(E_0;\Q)\tensor{\Q}\C \\
\downarrow \hspace{3cm} \downarrow \\
\h^p(\Eh;\Omega^q_{\Eh}) \Rightarrow 
\hDR{n}(\Eh;\C)
\end{gather*}
is indeed an isomorphism by the Base Change Proposition \ref{basechange},
the {\small GAGA}-theorem \ref{gagathm} and the Comparison Theorem \ref{comp}.
Hence the spectral sequence 
$$
\h^p(E_0;\Omega^q_{E_0/\Q}) \Rightarrow \hDR{n}(E_0;\Q) 
$$
splits as well  and we obtain the decomposition (\ref{hodgedecomp}).

Pairing $\h_1^{\rm sing}(\Eh;\Z)=\Z^2$ with some 
$\omega'\in\Gamma(\Eh;\Omega^1_{\Eh})=\Q\,\omega$ gives us
$\Lambda_{u\cdot\omega_1,u\cdot\omega_2}\subset\C$, where we write
$\omega'=u\cdot\omega$ with $u\in\Q^{\times}$.
That is, we can recover $\Lambda_{\omega_1,\,\omega_2}$ up to a rational
multiple 
from the integral singular homology $\h_1^{\rm sing}(\Eh;\Z)$ of $\Eh$
plus
the algebraic \deRham cohomology of $E_0$ equipped with its 
Hodge decomposition (\ref{hodgedecomp}),
i.e. we can recover the isomorphy-type of $E_0$ from this data.

\subsection{Fourth Example: A $\zeta$-value}
\label{exp4}
\index{zeta-value@$\zeta$-value}

Earlier we raised the question of how to write 
$\zeta(2)$ as a period and found the identity (cf. Proposition \ref{zetaprop})
$$
\zeta(2)=
\int_{0\,\le\,x\,\le\,y\,\le\,1} \frac{dx \wedge dy}{(1-x)\,y} .
$$
The problem was that this identity did not give us a valid
representation of $\zeta(2)$ as a na{\"\i}ve period (cf. Definition \ref{per3}),
since the pole locus of the integrand and the domain of integration
are not disjoint. We show how to circumvent this difficulty.

\begin{figure}[h]
\begin{center}
\epsfig{file=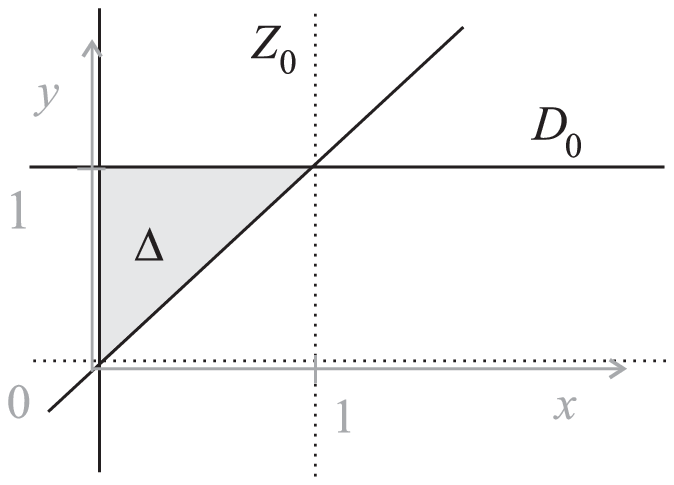,width=6cm}
\caption{The configuration $Z_0,D_0,\triangle$}
\label{exp4a}
\end{center}
\vspace{-5mm}
\end{figure}

First we define
\begin{align*}
&Y_0:=\A^2_{\Q}\quad\text{with coordinates $x$ and $y$,}\\
&Z_0:=\{x=1\}\cup\{y=0\},\\
&X_0:=Y_0\setminus Z_0,\\
&D_0:=(\{x=0\}\cup\{y=1\}\cup\{x=y\})\setminus Z_0,\\
&\triangle:=\{(x,y)\in\Yh\st x,y\in\R,\, 0\le x\le y\le 1\}
\quad\text{a triangle in \Yh},
\quad\text{and}\\
&\omega_0:=\frac{dx\wedge dy}{(1-x)\,y},
\end{align*}
thus getting
$$
\zeta(2)=\int_{\triangle}\omega,
$$
with  
$\omega_0\in\Gamma(X_0,\Omega^2_{X_0/\Q})$ and
$\del\triangle\subset\Dh\cup\{(0,0),(1,1)\}$,
see Figure \ref{exp4a}.

\begin{figure}[h]
\begin{center}
\epsfig{file=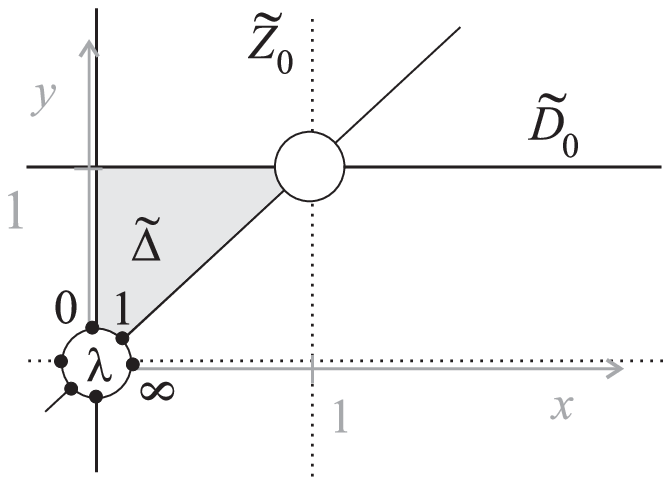,width=6cm}
\caption{The configuration $\wt{Z}_0,\wt{D}_0,\wt{\triangle}$}
\label{exp4b}
\end{center}
\vspace{-5mm}
\end{figure}

Now we blow up $Y_0$ in the points $(0,0)$ and $(1,1)$
obtaining $\pi_0 : \wt{Y}_0 \to Y_0$. We denote the strict transform
of $Z_0$ by $\wt{Z}_0$, $\pi_0^{\ast}\omega_0$ by $\wt{\omega}_0$ and
$\wt{Y}_0 \setminus \wt{Z}_0$ by $\wt{X}_0$. The ``strict transform''
$\ol{\pi^{-1}_{\rm an}(\triangle\setminus\{(0,0),(1,1)\})}$ will be called
$\wt{\triangle}$ and (being $\QR$-semi-algebraic hence triangulable ---
cf. Theorem \ref{triangulation})
gives rise to a singular chain
$$
\wt{\gamma}\in\h_2^{\rm sing}(\wt{X}^{\rm an}, \wt{D}^{\rm an};\Q).
$$
Since $\pi$ is an isomorphism away from the exceptional locus, this
exhibits $\zeta(2)$ as an (abstract) period (cf. Definition \ref{per2})
$$\zeta(2)=\int_{\triangle}\omega=\int_{\wt{\gamma}}\wt{\omega}\in\P_a=\P,
$$
see Figure \ref{exp4b}.

We will conclude this example by writing 
out $\wt{\omega}_0$ and $\wt{\triangle}$ more
explicitely. 
Note that $\wt{Y}_0$ can be described as the subvariety of
$$
\A^2_{\Q} \times \Q P^1 \times \Q P^1\quad
\text{with coordinates}\quad
(\wt{x},\wt{y};\lambda_0:\lambda_1;\mu_0:\mu_1)
$$
cut out by
$$
\wt{x}\lambda_0=\wt{y}\lambda_1\quad\text{and}\quad
(\wt{x}-1)\mu_0=(\wt{y}-1)\mu_1.
$$
With this choice of coordinates $\pi_0$ takes the form
$$
\begin{matrix}
\pi_0: & \wt{Y}_0 & \to & Y_0 \\
& (\wt{x},\wt{y};\lambda_0:\lambda_1;\mu_0:\mu_1) & \mapsto &
(\wt{x},\wt{y}) 
\end{matrix}
$$
and we have
$\wt{X}_0:=\wt{Y}_0\setminus(\{\lambda_0=0\}\cup\{\mu_1=0\})$.
We can embed $\wt{X}_0$ into affine space
\begin{align*}
\wt{X}_0 & \to \A_{\Q}^4 \\
(\wt{x},\wt{y};\lambda_0:\lambda_1;\mu_0:\mu_1) & 
\mapsto (\wt{x},\wt{y},\frac{\lambda_1}{\lambda_0},\frac{\mu_0}{\mu_1})
\end{align*}
and so have affine coordinates $\wt{x}$, $\wt{y}$,
$\lambda:=\frac{\lambda_1}{\lambda_0}$ and $\mu:=\frac{\mu_0}{\mu_1}$ 
on $\wt{X}_0$.

Now, on $\wt{X}_0\setminus\pi_0^{-1}(1,1)$, the form $\wt{\omega}_0$ is
given by
$$
\wt{\omega}_0 =
\frac{d\wt{x}\wedge d\wt{y}}{(1-\wt{x})\,\wt{y}} =
\frac{d(\lambda\wt{y})\wedge d\wt{y}}{(1-\wt{x})\,\wt{y}} =
\frac{d\lambda\wedge d\wt{y}}{1-\wt{x}},
$$
while on $\wt{X}_0\setminus\pi_0^{-1}(0,0)$
we have
$$
\wt{\omega}_0 =
\frac{d\wt{x}\wedge d\wt{y}}{(1-\wt{x})\,\wt{y}} =
\frac{d\wt{x}\wedge d(\wt{y}-1)}{(1-\wt{x})\,\wt{y}} =
\frac{d\wt{x}\wedge d(\mu(\wt{x}-1))}{(1-\wt{x})\,\wt{y}} =
\frac{-d\wt{x}\wedge d\mu}{\wt{y}}.
$$
The region $\wt{\triangle}$ is given by
$$
\wt{\triangle}=\{(\wt{x},\wt{y},\lambda,\mu)\in\wt{X}^{\rm an}\st
\wt{x},\wt{y},\lambda,\mu\in\R,\,\,\,
0\le \wt{x}\le \wt{y}\le 1,\,\,\,
0\le\lambda\le 1,\,\,\,
0\le\mu\le 1\}.
$$

\subsection{Fifth Example: The Double Logarithm Variation of Mixed
  Hodge Structures}

The previous example fits into a more general framework.

\subsubsection{Multiple Polylogarithms}
\index{Polylogarithm!multiple}
\label{polylogsubsec}

Define the {\em hyperlogarithm} 
\index{Hyperlogarithm}
as the iterated integral
$$
\ii_n(a_1,\ldots,a_n):=\int_{0\le t_1\le \cdots \le t_n\le 1}
\frac{dt_1}{t_1-a_1} \wedge \cdots \wedge \frac{dt_n}{t_n-a_n}
$$
with $a_1,\ldots,a_n\in\C$ 
(cf. \cite[p.~168]{zhao02}).
These integrals specialize to the
{\em multiple polylogarithm}
(cf. [{\em loc. cit.}]\footnote{%
Note that the factor $(-1)^n$ is missing in equation (1.2)
of \cite[p.~168]{zhao02}.})
$$
\Li_{m_1,\ldots,m_n}
  \left(\frac{a_2}{a_1},\cdots,\frac{a_n}{a_{n-1}},\frac{1}{a_n}\right):= 
(-1)^n \, \ii_{\sum  m_n}
(a_1,\underbrace{0,\ldots,0}_{m_1-1},
\ldots,
a_n,\underbrace{0,\ldots,0}_{m_n-1}),
$$
which is convergent if 
$1<|a_1|<\cdots<|a_n|$ (cf. \cite[2.3, p.~9]{goncharov:polylog}).
Alternatively, we can describe the multiple polylogarithm as
a power series
(cf. \cite[Thm.~2.2, p.~9]{goncharov:polylog})
\begin{equation}
\label{infinitesum}
\Li_{m_1,\ldots,m_n}(x_1,\ldots,x_n)=
\sum_{0<k_1<\cdots<k_n}
\frac{x_1^{k_1} \cdots x_n^{k_n}}{k_1^{m_1} \cdots k_n^{m_n}}
\quad\text{for}\quad\vert x_i\vert<1.
\end{equation}
Of special interest to us will be the {\em dilogarithm}
$\Li_2(x)=\sum_{k>0} \frac{x^k}{k^2}$ 
\index{Dilogarithm}
and the {\em double logarithm} 
\index{Logarithm!double}
$\Li_{1,1}(x,y)=\sum_{0<k<l} \frac{x^k y^l}{kl}$.
\label{dlogintro}

At first, the functions $\Li_{m_1,\ldots,m_n}(x_1,\ldots,x_n)$ only
make sense for $|x_i|<1$, but they can be analytically continued to
multivalued meromorphic functions on $\C^n$ (cf. \cite[p.~2]{zhao02}), for example
$\Li_1(x)=-\ln(1-x)$.

\paragraph{General Remarks.}
Polylogarithms have applications in Zagier's conjecture
(cf. \cite[p.~2]{zhao02}) and are used in the study 
of certain motivic fundamental groups 
(cf. \cite[p.~9]{huber_wildeshaus} and
\cite{goncharov:dihedral}).

Besides, equation (\ref{infinitesum}) would suggest a connection
between values of multiple polylogarithms and multiple $\zeta$-values
(cf. \cite[p.~617]{goncharov:dlog})
$$
\zeta(m_1,\ldots,m_n):=
\sum_{0<k_1<\ldots<k_m} \frac{1}{k_1^{m_1} \cdots k_n^{m_n}}.
$$
Indeed, we have formally
$$
-\zeta(2) = 
\frac{-\pi^2}{6}
= -\Li_2(1)
= \lim_{b\to 0} \lim_{a\to 1} \ii_{2}(a,b) 
= \lim_{b\to 0} \lim_{a\to 1} \Li_{1,1}\!\!\left(\frac{b}{a},\frac{1}{b}\right).
$$
However, such formal computations are problematic, since a
limit may diverge or at least depend on the direction in which the
limit is taken, because multiple polylogarithms are multivalued.
To deal with this problem, Zagier, Goncharov and others proposed
regularization procedures, for which we have to refer to
\cite[2.9--2.10, p.~18--24]{goncharov:polylog}.

\medskip
After these general remarks, 
we want to study a variation of mixed \Q-Hodge structures: 
the {\em double logarithm variation},
for which multiple polylogarithms appear as coefficients.
For that we will combine results in
\cite[p.~620f]{goncharov:dlog}, 
\cite[4, p.~32f]{kleinjung},
\cite{zhao02}, 
\cite{zhao03a} and 
\cite{zhao03b}.\footnote{Note that \cite{zhao03a} plus \cite{zhao03b} is just a detailed version of
\cite{zhao02}.}
By computing limit mixed \Q-Hodge structures, we will obtain
$-\zeta(2)$ as a ``period'' of a mixed \Q-Hodge structure.

\subsubsection{The Configuration}
\label{confsubsec}

\begin{figure}[h]
\begin{center}
\epsfig{file=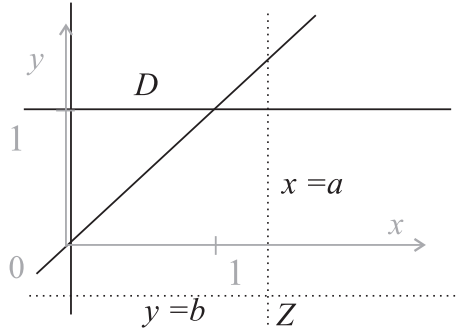,width=6cm}
\caption{The algebraic pair $(X,D)$}
\label{exp5a}
\end{center}
\vspace{-5mm}
\end{figure}

Let us consider the configuration
\begin{align*}
&Y:=\A^2_{\C}\quad\text{with coordinates $x$ and $y$,}\\
&Z:=\{x=a\}\cup\{y=b\}
  \quad\text{with}\quad a\ne 0,1
  \quad\text{and}\quad b\ne 0,1\\
&X:=Y\setminus Z\\
&D:=(\{x=0\}\cup\{y=1\}\cup\{x=y\})\setminus Z,
\end{align*}
see Figure \ref{exp5a}.

We denote the irreducible components of the divisor $D$ as follows:
\begin{align*}
& D_1:=\{x=0\}\setminus\{(0,b)\}, \\
& D_2:=\{y=1\}\setminus\{(a,1)\},\quad\text{and} \\
& D_3:=\{x=y\}\setminus\{(a,a),(b,b)\}.
\end{align*}
By projecting from $\Yh$ onto the $y$- or $x$-axis,
we get isomorphisms for the 
associated complex analytic spaces 
(cf. Subsection \ref{defcomplexanspace})
$$
\Dh_1\iso\C\setminus\{b\},\quad
\Dh_2\iso\C\setminus\{a\},\quad\text{and}\quad
\Dh_3\iso\C\setminus\{a,b\}.
$$

\subsubsection{Singular Homology} 
\label{exp5singsubsec}

\begin{figure}[b]
\begin{center}
\epsfig{file=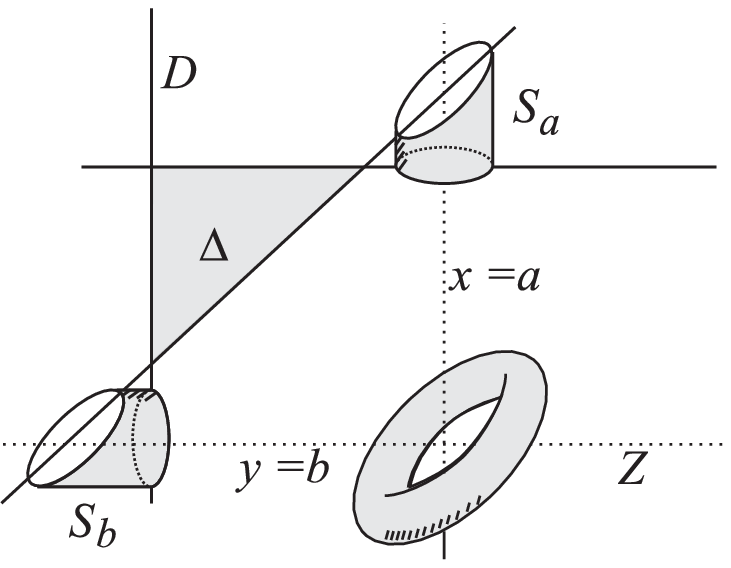,width=6cm}
\caption{Generators of  $\h^{\rm sing}_2(\Xh,\Dh;\Q)$}
\label{exp5b}
\end{center}
\end{figure}

We can easily give generators for the second singular homology of the
pair $(\Xh,\Dh)$, see Figure \ref{exp5b}.
\suppressfloats[t]
\begin{itemize}
\item
Let $\alpha:[0,1]\to\C$ be a smooth path, which does not meet $a$ or
$b$. We define a ``triangle''
$$
\triangle:=\{\bigl(\alpha(s),\alpha(t)\bigr)\st 0\le s\le t\le 1\}.
$$
\item
Consider the closed curve in \C
$$
C_b:=\left\{\frac{a}{b+\eps e^{2\pi i s}} \st s\in [0,1]\right\},
$$
which divides \C\ into two regions: an inner one containing
$\frac{a}{b}$ and an outer one.
We can choose $\eps>0$ small enough such that $C_b$
separates $\frac{a}{b}$ from $0$ to $1$, i.e. such that $0$ and $1$
are contained in the outer region.  This allows us to find a
smooth path $\beta:[0,1]\to \C$ from $0$ to $1$ not meeting $C_b$.
We define a ``slanted tube''
$$
S_b:=\bigl\{\bigl(\beta(t)\cdot(b+\eps e^{2\pi i s}),b+\eps e^{2\pi i s}\bigr) \st
s,t\in [0,1] \bigr\}
$$
which winds around $\{y=b\}$ and whose boundary components are
supported on $\Dh_1$ (corresponding to $t=0$) and $\Dh_3$
(corresponding to $t=1$).
The special choice of $\beta$ guarantees $S_b\cap\Zh=\emptyset$.
\item
Similarly, we choose $\eps>0$ such that the closed curve
$$
C_a:=\left\{\frac{b-1}{a-1-\eps e^{2\pi i s}} \st s\in [0,1]\right\}
$$
separates $\frac{b-1}{a-1}$ form $0$ and $1$. Let $\gamma:[0,1]\to\C$
be a smooth path from $0$ to $1$ which does not meet $C_a$. We have a
``slanted tube''
$$
S_a:=\bigl\{\bigl(a+\eps e^{2\pi i s},1+\gamma(t)\cdot(a+\eps e^{2\pi i s} -1)\bigr) 
\st s,t\in[0,1]\bigr\}
$$
winding around $\{x=a\}$ with boundary supported on $\Dh_2$ and
$\Dh_3$.
\item
Finally, we have a torus
$$
\T:=\{(a+\eps e^{2\pi i s},b+\eps e^{2\pi i t}) \st s,t\in[0,1]\}.
$$
\end{itemize}
The $2$-form $ds\wedge dt$ defines an orientation on the unit square
$[0,1]^2=\{(s,t)\st s,t\in[0,1]\}$. Hence the manifolds with boundary
$\triangle$, $S_b$, $S_a$, $\T$ inherit an orientation; and since they can
be triangulated, they give rise to smooth singular chains. By abuse
of notation we will also write $\triangle$, $S_b$, $S_a$, $\T$ for
these smooth singular chains. The homology classes of $\triangle$,
$S_b$, $S_a$ and $\T$ will be denoted by $\gamma_0$, $\gamma_1$,
$\gamma_2$ and $\gamma_3$, 
\label{gammas}
respectively.

\smallskip
An inspection of the long exact sequence in singular homology
will reveal that $\gamma_0,\ldots,\gamma_3$ form a system generators (see the following proof)
$$
\begin{CD}
\h^{\rm sing}_2(\Dh;\Q) @>>> 
\h^{\rm sing}_2(\Xh;\Q) @>>> 
\h^{\rm sing}_2(\Xh,\Dh;\Q) @>>> \\
\h^{\rm sing}_1(\Dh;\Q) @>i_1>> 
\h^{\rm sing}_1(\Xh;\Q) \, .
\end{CD}
$$

\begin{claim}
\label{exp5sing}

With notation as above, we have for the second singular homology of
the pair $(\Xh,\Dh)$
$$
\h^{\rm sing}_2(\Xh,\Dh;\Q)=
\Q\,\gamma_0\oplus
\Q\,\gamma_1\oplus
\Q\,\gamma_2\oplus
\Q\,\gamma_3.
$$
\end{claim}

\begin{proof}
For $c:=a$ and $c:=b$, the inclusion of the circle $\{c+\eps e^{2\pi i s}\st s\in[0,1]\}$
into $\C\setminus\{c\}$ is a homotopy equivalence,
hence the product map $\T\inclusion\Xh$ is also a homotopy equivalence.
This shows
$$
\h^{\rm sing}_2(\Xh;\Q)=\Q\,\T,
$$
while $\h^{\rm sing}_1(\Xh;\Q)$ has rank two with generators
\begin{itemize}
\item
one loop winding counterclockwise around $\{x=a\}$ once, but not
around $\{y=b\}$, thus being homologous to both $\del S_a \cap \Dh_2$
and $-\del S_a \cap \Dh_3$, and
\item
another loop winding counterclockwise around $\{y=b\}$ once, but not
around $\{x=a\}$, thus being homologous to $\del S_b\cap\Dh_1$ and
$-\del S_b\cap\Dh_3$.
\end{itemize}
In order to compute the Betti-numbers $b_i$ of $\Dh$, we use the spectral
sequence from Proposition \ref{specseqdivan} for the closed covering $\{\Dh_i\}$
\begin{center}
\vspace{-\baselineskip}
$$ 
E_2^{p,q}:=
\quad%
\begin{array}{cccccc}
\cdots & 0 & 0 & 0 & 0 & \cdots \\
\cdots & 0 & \bigoplus_{i=1}^3 \hDR{1}(\Dh_i;\C) & 0 & 0 & \cdots \\
\cdots & 0 & \ker\delta & \coker\delta  & 0 & \cdots\\
\cdots & 0 & 0 & 0 & 0 & \cdots 
\end{array}
\quad\Rightarrow\quad
E^n_{\infty}:=\hDR{n}(\Dh;\C),
$$
\vspace{-\baselineskip}
\end{center}
where
$$
\delta:\bigoplus_{i=1}^3 \hDR{0}(\Dh_i;\C) \longto
\bigoplus_{i<j}\hDR{0}(\Dh_{ij};\C).
$$
Note that this spectral sequence degenerates. Since $\Dh$ is
connected, we have $b_0=1$, i.e. 
$$
1=
b_0=
\rank_{\C} E^0_{\infty}=
\rank_{\C} E^{0,0}_2=
\rank_{\C}\ker\delta.
$$
Hence
\begin{align*}
\rank_{\C}\coker\delta&=
\rank_{\C}{\rm codomain}\,\delta-\rank_{\C}{\rm domain}\,\delta+\rank_{\C}\ker\delta\\
&=(1+1+1)-(1+1+1)+1=1,
\end{align*}
and so
\begin{align*}
b_1&=
\rank_{\C} E^1_{\infty}=
\rank_{\C} E^{1,0}_2 + \rank_{\C} E^{0,1}_2\\
&=\sum_{i=1}^3\,\rank_{\C}\hDR{1}(\Dh_i;\C)+\rank_{\C}\coker\delta\\
&=\rank_{\C} \h^1(\C\setminus\{b\};\C)+
\rank_{\C} \h^1(\C\setminus\{a\};\C)+
\rank_{\C} \h^1(\C\setminus\{a,b\};\C)+1\\
&=(1+1+2)+1=5.
\end{align*}
We can easily specify generators of $\h^{\rm sing}_1(\Dh;\Q)$ as
follows
$$
\h^{\rm sing}_1(\Dh;\Q)=
\Q\cdot(\del S_b\cap\Dh_1)\oplus\,
\Q\cdot(\del S_a\cap\Dh_2)\oplus\,
\Q\cdot(\del S_b\cap\Dh_3)\oplus\,
\Q\cdot(\del S_a\cap\Dh_3)\oplus\,
\Q\cdot\del\triangle.
$$
Clearly $b_2=\rank_{\C}\h^{\rm sing}_2(\Dh;\Q)=0$. 
Now we can compute $\ker i_1$ and obtain
$$
\ker i_1=\Q\cdot\del\triangle\oplus
\Q\cdot(\del S_b\cap\Dh_1+\del S_b\cap\Dh_3)\oplus
\Q\cdot(\del S_a\cap\Dh_2+\del S_a\cap\Dh_3).
$$
This shows finally
$$
\rank_{\Q}\h^{\rm sing}_2(\Xh,\Dh;\Q)=
\rank_{\Q}\h^{\rm sing}_2(\Xh;\Q)+\rank_{\Q}\ker i_1=
1+3=4.
$$
From these explicit calculations we also derive the linear
independence of $\gamma_0=[\triangle]$, $\gamma_1=[S_b]$,
$\gamma_2=[S_a]$, $\gamma_3=[\T]$ and Claim \ref{exp5sing} is proved.
\end{proof}

\subsubsection{Smooth Singular Homology}

\newcommand{\tuple}[7]{\ensuremath{(
\underset{\emptyset}{#1},
\underset{1}{#2},
\underset{2}{#3},
\underset{3}{#4},
\underset{12}{#5},
\underset{13}{#6},
\underset{23}{#7}
)}}
Recall the definition of smooth singular cohomology (cf. Definition \ref{defsmooth2}).
With the various sign conventions made so far 
(cf. Subsection \ref{smoothsubsec} and the appendix), 
the boundary map
$\delta:\c^{\infty}_2(\Xh,\Dh;\Q)\to\c^{\infty}_1(\Xh,\Dh;\Q)$ is
given by
$$
\delta: 
\c^{\infty}_2(\Xh;\Q)\oplus
\bigoplus_{i=1}^3\c^{\infty}_1(\Dh_i;\Q)\oplus
\bigoplus_{i<j} \c^{\infty}_0(\Dh_{ij};\Q)
\to
\c^{\infty}_1(\Xh;\Q)\oplus
\bigoplus_{i=1}^3(\Dh_i;\Q) 
$$
\vspace{-7mm}
\begin{multline*}
\hspace{3mm}
\tuple{\sigma_{\phantom{1}}\!\!\!}{\sigma_1}{\sigma_2}{\sigma_3}{\sigma_{12}}{\sigma_{13}}{\sigma_{23}}
\mapsto\\
\hspace{1.5cm}
(
\underset{\emptyset}{\del\sigma+\sigma_1+\sigma_2+\sigma_3},
\underset{1}{-\del\sigma_1+\sigma_{12}+\sigma_{13}},
\underset{2}{-\del\sigma_2-\sigma_{12}+\sigma_{23}},
\underset{3}{-\del\sigma_3-\sigma_{13}-\sigma_{23}}
).
\end{multline*}
Thus the following elements of $\c^{\infty}_2(\Xh,\Dh;\Q)$ are cycles
\begin{itemize}
\item
$\Gamma_0:=
\tuple{\triangle}{-\del\triangle\cap\Dh_1}{-\del\triangle\cap\Dh_2}{-\del\triangle\cap\Dh_3}%
{\Dh_{12}}{-\Dh_{13}}{\Dh_{23}}$,
\item
$\Gamma_1:=
\tuple{S_b}{-\del S_b\cap\Dh_1}{0}{-\del S_b\cap\Dh_3}{0}{0}{0}$,
\item
$\Gamma_2:=
\tuple{S_a}{0}{-\del S_a\cap\Dh_2}{0}{-\del S_a\cap\Dh_3}{0}{0}$
and
\item
$\Gamma_3:=
\tuple{\T}{0}{0}{0}{0}{0}{0}$.
\end{itemize}
Under the isomorphism
$\h^{\infty}_2(\Xh,\Dh;\Q)
\stackrel{\sim}{\longto}
\h^{\rm sing}_2(\Xh,\Dh;\Q)$ 
the classes of these cycles
$[\Gamma_0]$, $[\Gamma_1$], $[\Gamma_2]$, $[\Gamma_3]$ 
are mapped to 
$\gamma_0$, $\gamma_1$, $\gamma_2$, $\gamma_3$, 
respectively (cf. Proposition \ref{compplus}).

\subsubsection{Algebraic \DeRham Cohomology of $(\Xh,\Dh)$}
\label{exp5drsubsec}

Recall the definition of the complex $\wom{X,D/\C}$ on page \pageref{defwomrel}.
We consider 
$$
\Gamma(X;\wt{\Omega}^2_{X,D/\C})=
\Gamma(X;\Omega^2_{X/\C})\oplus
\bigoplus_{i=1}^3 \Gamma(D_i;\Omega^1_{D_i/\C})\oplus
\bigoplus_{i<j} \Gamma(D_{ij};\OO_{D_{ij}})
$$
together with the following cycles of $\Gamma(X;\wt{\Omega}^2_{X,D/\C})$ 

\begin{itemize}
\item
$\omega_0:=
\tuple{\frac{dx\wedge dy}{(x-a)(y-b)}}{0}{0}{0}{0}{0}{0}$,
\item
$\omega_1:=
\tuple{0}{\frac{-dy}{y-b}}{0}{0}{0}{0}{0}$,
\item
$\omega_2:=
\tuple{0}{0}{\frac{-dx}{x-a}}{0}{0}{0}{0}$, and
\item
$\omega_3:=
\tuple{0}{0}{0}{0}{0}{0}{1}$.
\end{itemize}
By computing the (transposed) period matrix
$P_{ij}:=\langle\Gamma_j,\omega_i\rangle$ 
(cf. Definition \ref{per1}) 
and checking its
non-degeneracy, we will show that $\omega_0$, $\ldots$, $\omega_3$
span $\hDR{2}(X,D/\C)$.

\begin{claim}
\label{exp5dr}
Let $X$ and $D$ be as in Subsection \ref{confsubsec}. Then the second
algebraic \deRham cohomology group 
$\hDR{2}(X,D/\C)$
of the pair $(X,D)$ is generated by the cycles
$\omega_0,\ldots,\omega_3$ considered above.
\end{claim}

\subsubsection{Period Matrix of $(X,D)$}
Easy calculations give us for the (transposed) period matrix
$P:=(\langle\Gamma_j,\omega_i\rangle)_{i,j=1}^3$
(cf. Definition \ref{per1})
\begin{center}
\begin{tabular}{c|cccc}
& $\Gamma_0$ & $\Gamma_1$ & $\Gamma_2$ & $\Gamma_3$ \\
\hline 
$\omega_0$ & $1$ & $0$ & $0$ & $0$ \\
$\omega_1$ & $\Li_1(\frac{1}{b})$ & $2\pi i$ & $0$ & $0$ \\
$\omega_2$ & $\Li_1(\frac{1}{a})$ & $0$ & $2\pi i$ & $0$ \\
$\omega_3$ & ? & $2\pi i\Li_1(\frac{b}{a})$ & 
  $2\pi i\ln\!\left(\frac{a-b}{1-b}\right)$ & $(2\pi i)^2$;
\end{tabular}
\end{center}
for example (cf. Subsection \ref{alteriso})
\begin{itemize}
\item
$P_{1,1}=\langle\Gamma_1,\omega_1\rangle
=\int_{-\del S_b\cap \Dh_1} \frac{-dy}{y-b} \\[1mm]
\phantom{P_{1,1}=\langle\Gamma_1,\omega_1\rangle}
=\int_{|y-b|=\eps}\frac{dy}{y-b} \\[1mm]
\phantom{P_{1,1}=\langle\Gamma_1,\omega_1\rangle}
=2\pi i$,
\item
$P_{3,3}=\langle\Gamma_3,\omega_3\rangle
=\int_{\T}\frac{dx}{x-a}\wedge\frac{dy}{y-b} \\[1mm]
\phantom{P_{3,3}=\langle\Gamma_3,\omega_3\rangle}
=\left(\int_{|x-a|=\eps}\frac{dx}{x-a}\right)\cdot
\left(\int_{|y-b|=\eps}\frac{dy}{y-b}\right)\quad\text{by Fubini} \\[1mm]
\phantom{P_{3,3}=\langle\Gamma_3,\omega_3\rangle}
=(2\pi i)^2$,
\item
$P_{1,0}=\langle\Gamma_0,\omega_1\rangle
=\int_{-\del\triangle\cap\Dh_1}\frac{-dy}{y-b} \\[1mm]
\phantom{P_{1,0}=\langle\Gamma_0,\omega_1\rangle}
=\int_0^1\frac{-\alpha(t)}{\alpha(t)-b} \\[1mm]
\phantom{P_{1,0}=\langle\Gamma_0,\omega_1\rangle}
=-[\ln(\alpha(t))-b]_0^1 \\
\phantom{P_{1,0}=\langle\Gamma_0,\omega_1\rangle}
=-\ln\!\left(\frac{1-b}{-b}\right) \\[1mm]
\phantom{P_{1,0}=\langle\Gamma_0,\omega_1\rangle}
=-\ln\!\left(1-\frac{1}{b}\right) \\[1mm]
\phantom{P_{1,0}=\langle\Gamma_0,\omega_1\rangle}
=\Li_1\!\left(\frac{1}{b}\right)$,\quad and
\item
$P_{3,1}=\langle\Gamma_1,\omega_3\rangle
=\int_{S_b}\frac{dx}{x-a}\wedge\frac{dy}{y-b} \\[1mm]
\phantom{P_3,1=\langle\Gamma_1,\omega_3\rangle}
=\int_{[0,1]^2}
  \frac{d\left(\beta(t)\cdot\left(b+\eps e^{2\pi i s}\right)\right)}
    {\beta(t)\cdot\left(b+\eps e^{2\pi i s}\right)-a}
  \wedge
  \frac{d\left(b+\eps e^{2\pi i s}\right)}
    {\eps e^{2\pi i s}} \\[1mm]
\phantom{P_3,1=\langle\Gamma_1,\omega_3\rangle}
=\int_{[0,1]^2}
  \frac{b+\eps e^{2\pi i s}}
    {\beta(t)\cdot\left(b+\eps e^{2\pi i s}\right)-a}
  d\beta(t)
  \wedge
  2\pi i
  ds \\[1mm]
\phantom{P_3,1=\langle\Gamma_1,\omega_3\rangle} 
=-{\displaystyle\int_0^1}
  \left[
  \frac{a\ln\left(\beta(t)
      \cdot
      \left(b+\eps e^{2\pi i s}\right)-a\right)
    -2\pi i\beta(t)b s}
  {\beta(t)\cdot\left(-\beta(t)b+a\right)}
  \right]_0^1
\,d\beta(t) \\[1mm]
\phantom{P_3,1=\langle\Gamma_1,\omega_3\rangle}
=-2\pi i\int_0^1
  \frac{d\beta(t)}
  {\beta(t)-\frac{a}{b}} \\
\phantom{P_3,1=\langle\Gamma_1,\omega_3\rangle}
=-2\pi i\left[\ln\!\left(\beta(t)-\frac{a}{b}\right)\right]_0^1 \\[1mm]
\phantom{P_3,1=\langle\Gamma_1,\omega_3\rangle}
=-2\pi i\ln\!\left(\frac{1-\frac{a}{b}}{-\frac{a}{b}}\right) \\[1mm]
\phantom{P_3,1=\langle\Gamma_1,\omega_3\rangle}
=-2\pi i\ln\!\left(1-\frac{a}{b}\right) \\[1mm]
\phantom{P_3,1=\langle\Gamma_1,\omega_3\rangle}
=2\pi i\,\Li_1\!\left(\frac{b}{a}\right)$.
\end{itemize}

Obviously the period matrix $P$ is non-degenerate and so Claim
\ref{exp5dr} is proved.

What about the entry $P_{3,0}$? We want to show that
$\langle\Gamma_0,\omega_3\rangle=\Li_{1,1}\!\left(\frac{b}{a},\frac{1}{b}\right)$, 
where $\Li_{1,1}(x,y)$ is an analytic continuation of the double
logarithm defined for $|x|, |y|<1$ in Subsection \ref{polylogsubsec}.

\subsubsection{Analytic Continuation  of the Double Logarithm
  $\Li_{1,1}(x,y)$}

We describe this analytic continuation in detail. 
Our approach is similar to the one taken in 
\cite[2.3, p.~9]{goncharov:polylog}, 
but differs from that in \cite[p.~7]{zhao03a}.

Let $\Bh:=(\C\setminus\{0,1\})^2$ be the parameter space and choose a
point
$(a,b)\in\Bh$. For $\eps>0$ we denote by $D_{\eps}(a,b)$ the
polycylinder
$$
D_{\eps}(a,b):=\{(a,b)\in\Bh\st|a'-a|<\eps, |b'-b|<\eps\}.
$$
If $\alpha:[0,1]\to\C$ is a smooth path from $0$ to $1$ passing
through neither $a$ nor $b$, then there exists an $\eps>0$ such that
$\im\alpha$ does not meet any of the discs
\begin{align*}
& D_{2\eps}(a):=\{a'\in\C\st |a'-a|<2\eps\},\quad\text{and} \\
& D_{2\eps}(b):=\{b^{\,\prime}\in\C\st |b^{\,\prime}-b|<2\eps\}.
\end{align*}
Hence the power series (\ref{powerseries}) below
\begin{align}
\frac{1}{\alpha(s)-a'}
\frac{1}{\alpha(t)-b^{\,\prime}}&=
\frac{1}{\alpha(s)-a}\frac{1}{1-\frac{a'-a}{\alpha(s)-a}}
\frac{1}{\alpha(t)-b}\frac{1}{1-\frac{b^{\,\prime}-b}{\alpha(t)-b}} \notag \\
&=\sum_{k,l=0}^{\infty}
\underbrace{\frac{1}{(\alpha(s)-a)^{k+1}(\alpha(t)-b)^{l+1}}}_{c_{k.l}}
(a'-a)^k(b^{\,\prime}-b)^l \label{powerseries}
\end{align}
has coefficients $c_{k,l}$ satisfying
$$
|c_{k,l}|<\left(\frac{1}{2\eps}\right)^{k+l+2}.
$$
In particular, (\ref{powerseries}) converges uniformly for $(a',b^{\,\prime})\in
D_{\eps}(a,b)$ and we see that the integral
\label{intalpha}
\begin{align*}
\ii^{\alpha}_2(a',b^{\,\prime}):=&
\int_{0\le s\le t\le 1}
\frac{d\alpha(s)}{\alpha(s)-a'}\wedge
\frac{d\alpha(t)}{\alpha(t)-b^{\,\prime}} \\
=&\sum_{k,l=0}
\left(\int_{0\le s\le t\le 1}
\frac{d\alpha(s)}{(\alpha(s)-a)^{k+1}}\wedge
\frac{d\alpha(t)}{(\alpha(t)-b)^{l+1}}\right)
(a'-a)^k(b^{\,\prime}-b)^l
\end{align*}
defines an analytic function of $D_{\eps}(a,b)$.
In fact, by the same argument we get an analytic function
$\ii^{\alpha}_2$ on all of $(\C\setminus\im\alpha)^2$.

Now let 
$\alpha_r:[0,1]\to\C\setminus\left(D_{2\eps}(a)\cup
D_{2\eps}(b)\right)$
with $r\in[0,1]$ be a smooth homotopy of paths from $0$ to $1$,
i.e. $\alpha_r(0)=0$ and $\alpha_r(1)=1$ for all $r\in[0,1]$. 
We show
$$
\ii^{\alpha_0}_2(a',b^{\,\prime})=\ii^{\alpha_1}_2(a',b^{\,\prime})\quad\text{for all}\quad
(a',b^{\,\prime})\in D_{\eps}(a,b).
$$
Define a subset $\Gamma\subset\C^2$
$$
\Gamma:=\{(\alpha_r(s),\alpha_r(t))\st 0\le s\le t\le 1, r\in[0,1]\}.
$$
The boundary of $\Gamma$ is built out of five components (each being a
manifold with boundary)
\begin{itemize}
\item
$\Gamma_{s=0}:=\{(0,\alpha_r(t))\st r,t\in[0,1]\}$,
\item
$\Gamma_{s=t}:=\{(\alpha_r(s),\alpha_r(s))\st r,s\in[0,1]\}$,
\item
$\Gamma_{t=1}:=\{(\alpha_r(s),1)\st r,s\in[0,1]\}$,
\item
$\Gamma_{r=0}:=\{(\alpha_0(s),\alpha_0(t)\st 0\le s\le t\le 1\}$,
\item
$\Gamma_{r=1}:=\{(\alpha_1(s),\alpha_1(t)\st 0\le s\le t\le 1\}$.
\end{itemize}
Let $(a',b^{\,\prime})\in D_{\eps}(a,b)$. Since the restriction of 
$\frac{dx}{x-a'}\wedge\frac{dy}{y-b^{\,\prime}}$ to $\Gamma_{s=0}$,
$\Gamma_{s=t}$ and $\Gamma_{t=1}$ is zero, we get by Stoke's theorem
\begin{align*}
0=\int_{\Gamma}0
&=\int_{\Gamma}d \frac{dx}{x-a'}\wedge\frac{dy}{y-b^{\,\prime}} \\
&=\int_{\del\Gamma} \frac{dx}{x-a'}\wedge\frac{dy}{y-b^{\,\prime}} \\
&=\int_{\Gamma_{r=1}-\Gamma_{r=0}}\frac{dx}{x-a'}\frac{dy}{y-b^{\,\prime}} \\
&=\ii^{\alpha_1}_2(a',b^{\,\prime})-\ii^{\alpha_0}_2(a',b^{\,\prime}).
\end{align*}

For each pair of smooth paths $\alpha_0,\alpha_1:[0,1]\to\C$ from $0$
to $1$, we can find a homotopy $\alpha_r$ relative to $\{0,1\}$
between both paths. Since $\im\alpha_r$ is compact, we also find a
point $(a,b)\in\Bh=(\C\setminus\{0,1\})^2$ 
and an $\eps>0$ such that $\im\alpha_r$ does not
meet $D_{2\eps}(a,b)$ or $D_{2\eps}(a,b)$. Then $\ii^{\alpha_0}_2$ and
$\ii^{\alpha_1}_2$ must agree on $D_{\eps}(a,b)$. By the identity
principle for analytic functions of several complex variables \cite[A,
  3, p.~5]{gunning}, the functions $\ii^{\alpha}_2(a',b^{\,\prime})$, each
defined on $(\C\setminus\im\alpha)^2$, patch together to give a
multivalued analytic function on $\Bh=(\C\setminus\{0,1\})^2$.

Now assume $1<|b|<|a|$, then we can take $\alpha=\id:[0,1]\to\C,\,
s\mapsto s$, and obtain
$$
\ii^{\id}_2(a,b)=\ii_2(a,b)=\Li_{1,1}\!\left(\frac{b}{a},\frac{1}{y}\right),
$$
where $\Li_{1,1}(x,y)$ is the double logarithm defined for
$|x|,|y|<1$ in Subsection \ref{polylogsubsec}.
Thus we have proved the following lemma.

\pagebreak

\begin{lemma}
\label{dlog}
The integrals 
$$
\ii^{\alpha}_2\!\left(\frac{1}{xy},\frac{1}{y}\right)=
\int_{0\le s\le t\le 1}
\frac{d\alpha(s)}{\alpha(s)-\frac{1}{xy}}\wedge
\frac{d\alpha(t)}{\alpha(t)-\frac{1}{y}}
$$
with $\alpha:[0,1]\to\C$ a smooth path from $0$ to $1$, and
$\frac{1}{xy},\frac{1}{b}\in\C\setminus\im\alpha$, 
defined above on page \pageref{intalpha}, provide a
genuine analytic continuation of 
$\Li_{1,1}\!\left(x,y\right)$ to a multivalued function
which is defined on 
$\{(x,y)\in\C^2\st x,y\ne0, xy\ne 1, y\ne 1\}$.
\end{lemma}

\begin{definition}[Double logarithm]
\label{li11}
We call the analytic continuation from Lemma \ref{dlog} the {\em
double logarithm} as well and continue to use the notation
$\Li_{1,1}(x,y)$.
\end{definition}

This gives us finally $P_{3,0}=\Li_{1,1}\!\left(\frac{b}{a},\frac{1}{b}\right).$

\subsubsection{Properties of the Double and the Dilogarithm}

Similarly, one can analytically continue the dilogarithm $\Li_2(x)$ to
a multivalued meromorphic function on \C.

Since analytic continuations are uniquely determined (as multivalued
functions) our Definition \ref{li11} must agree with the ones given by
Goncharov in \cite[p.~9]{goncharov:polylog} or by Zhao in \cite[4, p.~4f]{zhao03b}. Hence we are allowed to cite
various facts concerning the double and the dilogarithm.

\begin{lemma}[Monodromy of the double logarithm, {\cite[Thm. 5.3,
	p.~11]{zhao03a}}]
\label{dlogmono}

Let $\sigma=\{1+\eps e^{2\pi i s},b)\st s\in[0,1]\}$ 
be a loop in $\Bh=(\C\setminus\{0,1\})^2$ winding
counterclockwise around $\{a=1\}$ once, but not around $\{a=0\}$,
$\{b=0\}$ or $\{b=1\}$. Then
\begin{align*}
& \int_{\sigma} d\Li_{1}\!\left(\frac{1}{a}\right)=-2\pi
  i,\quad\text{and} \\
& \int_{\sigma}d\Li_{1,1}\!\left(\frac{b}{a},\frac{1}{b}\right)=
  -2\pi i\ln\!\left(\frac{1-b+\eps}{1-b}\right).
\end{align*}
\end{lemma}

\begin{lemma}[Relation with the dilogarithm,
  {\cite[p.~620]{goncharov:dlog}}]
\label{dlogdilog}
We have the following identity of multivalued functions
$$
\Li_{1,1}\!\left(b,\frac{1}{b}\right)=-\Li_2\!\left(\frac{1}{1-b}\right).
$$
\end{lemma}

\begin{lemma}[Monodromy of the dilogarithm, {\cite[Prop.~2.2, p.~7]{hain}}]
\label{dilogmono}
If $\sigma=\{1+\eps e^{2\pi i s}\st s\in[0,1]\}$ is a circle in
$\C\setminus\{1\}$ of radius $\eps$ winding counterclockwise
around $1$ once, then
$$
\int_{\sigma}d\Li_{2}(c)=2\pi i\,\Li_1(-\eps).
$$
\end{lemma}

\subsubsection{The Mixed \Q-Hodge Structure on
  $\h^2\sing(\Xh,\Dh;\Q)$}

According to Deligne (cf. \cite[8.3.8, p.~43]{hodge3}) the cohomology group
$\h^2\sing(\Xh,\Dh;\Q)$ carries a {\em mixed \Q-Hodge structure},
short \Q-MHS,%
\index{Mixed \Q-Hodge structure}%
\index{Q-MHS@\Q-MHS|see{Mixed \Q-Hodge structure}}%
\notation{Q}{\Q-MHS}{abbreviation  for mixed \Q-Hodge structure}%
consisting of 
\begin{itemize}
\item
an increasing {\em weight filtration} $W_{\Bul}$ on
$\h^2\sing(\Xh,\Dh;\Q)$ and 
\item
a decreasing {\em Hodge filtration}
$F\bul$ on $\h^2\sing(\Xh,\Dh;\C)$.
\end{itemize}
Deligne's {\em mixed Hodge theory} takes up several publications,
\cite{hodge2} and \cite{hodge3} being of central importance.
For a brief introduction, see \cite[Ch. 6]{gelfand_manin} and
\cite{durfee}.

Using the same ideas as in \cite[4.2, p.~38f]{kleinjung}, we can
compute this \Q-MHS and obtain the following lemma.

\begin{lemma}
\label{MHS}

Let $X$ and $D$ be as in Subsection \ref{confsubsec}
and
let $\Xh$, $\Dh$ be 
the associated complex analytic spaces 
(cf. Subsection \ref{defcomplexanspace}).

The cohomology group 
$\h^2\sing(\Xh,\Dh;\Q)$ 
car\-ries a \Q-MHS
with weight fil\-tra\-tion $W_{\Bul}$ gi\-ven by
\begin{align*}
&W_4  =\h^2\sing(\Xh,\Dh;\Q)
= \Q\gamma_0^{\ast}\oplus\Q\gamma_1^{\ast}\oplus\Q\gamma_2^{\ast}\oplus\Q\gamma_3^{\ast} \supset \\
&W_3  =W_2=\Q\gamma_0^{\ast}\oplus\Q\gamma_1^{\ast}\oplus\Q\gamma_2^{\ast}\supset \\
&W_1  =W_0=\Q\gamma_0^{\ast}\supset \\
&W_{-1}  =0,
\end{align*}
where the $\gamma_j^{\ast}$ are the duals of the $\gamma_j$,
$j=1,\ldots,4$ defined in Subsection \ref{exp5singsubsec},
and Hodge filtration $F\bul$ given by
\begin{align*}
F^0 & =\h^2\sing(\Xh,\Dh;\C)
= \C\omega_0\oplus\C\omega_1\oplus\C\omega_2\oplus\C\omega_3 \supset \\
F^1 & =\C\omega_1\oplus\C\omega_2\oplus\C\omega_3\supset \\
F^2 & =\C\omega_3\supset \\
F^3 & =0,
\end{align*}
where we used the cycles $\omega_0,\ldots,\omega_3$ from Subsection \ref{exp5drsubsec}
and identified 
$$
\hDR{2}(\Xh,\Dh;\C)\iso\h^2\sing(\Xh,\Dh;\C)
$$ 
using the Comparison Theorem \ref{comp}.
\end{lemma}

We illustrate the proof of Lemma \ref{MHS} with two exemplary
calculations.
\begin{itemize}
\item
We want to show that $\gamma_3^{\ast}\not\in
W_3\h^2\sing(\Xh,\Dh;\Q)$. From the inclusion $\T\inclusion\Xh$, we
get an element $[\T]^{\ast}\in\h^2\sing(\Xh,\Dh;\Q)$, which is the
preimage of $\gamma_3^{\ast}$ under the connecting morphism
$$
\h^2\sing(\Xh;\Q) \longto \h^2\sing(\Xh,\Dh;\Q).
$$
This map induces a map of \Q-MHS (cf. \cite[Ch.~6, Thm.~2.4,
p.~147]{gelfand_manin}), which is strictly compatible with the weight
and Hodge filtrations (cf. \cite[Ch.~6, 1.3, p.~143]{gelfand_manin}). In
particular, it will be sufficient to show
$$
[\T]^{\ast}\not\in W_3\h^2\sing(\Xh;\Q).
$$
We consider first
$$
U_c:=\A^1_{\C}\setminus\{c\}\quad\text{for}\quad c\in\C\quad
\text{with standard coordinate $t$}.
$$
Let $\sigma_c$ be the loop $\{c+\eps e^{2\pi i s}\st s\in[0,1]\}$
winding around $c$ and $\sigma_c^{\ast}$ its dual.
Since the weight-one-part of $\h^1\sing(U_c;\Q)$ is zero 
(cf. 
\cite[Ch.~6, 2.5.1, p.~147]{gelfand_manin})
\label{Uhb}
$$
W_1\h^1\sing(\Uh_c;\Q)=
\h^1\sing(\ol{U}_c^{\rm an};\Q)=
\h^1\sing(\C \Ph^1;\Q)=0,
$$
we see $\sigma_c^{\ast}\in W_2\h^1\sing(\Uh_c;\Q)$, but
$\sigma_c^{\ast}\not\in W_1\h^1\sing(\Uh_c;\Q)$,
i.e. $\sigma_c^{\ast}$ has weight two. 
We have $X=U_a\times U_b$ and the K{\"u}nneth isomorphism
(cf. \cite[Ch.~VI, Thm.~1.6, p.~320]{bredon})
$$
\h^2\sing(\Xh;\Q)=
\bigoplus_{k+l=2}
\h^k\sing(\Uh_a;\Q) \tensor{} \h^l\sing(\Uh_b;\Q)
$$
maps $[\T]^{\ast}$ to $\sigma_a^{\ast}\tensor{}\sigma_b^{\ast}$. The
K{\"u}nneth isomorphism is an isomorphism of \Q-MHS
(cf. \cite[Prop.~8.2.10, p.~40]{hodge3}, 
\cite[Ch.~6, Thm.~2.1, p.~145]{gelfand_manin}),
hence $[\T]^{\ast}$ has weight four,
i.e. $[\T]^{\ast}\not\in W_3\h^2\sing(\Xh;\Q)$.
\item
We want to show $[\omega_1]\in F^1\h^1\sing(\Xh,\Dh;\Q).$
The connecting morphism of the long exact sequence in analytic \deRham
cohomology \ref{longexseqandr}
$$
\delta: \hDR{1}(\Dh;\C) \longto \hDR{2}(\Xh,\Dh;\C)
$$
is indeed a morphism of \Q-MHS (see \cite[Prop.~8.3.9, p.~43]{hodge3}, 
\cite[Ch.~6, Thm.~2.4, p.~147]{gelfand_manin});
and so we have a strict morphism of filtered complexes 
(cf. \cite[Ch.~6, 1.3, p.~143]{gelfand_manin})
$$
\delta: (\hDR{1}(\Dh;\C),F\bul) \longto (\hDR{2}(\Xh,\Dh;\C),F\bul).
$$
Now the class of
$$
\omega'_1:=\left(\underset{1}{\frac{-dy}{y-b}},\underset{2}{0},\underset{3}{0}\right)\in
\Gamma(\Dh;\wt{\Omega}^1_{\Dh})=
\bigoplus_{i=1}^3\Gamma(\Dh_i;\Omega^1_{\Dh_i})
$$
is mapped to $[\omega_1]$ under $\delta$; hence it suffices to prove
$$
[\omega'_1]\in F^1\h^1\sing(\Dh;\C).
$$
Using a similar argument for the morphism of \Q-MHS 
(cf. \cite[Ch.~6, Thm.~2.1, p.~145]{gelfand_manin})
$$
i^{\ast}: 
(\h^1\sing(\Dh;\Q),W_{\Bul},F\bul) \longto
(\h^1\sing(\Dh_1;Q),W_{\Bul},F\bul)
$$ 
induced by the inclusion $i: \Dh_1\inclusion \Dh$,
we reduce to the claim
$$
\left[\frac{-dy}{y-b}\right]\in F^1 \h^1\sing(\Dh_1;\C).
$$
From the map
\begin{align*}
j: D_1 & \inclusion \C P^1 \\
(0,y) & \mapsto [1:y-b]
\end{align*}
we get an isomorphism
$$
j: D_1\iso\C P^1\setminus\{[1:0],[0:1]\}.
$$
The divisor $[1:0]+[0:1]$ on $\C P^1$ is just the inverse image of the
very ample divisor $[1:0:0]$ on $\C P^2$ under the $2$-uple embedding 
by \cite[Ex.~7.6.1, p.~155]{hartshorne},
hence very ample itself. This allows us to apply Proposition 9.2.4
from \cite[p.~49]{hodge3}, which states that:
\begin{quote}
\em For $Z$ a very ample divisor on a smooth projective variety $Y$ over
\C\ and $\omega$ a closed differential $q$-form on $Y\setminus Z$, 
we have $[\omega]\in F^p\h^q\sing(\Yh\setminus\Zh;\C)$ if
and only if $\omega$ has a pole of order less or equal $q-p+1$ along $Z$.
\end{quote}
Thus we obtain
$$
j_{\ast}\left[\frac{-dy}{y-b}\right]\in
F^1\h^1\sing(\C \Ph^1\setminus\{[1:0],[0:1]\};\C),
$$
i.e. 
$$
\left[\frac{-dy}{y-b}\right]\in F^1\h^1\sing(\Dh_1;\C),
$$
and our claim is proved.

Alternatively, we could have argued that $D_1$ is isomorphic to $U_b$
(as defined in the previous calculation)
\begin{align*}
j:D_1\stackrel{\sim}{\longto} &\,U_b \\
(0,y) \mapsto &\,y.
\end{align*}
Since $\h^1\sing(\Uh_b;\Q)$ carries a pure \Q-Hodge structure (see
page \pageref{Uhb}) and the corresponding Hodge decomposition enjoys
the symmetry 
(cf. \cite[p.~116]{griffiths_harris})
$$
\h^q(\Uh_b;\Omega^p_{\Uh_b}) \iso \h^p(\Uh_b;\Omega^q_{\Uh_b}),
$$
it follows from the fact that $\hDR{1}(\Uh_b;\C)$ is one-dimensional
with generator $j_{\ast}\left[\frac{-dy}{y-b}\right]$, that
$\frac{-dy}{y-b}$ has type $(1,1)$. 
\end{itemize}

For aesthetic reasons, we will dualize the \Q-MHS on
$\h^2\sing(\Xh,\Dh;\Q)$; thus getting a \Q-MHS on the homology group
$\h_2^{\rm sing}(\Xh,\Dh;\Q)$.
Recall that the dual weight and Hodge filtrations of a \Q-MHS
$(M_{\Q},W_{\Bul},F\bul)$ is defined as
(cf. \cite[Ch.~6, 1.5, p.~143]{gelfand_manin}, 
\cite[1.1.6, p.~7]{hodge2})
\begin{align*}
& W_p M_{\Q}^{\vee}:=
  \left(M_{\Q} / W_{-1-p}M_{\Q}\right)^{\vee},\quad\text{and} \\[1mm]
& F^p M_{\C}^{\vee}:=
  \left(M_{\C} / F^{1-p}M_{\C}\right)^{\vee},\
\quad\text{for}\quad p\in\Z \, .
\end{align*}

Let us summarize our results so far.

\begin{proposition}
\label{exp5prop1}

Let $X$, $D$ be as in Subsection \ref{confsubsec} and let $\Xh$, $\Dh$
be the associated complex analytic spaces (cf. Subsection \ref{defcomplexanspace}).
Then the homology group $\h_2^{\rm sing}(\Xh,\Dh;\Q)$ of the pair $(X,D)$
carries a \Q-MHS $(W_{\Bul},F\bul)$.
If $\gamma_0$, $\gamma_1$, $\gamma_2$, $\gamma_3$ is the basis of
$\h_2^{\rm sing}(\Xh,\Dh;\Q)$ defined in Subsection \ref{exp5singsubsec} and 
$\omega_0^{\ast}$, $\omega_1^{\ast}$, $\omega_2^{\ast}$,
$\omega_3^{\ast}$ the basis of $\h_2^{\rm sing}(\Xh,\Dh;\C)$ dual to
the one considered in Subsection \ref{exp5drsubsec},
we find
\begin{enumroman}
\item
a (transposed) period matrix $P$ expressing the $\gamma_j$ in terms of the
$\omega_i^{\ast}$

$$
(\gamma_0,\gamma_1,\gamma_2,\gamma_3)=
\left(
\begin{matrix}
\omega_0^{\ast} \\
\omega_1^{\ast} \\
\omega_2^{\ast} \\
\omega_3^{\ast}
\end{matrix}
\right)
\cdot
\left(
\begin{matrix}
1 & 0 & 0 & 0 \\
\Li_1(\frac{1}{b}) & 2\pi i & 0 & 0 \\
\Li_1(\frac{1}{a}) & 0 & 2\pi i & 0 \\
\Li_{1,1}\!\left(\frac{b}{a},\frac{1}{b}\right) 
  & 2\pi i\Li_1(\frac{b}{a}) 
  & 2\pi i\ln\left(\frac{a-b}{1-b}\right) 
  & (2\pi i)^2
\end{matrix}
\right),
$$
\item
the weight filtration in terms of the $\{\gamma_j\}$
$$
W_p\h_2^{\rm sing}(\Xh,\Dh;\Q)=
\begin{cases}
0 & \text{for}\quad p\le -5 \\
\Q\gamma_3 &\text{for}\quad p=-4,-3 \\
\Q\gamma_1\oplus\Q\gamma_2\oplus\Q\gamma_3 &\text{for}\quad p=-2,-1 \\
\Q\gamma_0\oplus\Q\gamma_1\oplus\Q\gamma_2\oplus\Q\gamma_3 
  &\text{for}\quad p\ge 0,
\end{cases}
$$
\item
the Hodge filtration in terms of the $\{\omega_i^{\ast}\}$
$$
F^p\h_2^{\rm sing}(\Xh,\Dh;\Q)=
\begin{cases}
\C\omega_0^{\ast}\oplus\C\omega_1^{\ast}\oplus
\C\omega_2^{\ast}\oplus\C\omega_3^{\ast}
&\text{for}\quad p\le -2 \\
\C\omega_0^{\ast}\oplus\C\omega_1^{\ast}\oplus
\C\omega_2^{\ast}
&\text{for}\quad p=-1 \\
\C\omega_0^{\ast}
&\text{for}\quad p=0 \\
0 
&\text{for}\quad p\ge 1.
\end{cases}
$$
\end{enumroman}
\end{proposition}

\begin{remark}
This \Q-MHS resembles very much the \Q-MHS considered in \cite[2.2,
p.~620]{goncharov:dlog} and \cite[3.2, p.~6]{zhao03a}. Nevertheless
a few differences are note-worthy:
\begin{itemize}
\item
Goncharov defines the weight filtration slightly different, for
example his lowest weight is $-6$.
\item
The entry $P_{3,2}=2\pi i\ln\!\left(\frac{a-b}{1-b}\right)$ of the period
matrix $P$ differs by $(2\pi i)^2$, or put differently, the basis
$\{\gamma_0,\gamma_1,\gamma_2-\gamma_3,\gamma_3\}$
is used.
\end{itemize}
\end{remark}

\subsubsection{A Variation of Mixed \Q-Hodge Structures}

Up to now, the parameters $a$ and $b$ of the configuration $(X,D)$
have been fixed.
By varying $a$ and $b$, we obtain a family of configurations:
Equip $\A^2_{\C}$ with coordinates $a$ and $b$ and let
$$
B:=\A^2_{\C}\setminus\left(\{a=0\}\cup\{a=1\}\cup\{b=0\}\cup\{b=1\}\right)
$$
be the parameter space. Take another copy of $\A^2_{\C}$ with
coordinates $x$ and $y$ and define total spaces
\begin{align*}
& \ul{X}:=(\underset{(a,b,x,y)}{B\times \A^2_{\C}})\setminus
\left(\{x = a\}\cup\{y = b\}\right),\quad\text{and} \\
& \ul{D}:=\text{``}B\times D\text{''}=
\ul{X}\cap\left(\{x=0\}\cup\{y=1\}\cup\{x=y\}\right).
\end{align*}
We now have a projection
$$
\begin{matrix}
\ul{D} & \inclusion & \ul{X} & \hspace{5mm} & (a,b,x,y) \\
& \searrow & \down[\pi] & & \down[] \\
& & B & & (a,b)
\end{matrix},
$$
whose fibre over a closed point $(a,b)\in B$ is precisely the
configuration
$(X,D)$ for the parameter choice $a$, $b$.

{\em $\pi$ is flat.}\/ Since open immersions and structure morphisms
to the spectrum of a field are flat and flatness is stable under base
extension and composition \cite[Prop.~9.2, p.~254]{hartshorne}, we see
that $\pi$ and $\pi_{|\ul{D}}$ define flat families.
Note that these families could be defined over $\Z$ or $\Q$ already.
This justifies the use of the term period on the preceding pages: For
$(a,b)\in\Bh$ with $a,b\in\Q$, the fibre over $(a,b)$ of $\pi$ is the
base change to \C\ of a pair of varieties $(X_0,D_0)$ defined over \Q{}
and the ``differential forms'' $\omega_0$, $\omega_1$, $\omega_2$
$\omega_3$ span $\hDR{2}(X_0,D_0/\Q)$. Hence the complex numbers
$P_{ij}$ computed above are indeed periods in the sense of definitions
\ref{per1}, \ref{per2}, \ref{per3}, if $a,b\in\Q$.

In the following, we identify $\h_2^{\rm sing}(\Xh,\Dh;\C)$ with
$\C^4$ by mapping the $\omega_i^{\ast}$ to a standard basis of
$\C^4$. 
For each choice of parameters $(a,b)\in\Bh$, this defines a \Q-lattice
of full rank in $\C^4$ together with weight and Hodge filtrations that
turn $\C^4$ into a \Q-MHS. Note that the \Q-lattice varies when we move $(a,b)$
inside $\Bh$. Assigning to each $(a,b)\in\Bh$ this \Q-MHS defines not
merely a ``family of \Q-MHS'', but a 
{\em good unipotent variation of mixed \Q-Hodge structures}. 
\index{Variation of mixed \Q-Hodge structures}
Again we have to refer to the literature for the definition of this concept.
See for example 
\cite{hain_zucker},
\cite[p.~37f]{steenbrink} or
\cite[3.1, p.~4f]{zhao03b}.
For explicitness, we pick a branch for each entry $P_{ij}$ of the
(transposed) period matrix $P$.

\begin{proposition}[{\cite[Thm.~4.1, p.~184]{zhao02}}]
\label{VMHS}
The assignment
$$
\Bh\ni(a,b) \mapsto (V_{\Q},W_{\Bul},F\bul),
$$
where
$$
V_{\Q}:=\linspan_{\Q}\{s_0,\ldots,s_3\},
$$
$$
V_{\C}:=\C^4\quad\text{with standard basis $e_0,\ldots,e_3$},
$$
$$
s_0:=\left(\begin{matrix}
1 \\ \Li\!\left(\frac{1}{b}\right) \\ \Li_1\!\left(\frac{1}{a}\right) \\
\Li_{1,1}\!\left(\frac{b}{a},\frac{1}{b}\right)
\end{matrix}\right),\,
s_1:=\left(\begin{matrix}
0 \\ 2\pi i \\ 0 \\ 2\pi i\Li_1\!\left(\frac{b}{a}\right)
\end{matrix}\right),\,
s_2:=\left(\begin{matrix}
0 \\ 0 \\ 2\pi i \\ 2\pi i\ln\!\left(\frac{a-b}{1-b}\right)
\end{matrix}\right),\,
s_3:=\left(\begin{matrix}
0 \\ 0 \\ 0 \\ (2\pi i)^2 
\end{matrix}\right),
$$
$$
W_p V_{\Q}=\begin{cases}
0 & \text{for}\quad p\le -5\\
\Q s_3 & \text{for}\quad p=-4,-3\\
\Q s_1\oplus \Q s_2\oplus\Q s_3 &\text{for}\quad p=-2,-1\\
V_{\Q} & \text{for}\quad p\ge 0,\quad\text{and}
\end{cases}
$$
$$
F^p V_{\C}=\begin{cases}
V_{\C} &\text{for}\quad p\le -2\\
\C e_0\oplus\C e_1\oplus \C e_2 &\text{for}\quad p=-1\\
\C e_0 &\text{for}\quad p=0\\
0 &\text{for}\quad p\ge 1
\end{cases}
$$
defines a good unipotent variation of \Q-MHS on $\Bh$.
\end{proposition}

Note that the Hodge filtration $F\bul$ does not depend on
$(a,b)\in\Bh$.

One of the main characteristics of good unipotent variations of \Q-MHS
is that they can be extended to a compactification of the base space (if the
complement is a divisor with normal crossings).

The algorithm for computing these extensions, so called {\em limit
mixed \Q-Hodge structures}, can be found for example in 
\cite[7, p.~24f]{hain} and 
\cite[4, p.~12]{zhao03b}.

\subsubsection{Limit Mixed \Q-Hodge Structures}

In a first step, we extend the variation from Proposition \ref{VMHS} to
the divisor $\{a=1\}$ minus the point $(1,0)$ and then in a second
step we extend it to the point $(1,0)$. We will use the notation of
Proposition \ref{VMHS}. In particular, we assume that a branch has
been picked for each entry $P_{ij}$ of $P$. We will follow \cite[4.1,
p.~14f]{zhao03b} very closely.

{\em First step:}\/ Let $\sigma$ be the loop from Lemma \ref{dlogmono} winding
counterclockwise around $\{a=1\}$ once, but not around $\{a=0\}$,
$\{b=0\}$ or $\{b=1\}$. If we analytically continue an entry $P_{ij}$
of $P$ along $\sigma$ we possibly get a second branch of the same
multivalued function. In fact, by Lemma \ref{dlogmono}, the matrix resulting
from analytic continuation of every entry along $\sigma$ will be of
the form
$$
P\cdot T_{\{a=1\}},
$$
where
$$
T_{\{a=1\}}=
\left(
\begin{matrix}
1 & 0 & 0 & 0 \\
0 & 1 & 0 & 0 \\
-1 & 0 & 1 & 0 \\
0 & 0 & 0 & 1 
\end{matrix}
\right)
$$
is the {\em monodromy matrix} corresponding to $\sigma$. The {\em local
monodromy logarithm} is defined as
\begin{align*}
N_{\{a=1\}} & =\frac{\ln T_{\{a=1\}}}{2\pi i} =
\frac{1}{2\pi i}\sum_{n=1}^{\infty} \frac{-1}{n} \left( \left(
\begin{smallmatrix} 1 &  &  &  \\  & 1 &  &  \\  &  & 1 &  \\
   &  &  & 1\end{smallmatrix}\right)-T_{\{a=1\}}\right)^n\\
& =
\left(
\begin{matrix}
0 & 0 & 0 & 0 \\
0 & 0 & 0 & 0 \\
\frac{-1}{2\pi i} & 0 & 0 & 0 \\
0 & 0 & 0 & 0 
\end{matrix}
\right).
\end{align*}
We want to extend our \Q-MHS along the tangent vector
$\frac{\del}{\del a}$,
i.e. we introduce a local coordinate $t:=a-1$ and compute the {\em
limit period matrix}
\begin{align*}
P_{\{a=1\}}
& :=\lim_{t\to 0}P\cdot e^{-\ln(t)\cdot N_{\{a=1\}}} \\
&=\lim_{t\to 0}
\left(
\begin{matrix}
1 \!&\! 0 \!&\! 0 \!&\! 0 \\
\Li_1\!\left(\frac{1}{b}\right) \!&\! 2\pi i \!&\! 0 \!&\! 0 \\
\Li_1\!\left(\frac{1}{1+t}\right) \!&\! 0 \!&\! 2\pi i \!&\! 0 \\
\Li_{1,1}\!\left(\frac{b}{1+t},\frac{1}{b}\right)
\!&\! 2\pi i\Li_1\!\left(\frac{b}{1+t}\right) 
\!&\! 2\pi i\ln\!\left(\frac{1-b+t}{1-b}\right)
\!&\! (2\pi i)^2 
\end{matrix}
\right)
\cdot 
\left(
\begin{matrix}
1 \!&\! 0 \!&\! 0 \!&\! 0 \\
0 \!&\! 1 \!&\! 0 \!&\! 0 \\
\frac{ln(t)}{2\pi i} \!&\! 0 \!&\! 1 \!&\! 0 \\
0 \!&\! 0 \!&\! 0 \!&\! 1
\end{matrix}
\right) \\
\!&\!=\lim_{t\to 0}
\left(
\begin{matrix}
1 \!&\! 0 \!&\! 0 \!&\! 0 \\
\Li_1\!\left(\frac{1}{b}\right) \!&\! 2\pi i \!&\! 0 \!&\! 0 \\
\Li_1\!\left(\frac{1}{1+t}\right)+\ln(t) \!&\! 0 \!&\! 2\pi i \!&\! 0 \\
\Li_{1,1}\!\left(\frac{b}{1+t},\frac{1}{b}\right)+
\ln\!\left(\frac{1-b+t}{1-b}\right)\cdot\ln(t) 
\!&\! 2\pi i\Li_1\!\left(\frac{b}{1+t}\right) 
\!&\! 2\pi i\ln\!\left(\frac{1-b+t}{1-b}\right)
\!&\! (2\pi i)^2 
\end{matrix}
\right) \\
\!&\!\stackrel{(\ast)}{=}
\left(
\begin{matrix}
1 \!&\! 0 \!&\! 0 \!&\! 0 \\
\Li_1\left(\frac{1}{b}\right) \!&\! 2\pi i \!&\! 0 \!&\! 0 \\
0 \!&\! 0 \!&\! 2\pi i \!&\! 0 \\
-\Li_2\!\left(\frac{1}{1-b}\right) 
\!&\! 2\pi i\Li_1(b) \!&\! 0 \!&\! (2\pi i)^2 
\end{matrix}
\right).
\end{align*}
Here we used at $(\ast)$
\begin{itemize}
\item
${P_{\{a=1\}}}_{2,0}=
\lim_{t\to 0}\Li_1\!\left(\frac{1}{1+t}\right)+\ln(t) \\[1mm]
\phantom{{P_{\{a=1\}}}_{2,0}}=
\lim_{t\to 0}-\ln\!\left(1-\frac{1}{1+t}\right)+\ln(t) \\[1mm]
\phantom{{P_{\{a=1\}}}_{2,0}}=
\lim_{t\to 0}-\ln(t)+\ln(1+t)+\ln(t) \\
\phantom{{P_{\{a=1\}}}_{2,0}}=0,\quad\text{and}$
\item
${P_{\{a=1\}}}_{3,0}=
\lim_{t\to 0}\Li_{1,1}\!\left(\frac{b}{1+t},\frac{1}{b}\right)+
\ln\!\left(\frac{1-b+t}{1-b}\right)\cdot\ln(t) \\[1mm]
\phantom{{P_{\{a=1\}}}_{3.0}}=
\Li_{1,1}\!\left(b,\frac{1}{b}\right)\quad\displaystyle\text{by L'Hospital} \\[1mm]
\phantom{{P_{\{a=1\}}}_{3.0}}=
-\Li_2\!\left(\frac{1}{1-b}\right)\quad\displaystyle\text{by Lemma \ref{dlogdilog}}$.
\end{itemize}
The vectors $s_0$, $s_1$, $s_2$, $s_3$ spanning the \Q-lattice of the
limit \Q-MHS on $\{a=1\}\setminus\{(1,0)\}$ are now given by the
columns of the limit period matrix
$$
s_0=
\left(
\begin{matrix}
1 \\
\Li_1\!\left(\frac{1}{b}\right) \\
0 \\
-\Li_2\!\left(\frac{1}{1-b}\right)
\end{matrix}
\right),\quad
s_1=
\left(
\begin{matrix}
0 \\ 2\pi i \\ 0 \\ 2\pi i \Li_1(b)
\end{matrix}
\right),\quad
s_2=
\left(
\begin{matrix}
0 \\ 0 \\ 2\pi i \\ 0 
\end{matrix}
\right),\quad
s_3=
\left(
\begin{matrix}
0 \\ 0 \\ 0 \\ (2\pi i)^2
\end{matrix}
\right).
$$
The weight and Hodge filtration of the limit \Q-MHS can be
expressed in terms of the $s_j$ and the standard basis vectors $e_i$
of $\C^4$ as in Proposition \ref{VMHS} using the same formulae as given
there.
This gives us a variation of \Q-MHS on the divisor
$\{a=1\}\setminus\{(1,0)\}$.
This variation is actually (up to signs) an extension of
Deligne's famous {\em dilogarithm variation} considered for example in
\cite[4.2, p.~38f]{kleinjung}. In loc. cit. the geometric origin of
this variation is explained in detail.

{\em Second step:}\/ We now extend this variation along the tangent
vector $\frac{-\del}{\del b}$ to the point $(1,0)$, i.e. we write
$b=-t$ with a local coordinate $t$.
Let $\sigma$ be the loop in $\{a=1\}\setminus\{(1,0)\}$ from Lemma
\ref{dilogmono} winding counterclockwise around $(1,0)$ once, but not around
$(1,1)$. Then by Lemma \ref{dilogmono} the monodromy matrix
corresponding to $\sigma$  is given by 
$$
T_{(1,0)}=
\left(
\begin{matrix}
1 & 0 & 0 & 0 \\
1 & 1 & 0 & 0 \\
0 & 0 & 1 & 0 \\
0 & 0 & 0 & 1
\end{matrix}
\right),
$$
hence the local monodromy logarithm is given by
$$
N_{(1,0)}=\frac{\ln T_{(1,0)}}{2\pi i}
=\left(
\begin{matrix}
0 & 0 & 0 & 0 \\
\frac{1}{2\pi i} & 0 & 0 & 0 \\
0 & 0 & 0 & 0 \\
0 & 0 & 0 & 0 
\end{matrix}
\right).
$$
Thus we get for the limit period matrix
\begin{align*}
P_{(1,0)}:=& 
\lim_{t\to 0}P_{\{a=1\}}\cdot e^{-\ln(t)\cdot N_{(1,0)}} \\
=&\lim_{t\to 0}
\left(
\begin{matrix}
1 & 0 & 0 & 0 \\
\Li_1\!\left(\frac{-1}{t}\right) & 2\pi i & 0 & 0 \\
0 & 0 & 2\pi i & 0 \\
-\Li_2\!\left(\frac{1}{1+t}\right) & 2\pi i\,\Li_1(-t) & 0 & (2\pi i)^2
\end{matrix}
\right)
\cdot
\left(
\begin{matrix}
1 & 0 & 0 & 0 \\
\frac{-\ln(t)}{2\pi i} & 1 & 0 & 0 \\
0 & 0 & 1 & 0 \\
0 & 0 & 0 & 1
\end{matrix}
\right) \\
=&\lim_{t\to 0}
\left(
\begin{matrix}
1 & 0 & 0 & 0 \\
\Li_1\!\left(\frac{-1}{t}\right)-\ln(t) & 2\pi i & 0 & 0 \\
0 & 0 & 2\pi i & 0 \\
-\Li_2\!\left(\frac{1}{1+t}\right) -
\Li_1(-t)\cdot\ln(t) & 0 &
0 & (2\pi i)^2
\end{matrix}
\right) \\
\stackrel{(\ast)}{=}&
\left(
\begin{matrix}
1 & 0 & 0 & 0 \\
0 & 2\pi i & 0 & 0 \\
0 & 0 & 2\pi i & 0 \\
-\zeta(2) & 0 & 0 & (2\pi i)^2
\end{matrix}
\right).
\end{align*}
Here we used at $(\ast)$
\begin{itemize}
\item
${P_{(1,0)}}_{1,0}=
\lim_{t\to 0}\Li_1\!\left(\frac{-1}{t}\right)-\ln(t) \\
\phantom{{P_{(1,0)}}_{1,0}}=
\lim_{t\to 0}-\ln\!\left(1+\frac{1}{t}\right)-\ln(t) \\
\phantom{{P_{(1,0)}}_{1,0}}=
\lim_{t\to 0}-\ln(1+t)+\ln(t)-\ln(t) \\
\phantom{{P_{(1,0)}}_{1,0}}=0,\quad\text{and}$
\item
${P_{(1,0)}}_{3,0}=
\lim_{t\to 0}-\Li_2\!\left(\frac{1}{1+t}\right)-\Li_1(-t)\cdot\ln(t) \\[1mm]
\phantom{{P_{(1,0)}}_{3,0}}=
\lim_{t\to 0}\Li_2\!\left(\frac{1}{1+t}\right)+\ln(1+t)\cdot\ln(t) \\[1mm]
\phantom{{P_{(1,0)}}_{3,0}}
=-\Li_2(1)\quad\text{by L'Hospital} \\
\phantom{{P_{(1,0)}}_{3,0}}=
-\zeta(2).$
\end{itemize}
As in the previous step, the vectors $s_0$, $s_1$, $s_2$, $s_3$
spanning the \Q-lattice of the limit \Q-MHS are given by the columns
of the limit period matrix $P_{(1,0)}$ and weight and Hodge
filtrations by the formulae in Proposition \ref{VMHS}.

So we obtained $-\zeta(2)$ as a ``period'' of a limit \Q-MHS.

\appendix
\section{Sign Conventions and Related Material}
\label{sign}

We adopt some sign conventions from \cite[A~3.12,
p. 657f]{eisenbud}. Let \AAA\ be an abelian category.
\begin{itemize}
\item
Any \AAA-object $\a$ gives rise to an \AAA-complex $\a[0]$
concentrated in degree zero with $(\a[0])^0=\a$.
\notation{0}{$\a[0]$}{complex concentrated in degree zero with
$(\a[0])^0=\a$} Similarly, we define $\a[-1]$.
\item
For \AAA-complexes $\a\bul$ with differential $d_{\a\bul}$, we have an
$i$-th {\em cohomology functor} 
$\hh^i\a\bul:=\ker d_{\a\bul}^i / \im d_{\a\bul}^{i-1}$.
\notation{H}{$\hh^i\a\bul$}{$i$-th cohomology of a complex $\a\bul$}
\item
If $\a\bul$ is a complex of \AAA-objects with differential
$d_{\a\bul}$, we write $\a\bul[j]$ for the {\em complex shifted by $j$}%
\notation{0}{$\a\bul[j]$}{complex shifted by $j$}%
(to the left) with $\a^i[j]:=\a^{i+j}$ and 
$d_{\a\bul[j]} := (-1)^j d_{\a\bul}$.
\item
If $f:\b\bul\to\a\bul$ is a morphism of \AAA-complexes, we define
the {\em mapping cone}
\index{Mapping cone}
$M_f$ of $f$ by
$$
M_f:=\a\bul[-1]\oplus\b\bul
$$
with differential
$$
d_{M_f}:=\left(\begin{matrix} d_{\a\bul[-1]} & f \\ 
0 & d_{\b\bul} \end{matrix}\right),
$$
i.e. 
\begin{align*}
d^{\,i}_{M_f} : \a^{i-1} \oplus \b^i & \to \a^i \oplus \b^{i+1} \\
(a,b) \quad & \mapsto (-d_{\a\bul}^{\,i-1}a+f(b),d_{\b\bul}^{\,i} b).
\end{align*}
\item
By a {\em double complex} 
\index{Double complex}
$\a^{\Bul,\Bul}$ of \AAA-objects, we mean a
doubly indexed family of \AAA-objects $\{\a^{p,q}\}_{p,q\in\Z}$
together with vertical differentials
$$
d^{p,q}_{\rm I}: \a^{p,q} \to \a^{p+1,q} 
$$
and horizontal differentials
$$
d^{p,q}_{\rm II}: \a^{p,q} \to \a^{p,q+1},
$$
which satisfy for all $p,q\in\Z$
\begin{align*}
& d^{p+1,q}_{\rm I} \circ d^{p,q}_{\rm I} = 0, \\
& d^{p,q+1}_{\rm II} \circ d^{p,q}_{\rm II} = 0\quad\text{and}\\
& d^{p,q+1}_{\rm I} \circ d^{p,q}_{\rm II} = 
d^{p+1,q}_{\rm II} \circ d^{p,q}_{\rm I}.
\end{align*}
\item
If $\a^{\Bul,\Bul}$ is a double complex, we define its {\em total complex}
\index{Total complex}
$$
\a\bul := \tot \a^{\Bul,\Bul}
$$
  \notation{T}{$\tot \a^{\Bul,\Bul}$}{total complex 
  of a double complex $\a^{\Bul,\Bul}$}
by $\a^n:=\bigoplus_{p+q=n} \a^{p,q}$ with differential
$$
d_{\a\bul}^n=
\bigoplus_{p+q=n}\left(d^{p,q}_{\rm I} + 
  (-1)^p d^{p,q}_{\rm II}\right):
\a^n \to \a^{n+1}.
$$
\end{itemize}

\bibliographystyle{arbeit}
\bibliography{arbeit}

\begin{thebibliography}{Se:gaga}
\providecommand{\bysame}{\leavevmode\hbox
  to4em{\hrulefill}\thinspace:\newblock}
\expandafter\ifx\csname url\endcsname\relax
  \def\url#1{\texttt{#1}}\fi
\expandafter\ifx\csname urlprefix\endcsname\relax\def\urlprefix{URL }\fi

\bibitem[AH]{atiyah_hodge}
M.~Atiyah and W.~Hodge: \emph{Integrals of the second kind on an algebraic
  variety}.
\newblock Annals Math. \textbf{62} (1955), pp.~56--91.

\bibitem[Br]{bredon}
G.~Bredon: \emph{Topology and {G}eometry}.
\newblock Springer, 1993.

\bibitem[D2]{hodge2}
P.~Deligne: \emph{Th{\'e}orie des {H}odge {I}{I}}.
\newblock Publication Math{\'e}matique IHES \textbf{40} (1971), pp.~5--57.

\bibitem[D3]{hodge3}
\bysame \emph{Th{\'e}orie des {H}odge {I}{I}{I}}.
\newblock Publication Math{\'e}matique IHES \textbf{44} (1974), pp.~5--77.

\bibitem[Dur]{durfee}
A.~Durfee: \emph{A naive guide to mixed {H}odge theory}.
\newblock Proc. Symp. Pure Math. \textbf{40} (1983), no.~1, pp.~313--320.

\bibitem[Eb]{eisenbud}
D.~Eisenbud: \emph{Commutative {A}lgebra with a {V}iew {T}oward {A}lgebraic
  {G}eometry}.
\newblock Springer, 1994.

\bibitem[Fi]{fichtenholz}
G.~Fichtenholz: \emph{Differential- und {I}ntegralrechnung}.
\newblock VEB Deutscher Verlag der Wissenschaften, 1990.

\bibitem[GM]{gelfand_manin}
S.~Gelfand and Y.~Manin: \emph{Homological {A}lgebra}.
\newblock Springer, 1999.

\bibitem[G1]{goncharov:dlog}
A.~B. Goncharov: \emph{The double logarithm and {M}anin's complex for modular
  curves}.
\newblock Math. Res. Let. \textbf{4} (1997), pp.~617--636.

\bibitem[G2]{goncharov:dihedral}
\bysame \emph{The dihedral lie algebra and {G}alois symmetries of
  $\pi^{l}_1({\P}^1-(\{0,\infty\}\cup\mu_n))$}.
\newblock Duke Math. J. \textbf{110} (2001), no.~3, pp.~397--487.

\bibitem[G3]{goncharov:polylog}
\bysame \emph{Multiple polylogarithms and mixed {T}ate motives}, 2001.
\newblock ArXiv:math.AG/0103059 v4.

\bibitem[GR]{grauert_remmert}
H.~Grauert and R.~Remmert: \emph{Theorie der Steinschen R{\"a}ume}.
\newblock Springer, 1977.

\bibitem[GH]{griffiths_harris}
P.~Griffiths and J.~Harris: \emph{Principles of {A}lgebraic {G}eometry}.
\newblock Wiley, 1978.

\bibitem[Gro]{grothendieck_66}
A.~Grothendieck: \emph{On the de {R}ham cohomology of algebraic varieties}.
\newblock Publication Math{\'e}matique IHES \textbf{29} (1966), pp.~95--103.

\bibitem[Gun]{gunning}
R.~Gunning: \emph{Introduction to {H}olomorphic {F}unctions of {S}everal
  {V}ariables}, vol.~1.
\newblock Wadsworth, 1990.

\bibitem[H]{hain}
R.~M. Hain: \emph{Classical polylogarithms}.
\newblock Proc. Symp. Pure Math. \textbf{55} (1994), no.~2, pp.~3--42.

\bibitem[HZ]{hain_zucker}
R.~M. Hain and S.~Zucker: \emph{A guide to unipotent variations of mixed
  {H}odge structures}.
\newblock LN 1246. Springer, 1987.

\bibitem[Ha:dR]{hartshorne_71}
R.~Hartshorne: \emph{On the de {R}ham cohomology of algebraic varieties}.
\newblock Publication Math{\'e}matique IHES \textbf{45} (1971), pp.~1--99.

\bibitem[Ha]{hartshorne}
\bysame \emph{Algebraic {G}eometry}.
\newblock Springer, 1977.

\bibitem[Hi1]{hironaka_64}
H.~Hironaka: \emph{Resolution of singularities}.
\newblock Annals Math. \textbf{79} (1964), pp.~109--326.

\bibitem[Hi2]{hironaka_74}
\bysame \emph{Triangulation of algebraic sets}.
\newblock In: \emph{Algebraic Geometry, Arcata 1974}, vol.~29 of \emph{Proc.
  Symp. Pure Math.}, pp. 165--185. AMS, 1975.

\bibitem[HW]{huber_wildeshaus}
A.~Huber and J.~Wildeshaus: \emph{The classical polylogarithm, abstract of a
  series of lectures given at the workshop on polylogs in {E}ssen, 1997}.
\newblock Available online at {\tt http://
  www.mathematik.uni-leipzig.de/\~{}huber/preprints/essen.dvi}.

\bibitem[Hui]{huisman}
J.~Huisman: \emph{Real algebraic differential forms on complex algebraic
  varieties}.
\newblock Indag. Math. (N.S.) \textbf{11} (2000), no.~1, pp.~63--71.

\bibitem[Kj]{kleinjung}
T.~Kleinjung: \emph{{G}erahmte gemischte {T}ate-{M}otive und die {W}erte von
  {Z}etafunktionen zu {Z}ahlk{\"o}rpern an den {S}tellen $2$ und $3$}.
\newblock Ph.D. thesis, Rheinische Friedrich-Wilhelms-Universit{\"a}t Bonn,
  Bonner Mathematische Schriften, Nr. {\bf 340}, ISSN 0524-045X, 2000.

\bibitem[K]{kontsevich}
M.~Kontsevich: \emph{Operads and motives in deformation quantisation}.
\newblock Letters in Mathematical Physics \textbf{48} (1999), pp.~35--72.

\bibitem[KZ]{kontsevich_zagier}
M.~Kontsevich and D.~Zagier: \emph{Periods}.
\newblock In: B.~Engquis and W.~Schmid, editors, \emph{Mathematics unlimited --
  2001 and beyond}, pp. 771--808. Springer, 2001.

\bibitem[KB]{koopman_brown}
B.~Koopman and A.~Brown: \emph{On the covering of analytic loci by complexes}.
\newblock Trans. Amer. Math. Soc. \textbf{34} (1932), pp.~231--251.

\bibitem[Ma]{matsumura}
H.~Matsumura: \emph{Commutative Algebra}.
\newblock W.A.~Benjamin, 1970.

\bibitem[Sb]{seidenberg}
A.~Seidenberg: \emph{A new decision method for elementary algebra}.
\newblock Annals Math. \textbf{60} (1954), no.~2, pp.~365--374.

\bibitem[Se:gaga]{gaga}
J.-P. Serre: \emph{G{\'e}om{\'e}trie alg{\'e}brique et g{\'e}om{\'e}trie
  analytic}.
\newblock Ann. Inst. Four. \textbf{6} (1956), pp.~1--42.

\bibitem[Se]{serre}
\bysame \emph{A {C}ourse in {A}rithmetic}.
\newblock Springer, 1973.

\bibitem[Sil]{silverman}
J.~Silverman: \emph{The {A}rithmetic of {E}lliptic {C}urves}.
\newblock Springer, 1986.

\bibitem[St]{steenbrink}
J.~Steenbrink: \emph{A summary of mixed {H}odge theory}.
\newblock Proc. Symp. Pure Math. \textbf{55} (1994), no.~1, pp.~31--41.

\bibitem[W]{warner}
F.~Warner: \emph{{F}oundations of {D}ifferentiable {M}anifolds and {L}ie
  {G}roups}.
\newblock Springer, 1983.

\bibitem[Z1]{zhao02}
J.~Zhao: \emph{Multiple polylogarithms: analytic continuation, monodromy, and
  variations of mixed {H}odge structures}.
\newblock In: \emph{Contemporary trends in algebraic geometry and algebraic
  topology (Tianjin, 2000)}, vol.~5 of \emph{Nankai Tracts Math.}, pp.
  167--193. World Sci. Publishing, River Edge, NJ, 2002.

\bibitem[Z2a]{zhao03a}
\bysame \emph{Analytic continuation of multiple polylogarithms}, 2003.
\newblock ArXiv:math.AG/0302054 v2.

\bibitem[Z2b]{zhao03b}
\bysame \emph{Variations of {M}ixed {H}odgestructures of {M}ultiple
  {P}olylogarithms}, 2003.
\newblock ArXiv:math.AG/0302055 v3.

\end{thebibliography}

\printglossary

\printindex

\section*{Erkl\"arung}

\vspace{0.5cm}

Ich versichere, dass ich die vorliegende Arbeit selbstst\"andig und nur
unter Verwendung der angegebenen 
Quellen und Hilfsmittel angefertigt habe.

Diese Arbeit wurde in dieser oder 
\"ahnlicher Form noch keiner anderen
Pr\"u\-fungs\-be\-h\"or\-de vorgelegt.

\vspace{2.5cm}

\begin{tabular}{p{5cm}p{5cm}p{5cm}}
Ort & Datum & Unterschrift 
\end{tabular}

\end{document}